\documentclass[11pt,a4paper]{article}
\usepackage{amsfonts}
\usepackage{latexsym}
\usepackage{cite}
\usepackage{amsmath,amsfonts,latexsym,amssymb}
\usepackage[mathscr]{eucal}
\usepackage{cases,color}
\usepackage{amsthm}

\usepackage[bf,small]{caption2}
\usepackage{float}
\usepackage{graphicx}
\usepackage{amsmath}
\usepackage{amssymb}
\usepackage[all]{xy}

\newtheorem{theorem}{theorem}[section]

\newtheorem{conj}[theorem]{Conjecture}

\newtheorem{cor}[theorem]{Corollary}

\newtheorem{lem}[theorem]{Lemma}
\newtheorem{nota}[theorem]{Notation}

\newtheorem{rmk}[theorem]{Remark}
\newtheorem{thm}[theorem]{Theorem}

\begin{document}

\title{\vspace{-2cm}\textbf{Kauffman bracket skein module of the $(3,3,3,3)$-pretzel link exterior}}
\author{\Large Haimiao Chen}

\date{}
\maketitle

\begin{abstract}
  We show that the Kauffman bracket skein module of the $(3,3,3,3)$-pretzel link exterior over $\mathbb{Q}(q^{\frac{1}{2}})$
  is not finitely generated as a module over $\mathbb{Q}(q^{\frac{1}{2}})[t_1,t_2]$, where $t_1,t_2$ are the meridians of two components.
  This disproves a finiteness conjecture of Detcherry proposed in 2021.

\end{abstract}

\tableofcontents

\section{Introduction}

Let $R$ be a commutative ring with a distinguished invertible element $q^{\frac{1}{2}}$. Given an oriented 3-manifold $M$, its {\it Kauffman bracket skein module} over $R$ is defined as $\mathcal{S}(M;R):=R\mathfrak{F}(M)/{\rm Sk}(M)$, where $R\mathfrak{F}(M)$ is the free $R$-module on the set $\mathfrak{F}(M)$ of isotopy classes of framed links embedded in $M$, and ${\rm Sk}(M)$ is
the submodule of $R\mathfrak{F}(M)$ generated by all instances of the following {\it skein relations}:
\begin{figure}[h]
  \centering
  \includegraphics[width=9cm]{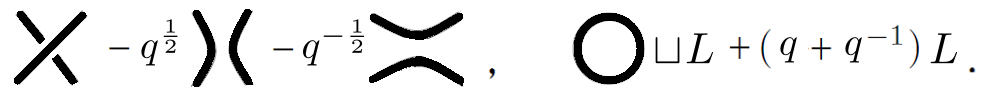}
\end{figure}
\\
As a convention, $R$ is identified with $R\emptyset\subseteq\mathcal{S}(M;R)$, via $1\mapsto\emptyset$.

When $\Sigma$ is an oriented surface, $\mathcal{S}(\Sigma;R):=\mathcal{S}(\Sigma\times[0,1];R)$ is a $R$-algebra, called the {\it Kauffman bracket skein algebra} of $\Sigma$. Given links $L_1,L_2$, their product $L_1L_2$ is defined by stacking $L_1$ over $L_2$ in the $[0,1]$ direction.

When $\partial M\ne\emptyset$, a collar neighborhood of $\partial M$ equips $\mathcal{S}(M;R)$ with the structure of a module over $\mathcal{S}(\partial M;R)$, via stacking.

In the case $R=\mathbb{C}$ and $q^{\frac{1}{2}}=-1$, by \cite{Bu97} there is a ring isomorphism
\begin{align}
\epsilon: \mathcal{S}(M;\mathbb{C})/{\rm nilradical}\stackrel{\cong}\to\mathbb{C}[\mathcal{X}(M)],  \label{eq:classicalization}
\end{align}
where $\mathcal{X}(M)=\mathcal{X}(\pi_1(M))$ is the ${\rm SL}(2,\mathbb{C})$-character variety of $\pi_1(M)$.

In recent years, skein modules have become central objects in quantum topology. One important reason is that, in light of (\ref{eq:classicalization}), skein module is regarded as a {\it quantization} of character variety. More deep studies in this vein were
carried out in \cite{DKS23,DKS25,Ko25}.
Computational results were seen in \cite{AF22,BP22,Ca17,DW21,Gi18,GM19,Le06,LT14,Ma10} and the references therein. However, skein modules still appear rather mysterious; there are many open problems about structures and properties of skein modules, as well as connections to classical invariants.

From now on, let $R=\mathbb{Q}(q^{\frac{1}{2}})$, and abbreviate $\mathcal{S}(M;R)$ to $\mathcal{S}(M)$.

As a remarkable achievement, Gunningham, Jordan and Safronov \cite{GJS23} confirmed Witten's {\it Finiteness Conjecture} which asserts that $\dim_R\mathcal{S}(M)<\infty$ for each closed $3$-manifold $M$.
Some generalizations of the finiteness conjecture to bounded $3$-manifolds were proposed by Detcherry \cite{De21} in 2021.
Suppose $L\subset S^3$ is a link with $n$-components. Let $E_L$ denote its exterior, and let $t_i\in\mathcal{S}(E_L)$ denote the meridian of the $i$-th component.
Then $R[t_1,\ldots,t_n]$ is a subring of $\mathcal{S}(\partial E_L)$.
The following was stated as \cite[Conjecture 3.2]{De21}:
\begin{conj}\label{conj:finiteness}
$\mathcal{S}(E_L)$ is finitely generated as a $R[t_1,\ldots,t_n]$-module.
\end{conj}

\begin{figure}[H]
  \centering
  \includegraphics[width=9cm]{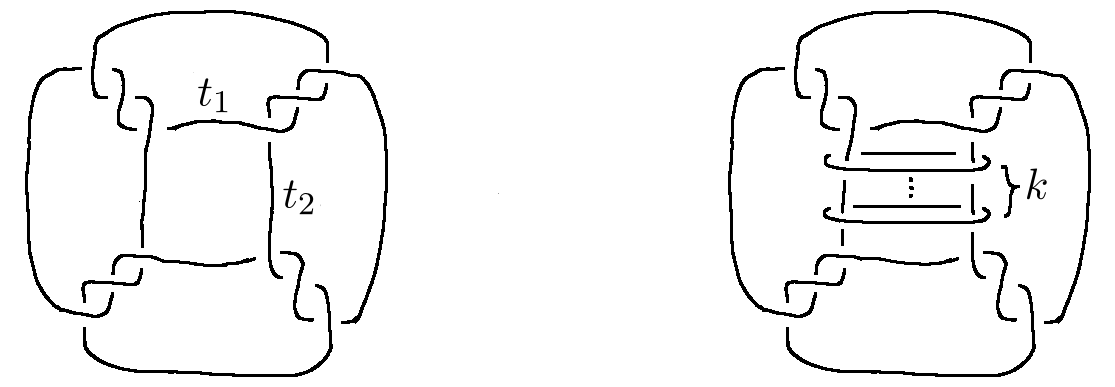}\\
  \caption{Left: the $(3,3,3,3)$-pretzel link $L$; the tangles are numerated by $1,2,3,4$. Right: the element $x_1^k\in \mathcal{S}(E_L)$, made of $k$ parallel copies of a certain knot.}\label{fig:pretzel}
\end{figure}

In this paper, let $L$ denote the $(3,3,3,3)$-pretzel link, as displayed in the left part of Figure \ref{fig:pretzel}. Let $t_1,t_2$ be the meridians of the two components of $L$.
For each $k\ge 1$, let $x_1^k$ denote the element of $\mathcal{S}(E_L)$ given in the right part of Figure \ref{fig:pretzel}.
Based on partial computation, we will show
\begin{thm}\label{thm:main}
For any finite subset $\mathcal{N}\subset\mathcal{S}(E_L)$, there exists $k_0\in\mathbb{N}$ such that the $R[t_1,t_2]$-submodule generated by $\mathcal{N}$ does not contain $x_1^k$ for $k>k_0$. Consequently, $\mathcal{S}(E_L)$ is not a finitely generated $R[t_1,t_2]$-module.
\end{thm}

Thus, $L$ is a counterexample to Conjecture \ref{conj:finiteness}.

The motivation of this work is that, in 2022 the author \cite{Ch22-0} found that for a $4$-strand Montesinos knot $K$, usually $\mathcal{X}(E_K)$ admits a high-dimensional component, i.e., one with dimension larger than $1$. Similarly, $\mathcal{X}(E_L)$ admits a component with dimension larger than $2$.
To be precise, when $\chi(t_1),\chi(t_2)$ are fixed generic values, $\chi(x_1)$ can take an arbitrary generic value. From the viewpoint of (\ref{eq:classicalization}), there would be something strange if Conjecture \ref{conj:finiteness} holds for $L$. So we are curious to seek the relations among the $x_1^k$'s.

To study $\mathcal{S}(E_L)$, we will decompose $E_L$ as
$$E_L=U\cup H_1\cup H_2\cup H_3\cup H_4\cup B^3,$$
where $U$ is a handlebody of genus $4$, the $H_i$'s are $2$-handles, and $B^3$ is a $3$-ball. Basically, $\mathcal{S}(E_L)\cong\mathcal{S}(U)/{\rm Sl}$, where ${\rm Sl}$ is the submodule generated by ``sliding relators". We will find a free basis $\mathcal{B}$ for $\mathcal{S}(U)$, and show that when the relations given by ${\rm Sl}$ are used up, still so many elements of $\mathcal{B}$ survive that $\mathcal{S}(E_L)$ cannot be finitely generated as a $R[t_1,t_2]$-module.

As we will see, ``relations among relations" in $\mathcal{S}(E_L)$ are much richer than those in the classical counterpart $\mathbb{C}[\mathcal{X}(E_L)]$. The non-commutativity and the lack of ring structure of $\mathcal{S}(E_L)$ cause much complexity.
During the process of overcoming various difficulties, we develop some methodologies for studying skein modules.


\medskip

{\bf Notations and conventions}.

Denote $q^{-1}$ by $\overline{q}$, denote $q^{-\frac{1}{2}}$ by $\overline{q}^{\frac{1}{2}}$, and so on. Let $\alpha=q+\overline{q}$.

For a finite set $A$, let $\#A$ denote its cardinality.

Given a subset $\mathcal{X}$ of a left (resp. right) module over a ring $Q$, denote the submodule generated by $\mathcal{X}$ by $Q\mathcal{X}$
(resp. $\mathcal{X}Q$).

Let $\Sigma_{0,n+1}$ denote the $n$-holed disk. Let $\mathcal{S}_n=\mathcal{S}(\Sigma_{0,n+1})$.

For each manifold, its corners are smoothed whenever necessary.

Given manifolds $M_1,M_2$, denote $M_1\cong M_2$ if they are diffeomorphic.

From now on, we use $\mathbf{x},\mathbf{k},\mathbf{p}$, etc, to denote $1$-manifolds,
except for the pretzel link $L$ we work with.

Each framed link is presented as a diagram in the plane of the screen, and the framing vectors are perpendicular to the plane pointing to the reader.

For a 3-manifold $M$ with $\partial M\ne\emptyset$, we always considered $\mathcal{S}(M)$ as a left module over $\mathcal{S}(\partial M)$, unless otherwise specified.

\section{Preparation}

\subsection{Chopping up a $1$-manifold}

\begin{figure}[H]
  \centering
  \includegraphics[width=4cm]{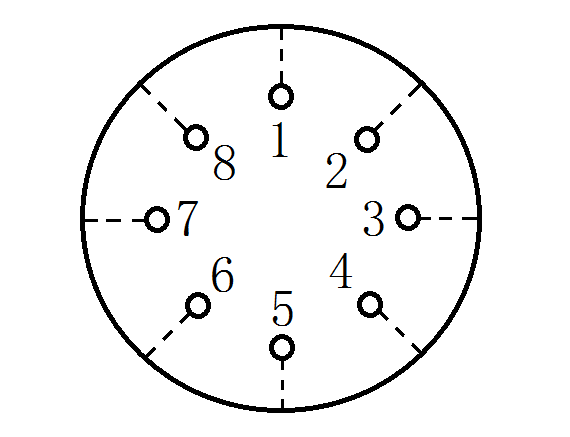}\\
  \caption{$\Sigma_{0,9}$; the dotted lines are $\mathbf{z}_1,\ldots,\mathbf{z}_8$.}\label{fig:disk}
\end{figure}

Let $\Sigma=\Sigma_{0,n+1}$. For $1\le i\le n$, let $\mathbf{z}_i$ be an arc connecting the boundary of $i$-th hole to the outer boundary of $\Sigma$. See Figure \ref{fig:disk} for the case $n=8$.

Fix $\mathsf{a},\mathsf{b}\in\partial(\Sigma\times[0,1])$.
Adopting the notion of \cite{PS00}, by a {\it relative link} we mean a $1$-submanifold $\mathbf{x}=\mathbf{c}\sqcup\mathbf{l}\subset\Sigma\times[0,1]$, where $\mathbf{c}$ is an arc connecting $\mathsf{a}$ to $\mathsf{b}$ and $\mathbf{l}$ is a link in $\Sigma\times(0,1)$. Always assume that $\mathbf{x}$ intersects $\mathbf{z}_i\times[0,1]$ transversally.
Let $\mathcal{S}(\mathsf{a},\mathsf{b})$ denote the $R$-module generated by isotopy classes of relative links modulo skein relations, where isotopies are required to fix $\mathsf{a},\mathsf{b}$.
Let $[\mathbf{x}]\in\mathcal{S}(\mathsf{a},\mathsf{b})$ denote the element represented by $\mathbf{x}$.

Given a relative link $\mathbf{x}$, let $|\mathbf{x}|_i=\#(\mathbf{x}\cap(\mathbf{z}_i\times[0,1]))$, and $|\mathbf{x}|=\sum_{i=1}^n|\mathbf{x}|_i$.

\begin{lem}\label{lem:chop}
Suppose $\mathsf{a},\mathsf{b}\in\Sigma\times\{1\}$. Equip $\mathcal{S}(\mathsf{a},\mathsf{b})$ with the structure of a right $\mathcal{S}_n$-module via stacking. Let $\mathcal{P}\subset\mathcal{S}(\mathsf{a},\mathsf{b})$ denote the subset consisting of elements represented by simple arcs $\mathbf{a}$ with $|\mathbf{a}|\le 2$. Then $\mathcal{S}(\mathsf{a},\mathsf{b})=\mathcal{P}\mathcal{S}_n$.
\end{lem}

\begin{proof}
For a relative link $\mathbf{x}$, we use induction on $|\mathbf{x}|$ to prove $[\mathbf{x}]\in\mathcal{P}\mathcal{S}_n$, which is clear if $|\mathbf{x}|\le 2$.

Suppose $|\mathbf{x}|>2$ and that $[\mathbf{x}']\in\mathcal{P}\mathcal{S}_n$ for any $|\mathbf{x}'|<|\mathbf{x}|$.
Resolve crossings whenever necessary, we may just assume that $\mathbf{x}$ has no crossing.
\begin{enumerate}
  \item Assume $\mathbf{x}$ is a simple arc. Take a subarc $\mathbf{a}\subset\mathbf{x}$ with $|\mathbf{a}|=3$. By \cite[Lemma 3.1]{Ch22},
        $[\mathbf{a}]=\sum_v[\mathbf{d}_v]b_v$ for some $[\mathbf{d}_v]\in\mathcal{P}$, $b_v\in\mathcal{S}_n$. Let $\mathbf{q}_v$
        denote the arc obtained from $\mathbf{x}$ by replacing $\mathbf{a}$ with $\mathbf{d}_v$, then $[\mathbf{x}]=\sum_v[\mathbf{q}_v]b_v$. By the inductive hypothesis, each $[\mathbf{q}_v]\in\mathcal{P}\mathcal{S}_n$. Thus, $[\mathbf{x}]\in\mathcal{P}\mathcal{S}_n$.
  \item If $\mathbf{x}=\mathbf{c}\mathbf{s}$, where $\mathbf{c}$ is an arc, and $\mathbf{s}$ is a simple curve with $|\mathbf{s}|=0$,
        then $[\mathbf{x}]=-\alpha[\mathbf{c}]\in\mathcal{P}\mathcal{S}_n$ by the above paragraph.
  \item If $\mathbf{x}=\mathbf{c}\mathbf{s}$, where $\mathbf{c}$ is an arc and $\mathbf{s}$ is a simple curve with $|\mathbf{s}|>0$,
        then $|\mathbf{c}|<|\mathbf{x}|$. By the inductive hypothesis, $[\mathbf{c}]\in\mathcal{P}\mathcal{S}_n$, hence $[\mathbf{x}]\in\mathcal{P}\mathcal{S}_n$.
\end{enumerate}
\end{proof}

\subsection{On $\mathcal{S}_3$}

\begin{figure}[h]
  \centering
  \includegraphics[width=11.5cm]{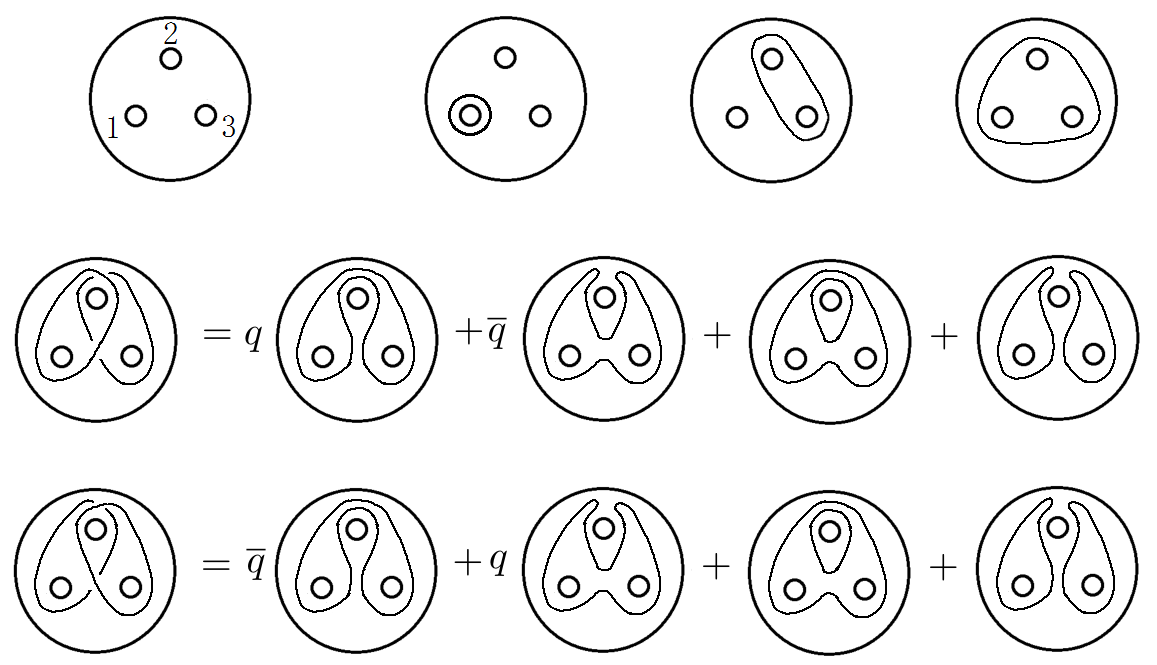}
  \caption{First row (from left to right): $\Sigma_{0,4}$, $t_1$, $t_{23}$, $t_{13}$, $t_{123}$.
  Second row: using skein relations to expand $t_{12}t_{23}$. Third row: expanding $t_{23}t_{12}$.}\label{fig:Sigma1}
\end{figure}

Display $\Sigma_{0,4}$ and introduce some elements of $\mathcal{S}_3$ in the first row of Figure \ref{fig:Sigma1}.

Comparing the second and third rows of Figure \ref{fig:Sigma1}, we see
\begin{align}
t_{12}t_{23}=q^2t_{23}t_{12}+(\overline{q}-q^3)t_{13}+(1-q^2)(t_1t_3+t_2t_{123}).  \label{eq:basic-triple-1}
\end{align}
Rotating by $2\pi/3$ transforms (\ref{eq:basic-triple-1}) into
$$t_{13}t_{12}=q^2t_{12}t_{13}+(\overline{q}-q^3)t_{23}+(1-q^2)(t_2t_3+t_1t_{123}),$$
which is equivalent to
\begin{align}
t_{12}t_{13}=\overline{q}^2t_{13}t_{12}+(q-\overline{q}^3)t_{23}+(1-\overline{q}^2)(t_2t_3+t_1t_{123}).  \label{eq:basic-triple-2}
\end{align}

\subsection{On $\mathcal{S}_4$}

When drawing $\Sigma_{0,5}$, we omit the outer boundary.

Clearly, $R[t_1,t_2,t_3,t_4]$ is a central subalgebra of $\mathcal{S}_4$.
By \cite[Theorem 1.3]{Ch24}, $\mathcal{S}_4$ is a $R[t_1,t_2,t_3,t_4]$-module freely generated by
$$\big\{t_{12}^{j_1}t_{23}^{j_2}t_{34}^{j_3}t_{14}^{j_4}a\colon j_1,\ldots,j_4\ge 0, \ a\in\mathcal{A}\big\},$$
where $\mathcal{A}$ is the union of $\{1,t_0,t_{123},t_{124},t_{134},t_{234}\}$ and
\begin{align*}
\{t_{13}^k, t_{13}^kt_0, t_{13}^kt_{123}, t_{13}^kt_{134}, t_{24}^k, t_{24}^kt_0, t_{24}^kt_{124}, t_{24}^kt_{234}\colon k\ge 1\}.
\end{align*}

\begin{figure}[H]
  \centering
  \includegraphics[width=11cm]{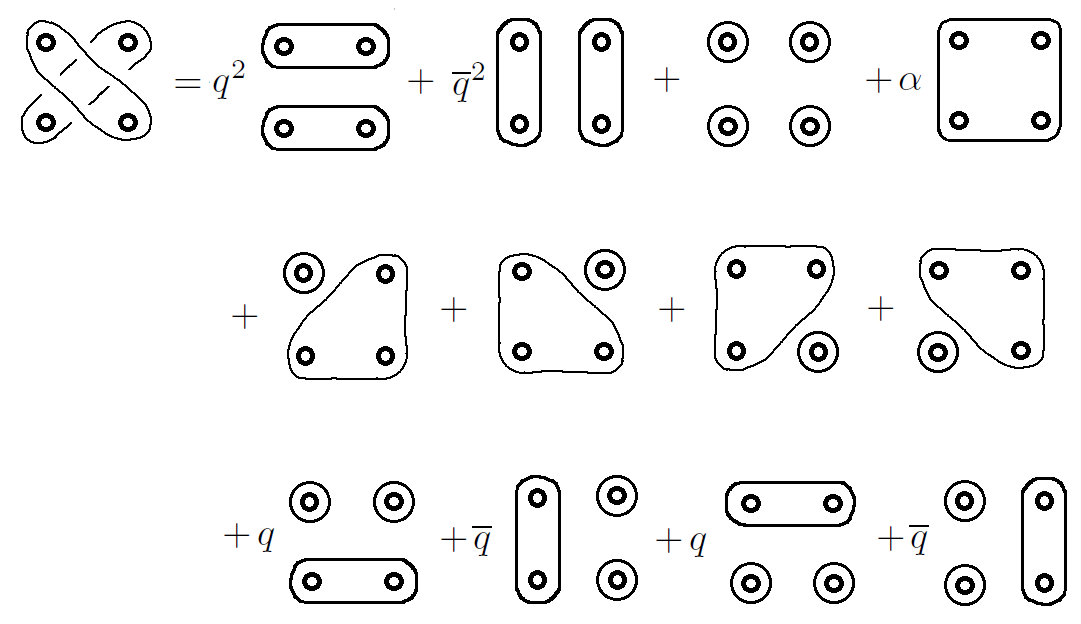}\\
  \caption{Expanding $t_{13}t_{24}$. See \cite{Ch24} for intermediate steps.}\label{fig:13-times-24}
\end{figure}

\begin{figure}[H]
  \centering
  \includegraphics[width=12cm]{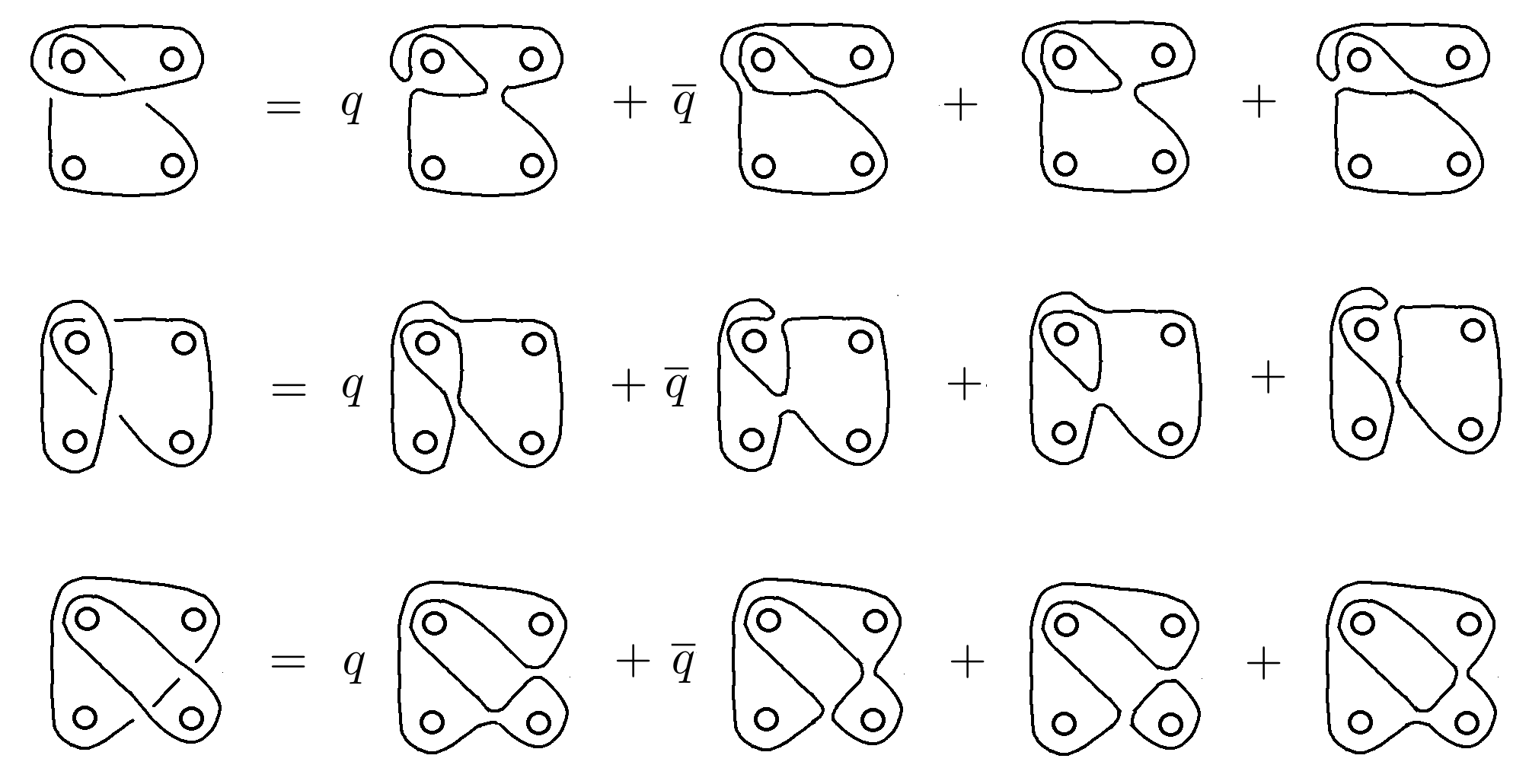}\\
  \caption{Expanding $t_{12}t_{134}$, $t_{14}t_{123}$, $t_{13}t_{124}$.}\label{fig:product}
\end{figure}

\begin{lem}\label{lem:basic-reduce}
The following holds in $\mathcal{S}_4$:
\begin{align*}
\alpha t_{1234}&\equiv t_{13}t_{24}-q^2t_{12}t_{34}-\overline{q}^2t_{14}t_{23}\equiv t_{24}t_{13}-\overline{q}^2t_{12}t_{34}-q^2t_{14}t_{23},  \\
t_{13}t_{124}&\equiv q^2t_{12}t_{134}+\overline{q}^2t_{14}t_{123}\equiv \overline{q}^2t_{134}t_{12}+\overline{q}^2t_{14}t_{123}   \\
&\equiv\overline{q}^2t_{134}t_{12}+q^2t_{123}t_{14}\equiv q^2t_{12}t_{134}+q^2t_{123}t_{14}.
\end{align*}
\end{lem}

\begin{proof}
The first row is seen from Figure \ref{fig:13-times-24}.

As shown in Figure \ref{fig:product},
\begin{align*}
t_{12}t_{134}&=qt_{234}+\overline{q}t_{12\overline{1}34}+t_1t_{1234}+t_2t_{34},  \\
t_{14}t_{123}&=qt_{123\overline{1}4}+\overline{q}t_{234}+t_1t_{1234}+t_4t_{23},  \\
t_{13}t_{124}&=qt_{12\overline{1}34}+\overline{q}t_{123\overline{1}4}+t_1t_{1234}+t_3t_{12\overline{1}4}.
\end{align*}
Eliminating $t_{12\overline{1}34}$ and $t_{123\overline{1}4}$, we obtain
$$t_{13}t_{124}\equiv q^2t_{12}t_{134}+\overline{q}^2t_{14}t_{123}.$$
The other relations can be deduced similarly.
\end{proof}

\subsection{Relations given by handle sliding}

\begin{lem}\label{lem:sliding}
Suppose $M,M'$ are $3$-manifolds such that $M'$ results from attaching $2$-handles $H_1,\ldots,H_r$ to $M$.
Let ${\rm Sl}(M)$ denote the submodule of $R\mathfrak{F}(M)$ generated by $\mathfrak{s}_i(\mathbf{l})-\mathbf{l}$ for links
$\mathbf{l}\hookrightarrow M$ and $i=1,\ldots,r$, where $\mathfrak{s}_i(\mathbf{l})$ denotes the link obtained by sliding $\mathbf{l}$ along the attaching circle of $H_i$.
Then
$$\mathcal{S}(M')\cong\mathcal{S}(M)/{\rm Sl}(M).$$
\end{lem}

This was stated as \cite[Proposition 2.2 (5)]{Pr99}, and also as \cite[Lemma 2.2]{BP22}. For self-containedness, we provide a proof here.

\begin{proof}
For each link $\mathbf{l}\hookrightarrow M'$, employing an isotopy if necessary, we can keep $\mathbf{l}$ away from the $H_i$'s.
Hence the inclusion $M\subset M'$ induces a surjection $\mathfrak{F}(M)\twoheadrightarrow\mathfrak{F}(M')$, and also one
${\rm Sk}(M)\twoheadrightarrow{\rm Sk}(M')$. Choose
a section $\mathfrak{F}(M')\to\mathfrak{F}(M)$ for $\mathfrak{F}(M)\twoheadrightarrow\mathfrak{F}(M')$ and extend it to
$f:R\mathfrak{F}(M')\to R\mathfrak{F}(M)$.
Note that $f(\mathbf{l})-\mathbf{l}\in{\rm Sl}(M)$ for each $\mathbf{l}$.

If $\sum_s\beta_s\mathbf{l}_s=0$ in $\mathcal{S}(M')$, with $\beta_s\in R$ and $\mathbf{l}_s\in\mathfrak{F}(M)$, then due to the surjection
${\rm Sk}(M)\twoheadrightarrow{\rm Sk}(M')$, there exists $\sum_t\gamma_t\mathbf{r}_t\in{\rm Sk}(M)$ (with $\gamma_t\in R$ and $\mathbf{r}_t\in\mathfrak{F}(M)$) given by skein relations such that
$${\sum}_s\beta_s\mathbf{l}_s={\sum}_t\gamma_t\mathbf{r}_t \qquad \text{in\ \ } R\mathfrak{F}(M').$$
It follows that ${\sum}_s\beta_sf(\mathbf{l}_s)={\sum}_t\gamma_tf(\mathbf{r}_t)$ in $R\mathfrak{F}(M).$
Hence
$${\sum}_s\beta_s\mathbf{l}_s={\sum}_t\gamma_t\mathbf{r}_t+{\sum}_t\gamma_t(f(\mathbf{r}_t)-\mathbf{r}_t)
+{\sum}_s\beta_s(\mathbf{l}_s-f(\mathbf{l}_s))\in{\rm Sk}(M)+{\rm Sl}(M).$$
\end{proof}

\begin{rmk}\label{rmk:slide}
\rm For a link $\mathbf{l}$ to be slid along the attaching circle $\mathbf{o}_i$ of $H_i$, employing an isotopy if necessary, we may assume $\mathbf{l}=\mathbf{c}\cup\mathbf{f}\cup\mathbf{l}^0$ and $\mathfrak{s}_i(\mathbf{l})=\mathbf{c}'\cup\mathbf{f}\cup\mathbf{l}^0$, for some link $\mathbf{l}^0$ and some arcs $\mathbf{c},\mathbf{c}',\mathbf{f}$ such that $\mathbf{c}\cup\mathbf{c}'=\mathbf{o}_i$.
In this case, denote $\mathbf{f}$ by ${\rm fa}(\mathbf{l})$, to stand for the ``fixed arc" of $\mathbf{l}$.
\end{rmk}

\section{Generators}

\begin{figure}[h]
  \centering
  \includegraphics[width=12.8cm]{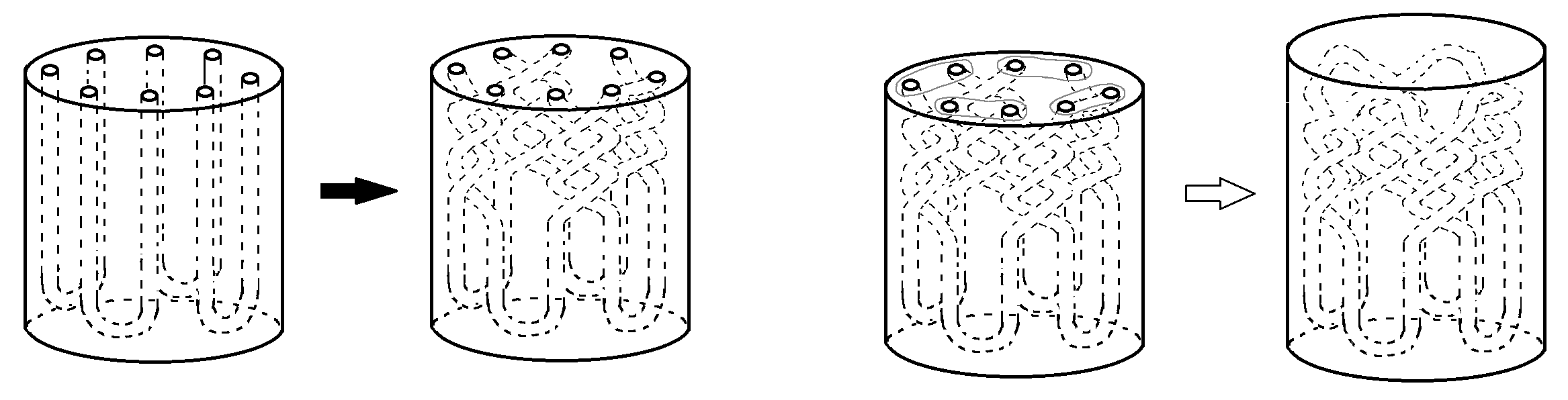}\\
  \caption{Left: twisting the upper half of $U_0$ in a certain way yields $U$.
  Right: attaching four $2$-handles to the top of $U$ along circles (in gray) makes $Y$.}\label{fig:U}
\end{figure}

Let $Y$ denote the exterior of $L$ in a sufficiently large solid cylinder. Clearly, $Y\cong E_L\setminus B^3$ (where $B^3$ is a ball), so the inclusion $Y\hookrightarrow E_L$ induces an isomorphism $\mathcal{S}(Y)\cong\mathcal{S}(E_L)$. A construction of $Y$ is illustrated in Figure \ref{fig:U}, namely, attaching four $2$-handles to $U$ along the dotted circles, where $U$ is the manifold obtained by appropriately twisting the upper half of $U_0$.

\begin{figure}[h]
  \centering
  \includegraphics[width=10.5cm]{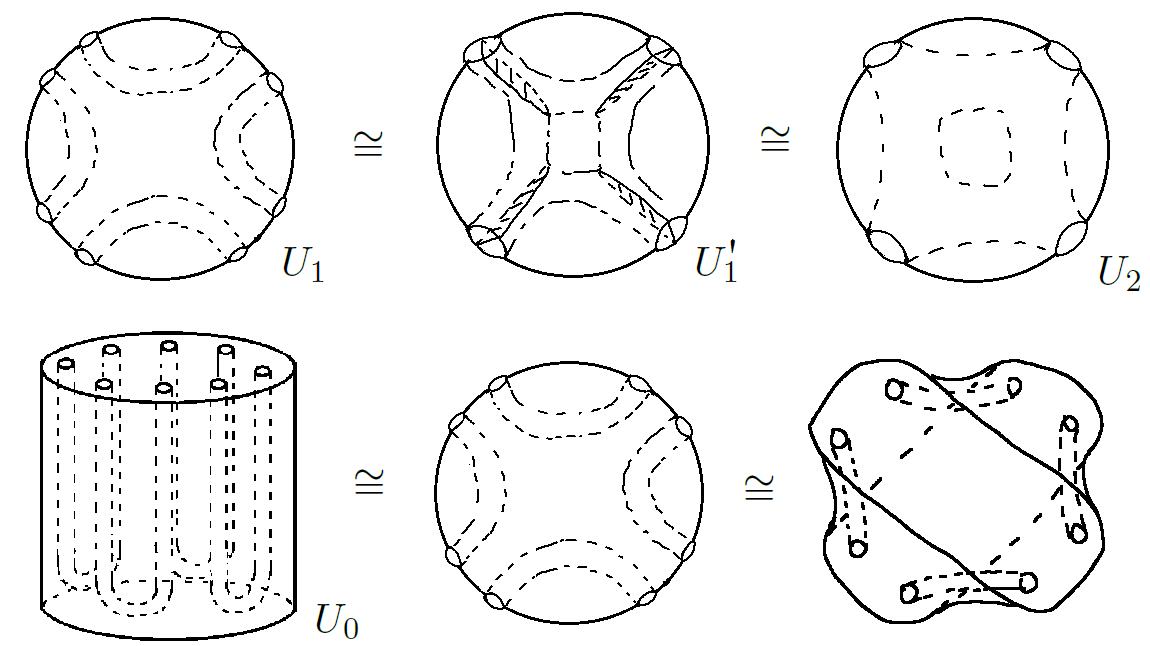}\\
  \caption{$U_0,U_1,U'_1,U_2$ are all diffeomorphic to $\Sigma_{0,5}\times[0,1]$. In the upper row, the diffeomorphism $U'_1\cong U_2$ is given by shrinking the shaded disks towards the center.}\label{fig:handlebody}
\end{figure}

Introduce another three manifolds $U_1,U'_1,U_2$ in Figure \ref{fig:handlebody}.
Then
$$U_0\cong U_2\cong U'_1\cong U_1\cong\Sigma_{0,5}\times[0,1].$$
It is easy to see
\begin{align}
U_0\cong(\Sigma_{0,9}\times[0,1])\cup H'_1\cup \cdots\cup H'_4,  \label{eq:U}
\end{align}
where $H'_1,\ldots,H'_4$ are $2$-handles attached onto $\Sigma_{0,9}\times\{0\}$.

\begin{figure}[h]
  \centering
  \includegraphics[width=3.6cm]{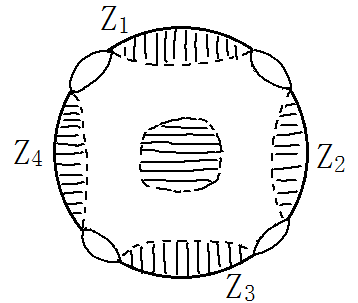}\\
  \caption{The shaded regions stand for ``horizontal" disks embedded in $U_2$. The central disk is the core of a $2$-handle attached to the thickened $\Sigma_{1,4}$; the other ones are denoted by $Z_1,Z_2,Z_3,Z_4$. Cutting $U_2$ along $Z:=\cup_{i=1}^4Z_i$ yields in a $3$-ball.} \label{fig:handlebody2}
\end{figure}

\begin{figure}[h]
  \centering
  \includegraphics[width=13cm]{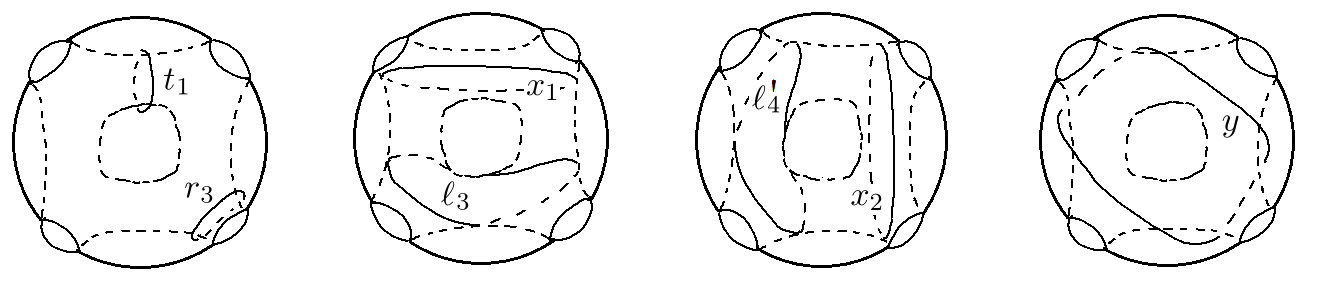}\\
  \caption{Some elements of $\mathcal{S}(\Sigma_{1,4})$, which are also regarded as elements of $\mathcal{S}(U_2)$. The meanings of notations such as $t_3,\ell_2,\ell'_1$ are self-explanatory.} \label{fig:elements}
\end{figure}

As illustrated in Figure \ref{fig:handlebody2}, $U_2$ can be obtained from $\Sigma_{1,4}\times[0,1]$ by attaching a $2$-handle whose core is the central disk, so $\mathcal{S}(U_2)$ is a module over $\mathcal{S}(\Sigma_{1,4})$.
Given a link $\mathbf{l}\subset U$, put $|\mathbf{l}|_i=\#(\mathbf{l}\cap Z_i)$ for $1\le i\le 4$, and define the {\it degree} by $|\mathbf{l}|:=\sum_{i=1}^4|\mathbf{l}|_i$.

Introduced in Figure \ref{fig:elements} are some elements of $\mathcal{S}(\Sigma_{1,4})$, which are also regarded as elements of $\mathcal{S}(U_2)$.
Clearly, $\mathcal{S}(\Sigma_{1,4})$ contains the commutative $R$-subalgebra
$$\mathcal{R}:=R[t_1,t_2,t_3,t_4,r_1,r_2,r_3,r_4].$$

\begin{rmk}
\rm Although $U_0$, $U_1$, $U_2$, $U$ are mutually diffeomorphic, each of them is advantageous in certain specific occasion.
\end{rmk}

\begin{thm}\label{thm:generator}
As a $\mathcal{R}$-module, $\mathcal{S}(U_2)$ is freely generated by
$$\mathcal{X}:=\{1,y\}\cup\{x_1^k,x_2^k\colon k\ge 1\}
\cup\{\ell_1x_1^k,\ell_3x_1^k, \ell_1\ell_3x_1^k, \ell_2x_2^k, \ell_4x_2^k, \ell_2\ell_4x_2^k \colon k\ge 0\}.$$
\end{thm}

\begin{figure}[H]
  \centering
  \includegraphics[width=11.5cm]{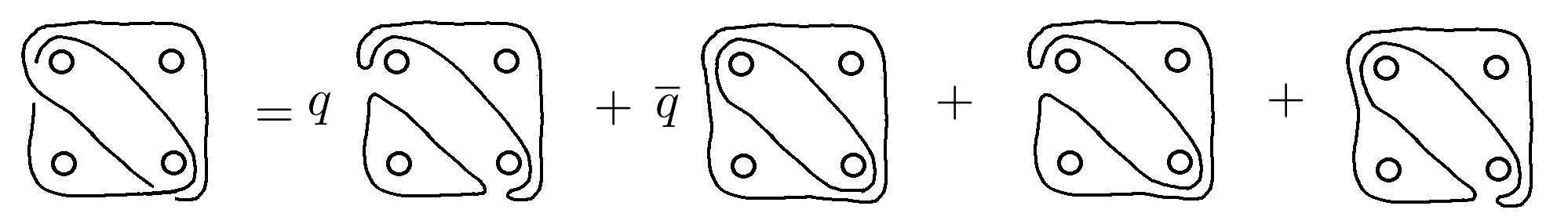}\\
  \caption{$\ell_1\ell_3\approx yx_1$.}\label{fig:l1-times-l3}
\end{figure}

\begin{proof}

Note that $x_1y\approx\ell_1\ell_3$, as illustrated in Figure \ref{fig:l1-times-l3}; similarly, $x_2y\approx\ell_2\ell_4$. Thus $\mathcal{R}\mathcal{X}=\mathcal{R}\mathcal{X}'$, where
$$\mathcal{X}'=\{1,y,\ell_1,\ell_2,\ell_3,\ell_4\}\cup
\{x_1^k, x_1^ky, \ell_1x_1^k, \ell_3x_1^k, x_2^k, x_2^ky, \ell_2x_2^k, \ell_4x_2^k \colon k\ge 1\}.$$

Let $\mathcal{G}=\{t_{i_1\cdots i_r}\colon 1\le i_1<\cdots<i_r\le 4\}$.
As noted in the proof of \cite[Theorem 1.3]{Ch24}, modulo ``simpler" terms and invertible factors, each of $t_{12},t_{23},t_{34},t_{14}$ commutes with any element of $\mathcal{G}$.
Hence as a $R[t_1,t_2,t_3,t_4]$-module, $\mathcal{S}_4$ is generated by elements of the form $u=t_{14}^{j_1}t_{23}^{j_2}at_{12}^{j_3}t_{34}^{j_4}$, with
\begin{align*}
a\in\mathcal{A}=\{t_{13}^k, t_{13}^kt_0, t_{13}^kt_{123}, t_{13}^kt_{134}, t_{24}^k, t_{24}^kt_0, t_{24}^kt_{124}, t_{24}^kt_{234}\colon k\ge 0\}.
\end{align*}

Let $\mathfrak{h}:\mathcal{S}_4\cong\mathcal{S}(U_2)$ denote the $R$-module isomorphism induced by the diffeomorphism $\Sigma_{0,5}\times[0,1]\cong U_2$ (given in Figure \ref{fig:handlebody}). Then
$$\mathfrak{h}(t_{14})=r_1, \qquad  \mathfrak{h}(t_{12})=r_2, \qquad  \mathfrak{h}(t_{23})=r_3, \qquad \mathfrak{h}(t_{34})=r_4.$$
Hence $\mathfrak{h}(u)=r_1^{j_1}r_2^{j_3}r_3^{j_2}r_4^{j_4}\mathfrak{h}(a)$. It suffices to show $\mathfrak{h}(a)\in\mathcal{R}\mathcal{X}'$.

For $f\in\mathcal{S}(U_2)$ and $b\in\mathcal{X}'$, write $f\approx b$ if there exists $0\ne\xi\in R$ such that $f-\xi b=\sum_i\lambda_ic_i$ for some $\lambda_i\in\mathcal{R}$ and $c_i\in\mathcal{X}'$ with $|c_i|<|b|$.

For all $k\ge 0$, we have
\begin{alignat*}{4}
\mathfrak{h}(t_{24}^k)&\approx x_1^k,  \  &\mathfrak{h}(t_{13}^k)&\approx x_2^k, \qquad
&\mathfrak{h}(t_{24}^kt_0)&\approx x_1^ky, \quad  &\mathfrak{h}(t_{13}^kt_0)&\approx x_2^ky;   \\
&\  &\mathfrak{h}(t_{24}^kt_{124})&\approx \ell_1x_1^k, &\ &\  &\mathfrak{h}(t_{13}^kt_{123})&\approx \ell_2x_2^k,  \\
&\  &\mathfrak{h}(t_{24}^kt_{234})&\approx \ell_3x_1^k, &\ &\  &\mathfrak{h}(t_{13}^kt_{134})&\approx \ell_4x_2^k.
\end{alignat*}
Thus, $\mathfrak{h}(a)\in\mathcal{R}\mathcal{X}'$ for each $a\in\mathcal{A}$.
The proof is completed.
\end{proof}

\begin{rmk} \label{rmk:l1-times-l3}
\rm By \cite[Theorem 1.3]{Ch24}, $g^2\sim 0$ for $g\in\{t_{123},t_{124},t_{134},t_{234}\}$. Hence $\ell_i^2\sim 0$ for $1\le i\le 4$.
\end{rmk}

The basis $\mathcal{B}$ (announced in Section 1) for the free $R$-module $\mathcal{S}(U_2)$ is
$$\{t_1^{i_1}t_2^{i_2}t_3^{i_3}t_4^{i_4}r_1^{j_1}r_2^{j_2}r_3^{j_3}r_4^{j_4}b\colon i_1,\ldots,i_4,j_1,\ldots,j_4\in\mathbb{Z}_{\ge0},\ b\in\mathcal{X}\}.$$
Call each element of $\mathcal{B}$ a {\it standard generator}.

\section{Finding all relations}

\begin{nota}
\rm Suppose $\mathbf{b}$ is an oriented arc in $U_2$ with $\partial_{\pm}\mathbf{b}=\mathsf{x}_{\pm}$. Starting at $\mathsf{x}_-$, walk along $\mathbf{b}$ towards $\mathsf{x}_+$; record $i^\ast=i$ (resp. $i^\ast=\overline{i}$) whenever passing through $Z_i$ upwards (resp. downwards).
If $i_1^\ast,\ldots,i_m^\ast$ are all the recordings listed in order, then we denote $\mathbf{b}$ as $\mathbf{b}_{i_1^\ast\cdots i_m^\ast}$.
\end{nota}

For a knot $\mathbf{k}$ with ${\rm fa}:={\rm fa}(\mathbf{k})=\mathbf{b}_{j_1\cdots j_r}$ and a standard generator $u$, denote $u(\mathfrak{s}_i(\mathbf{k})-\mathbf{k})$ by ${\rm Sl}^i_{j_1\cdots j_r}(u)$, and call it a {\it sliding relator of f-degree $r$}.
We allow $r=0$, in which case the notation is ${\rm Sl}^i(u)$.

\begin{lem}\label{lem:bound}
${\rm Sl}(U)$ is generated by sliding relators of f-degree at most $2$.
\end{lem}

\begin{figure}[H]
  \centering
  \includegraphics[width=7.8cm]{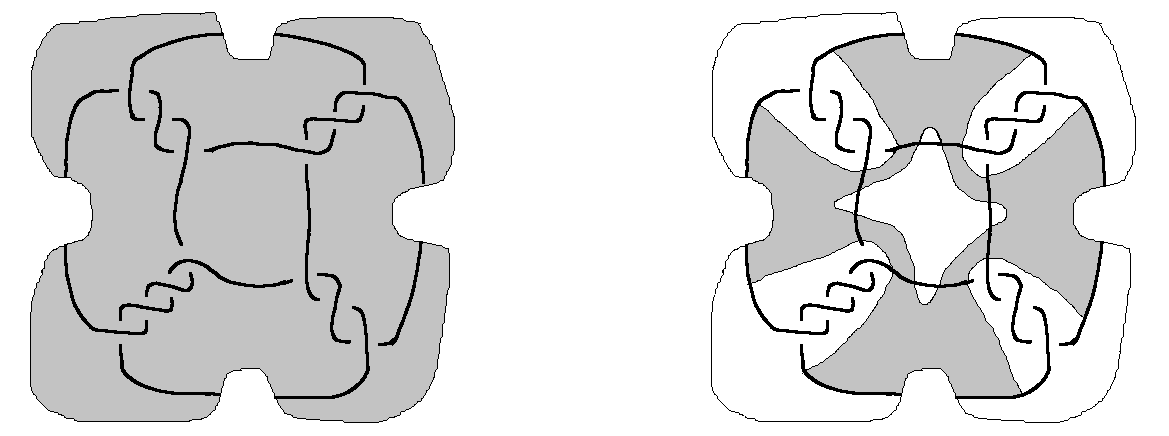}\\
  \caption{Left: the shaded region stands for $U$. Right: the shaded region stands for the submanifold $W$ of $U$, which is diffeomorphic to $\Sigma_{0,9}\times[0,1]$.}\label{fig:W}
\end{figure}

\begin{proof}
By Remark \ref{rmk:slide}, for a link $\mathbf{l}\hookrightarrow U$ to be slid along $\mathbf{o}_i$, we may assume $\mathbf{l}=\mathbf{c}\cup\mathbf{f}\cup\mathbf{l}_1$ and $\mathfrak{s}_i(\mathbf{l})=\mathbf{c}'\cup\mathbf{f}\cup\mathbf{l}_1$, with $\mathbf{c}\cup\mathbf{c}'=\mathbf{o}_i$.

Let $W$ denote the submanifold of $U$ depicted as the shaded region in Figure \ref{fig:W}; it is easy to see $W\cong\Sigma_{0,9}\times[0,1]$.
Observe that with $\partial\mathbf{f}$ fixed, $\mathbf{f}\cup\mathbf{l}_1$ can be isotoped into $W$.
By Lemma \ref{lem:chop}, we may replace $\mathbf{f}\cup\mathbf{l}_1$ by $\sum_t\mathbf{d}_tb_t$, for some $b_t\in\mathcal{S}_8$ and arcs $\mathbf{d}_t$ with $|\mathbf{d}_t|\le 2$, such that in $\mathcal{S}(W)$,
$${\sum}_t(\mathbf{c}\cup\mathbf{d}_t)b_t=\mathbf{l}, \qquad  {\sum}_t(\mathbf{c}'\cup\mathbf{d}_t)b_t=\mathfrak{s}_i(\mathbf{l}).$$
In view of (\ref{eq:U}), we can transform $b_t$ into a $R$-linear combination of standard generators, each of which can be written as a product of elements of $\mathcal{S}(\Sigma_{1,4})$.

For each $t$, isotope $b_t$ into $\Sigma_{1,4}\times[0,1]$, regarding $b_t$ as in $\mathcal{S}(\Sigma_{1,4})$,
then in $\mathcal{S}(U_2)$, whenever necessary we can take an isotopy fixing $\partial \mathbf{d}_t$ to transform $\mathbf{d}_t$ into $\mathbf{b}_{\underline{k}(t)}$, for some $\underline{k}(t)$ of the form $j_1\cdots j_r$ with $r\le 2$. Thus,
$$\mathfrak{s}_i(\mathbf{l})-\mathbf{l}={\sum}_t{\rm Sl}^i_{\underline{k}(t)}(b_t).$$
\end{proof}

\begin{rmk}
\rm It turns out that sometimes it is still more convenient to work with sliding relators of f-degree $\ge 3$.
\end{rmk}

\begin{figure}[h]
  \centering
  \includegraphics[width=9.8cm]{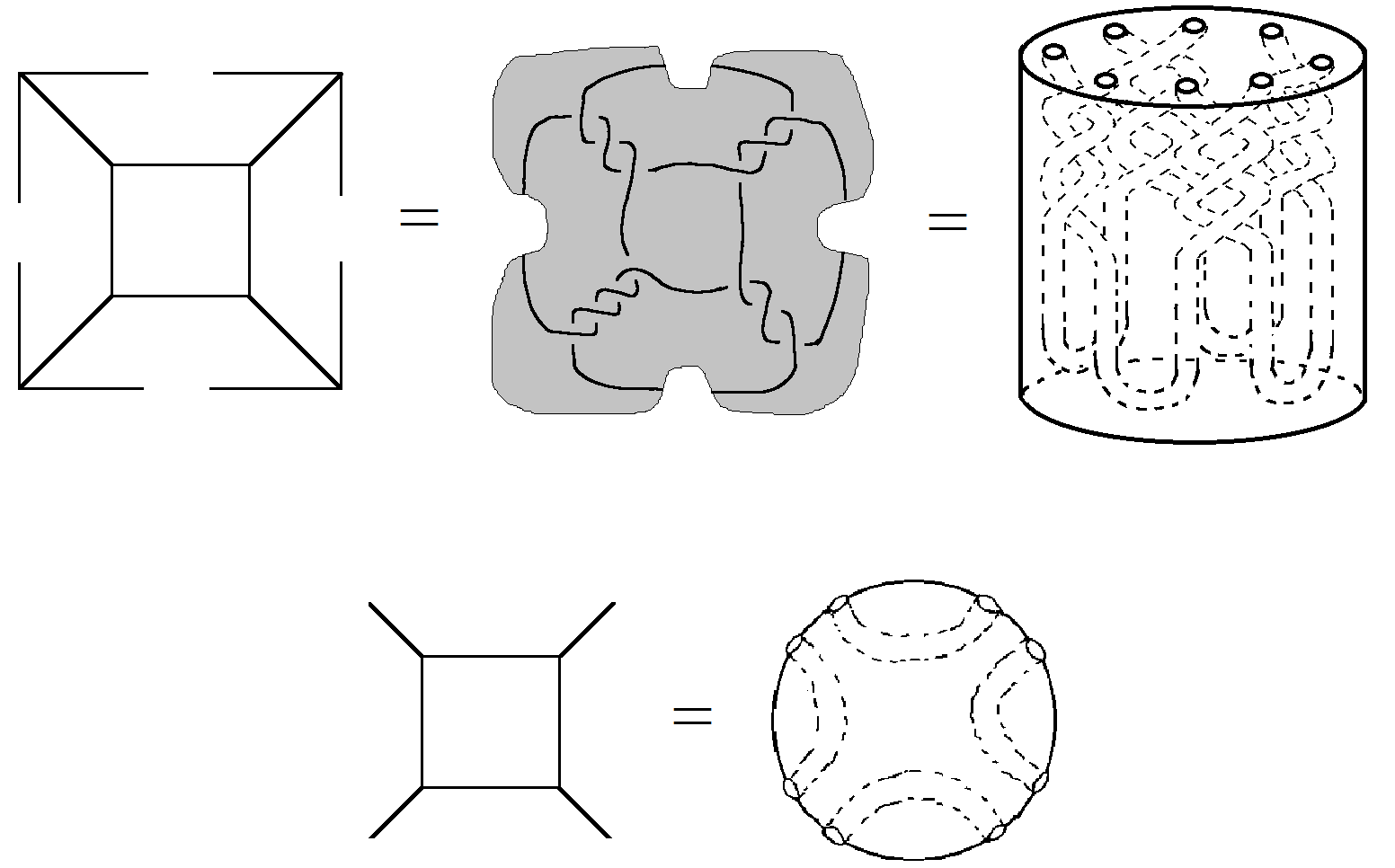}\\
  \caption{Upper: a notation for $U$. Lower: a notation for $U_0$.}\label{fig:symbol}
\end{figure}

\begin{figure}[h]
  \centering
  \includegraphics[width=11cm]{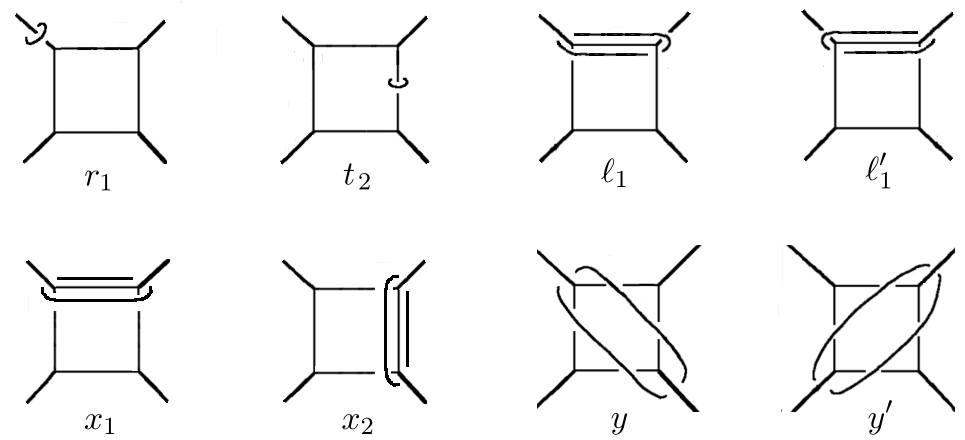}\\
  \caption{The meanings for $x_3$, $\ell_4$, $\ell_2$, etc. should be self-evident.}\label{fig:notation}
\end{figure}

\begin{figure}[h]
  \centering
  \includegraphics[width=8cm]{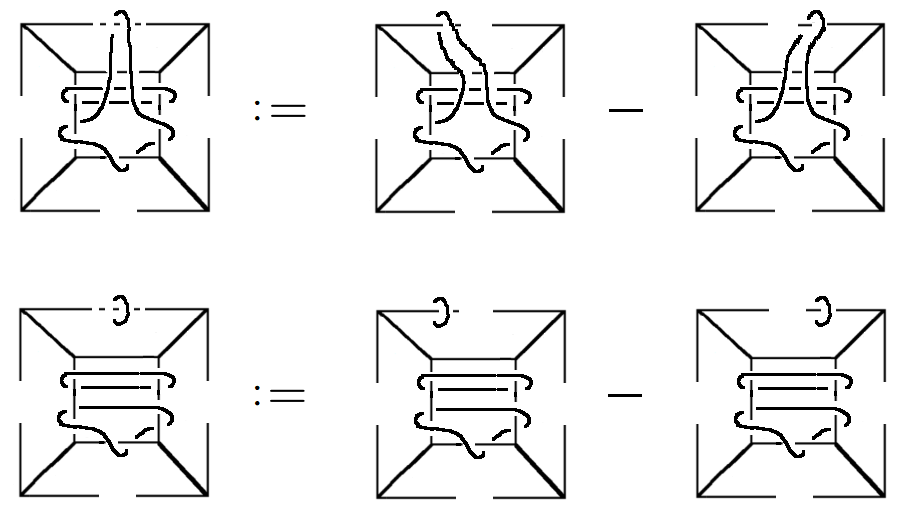}
  \caption{}\label{fig:notation-relator}
\end{figure}

\begin{nota}
\rm Adopt the symbolic notations for $U$ and $U_0$ in Figure \ref{fig:symbol}, and those for elements of $\mathcal{S}(U)$ in Figure \ref{fig:notation}.

Display ${\rm Sl}^i_{j_1\cdots j_r}(u)$ by the manner shown in Figure \ref{fig:notation-relator};
always remember that ${\rm Sl}^i_{j_1\cdots j_r}(u)\in\mathcal{S}(U)$.
\end{nota}


\begin{thm}\label{thm:basic-relator}
As a $\mathcal{R}$-module, ${\rm Sl}(U)$ is generated by ${\rm Sl}^1(u)$, ${\rm Sl}^2(u)$, ${\rm Sl}^3(u)$, ${\rm Sl}^4(u)$ for all $u$, and sliding relators in the following families:
\begin{enumerate}
  \item[\rm(L)] For $k\ge 1$,
        \begin{alignat*}{7}
        &{\rm Sl}^1_1(u), \quad &{\rm Sl}^1_4(u), \quad &{\rm Sl}^1_2(u), \quad &{\rm Sl}^1_{\overline{4}1}(u), \quad
        &{\rm Sl}^1_{\overline{1}2}(u), \quad  &{\rm Sl}^1_{\overline{4}2}(u), \quad  &u=x_1^k, \ell_3x_1^k;  \\
        &{\rm Sl}^3_3(u), \quad &{\rm Sl}^3_2(u),  \quad &{\rm Sl}^3_4(u), \quad &{\rm Sl}^3_{\overline{2}3}(u), \quad
        &{\rm Sl}^3_{\overline{3}4}(u), \quad  &{\rm Sl}^3_{\overline{2}4}(u), \quad  &u=x_1^k, \ell_1x_1^k;  \\
        &{\rm Sl}^2_2(u), \quad &{\rm Sl}^2_1(u),  \quad &{\rm Sl}^2_3(u), \quad &{\rm Sl}^2_{\overline{1}2}(u), \quad
        &{\rm Sl}^2_{\overline{2}3}(u), \quad  &{\rm Sl}^2_{\overline{1}3}(u), \quad  &u=x_2^k, \ell_4x_2^k;  \\
        &{\rm Sl}^4_4(u), \quad &{\rm Sl}^4_3(u), \quad &{\rm Sl}^4_1(u), \quad &{\rm Sl}^4_{\overline{3}4}(u), \quad
        &{\rm Sl}^4_{\overline{4}1}(u), \quad  &{\rm Sl}^4_{\overline{3}1}(u), \quad  &u=x_2^k, \ell_2x_2^k.
        \end{alignat*}
  \item[\rm(G)] For $k\ge 1$,
        \begin{alignat*}{3}
        &{\rm Sl}^1_1(u), \qquad &{\rm Sl}^1_3(u), \qquad &u=x_2^k,\ \ \ell_2x_2^k,\ \ \ell_4x_2^k, \ \ \ell_2\ell_4x_2^k;   \\
        &{\rm Sl}^2_2(u), \qquad &{\rm Sl}^2_4(u), \qquad &u=x_1^k, \ \ \ell_1x_1^k,\ \ \ell_3x_1^k, \ \ \ell_1\ell_3x_1^k.
        \end{alignat*}
  \item[\rm(E)] Sliding relators of low degrees:
        \begin{align*}
        {\rm Sl}^1_3(u), \ \ {\rm Sl}^3_1(u), \ \ {\rm Sl}^2_1(u), \ \ {\rm Sl}^2_3(u), \ \ {\rm Sl}^4_1(u), \ \ {\rm Sl}^4_3(u), \ \
        {\rm Sl}^2_{\overline{1}3}(u), \ \ {\rm Sl}^4_{\overline{1}3}(u),  \\
        u\in\{x_1^k, \ell_1x_1^k, \ell_3x_1^k, \ell_1\ell_3x_1^k \colon k\le 1\};   \\
        {\rm Sl}^2_4(u), \ \ {\rm Sl}^4_2(u), \ \ {\rm Sl}^1_2(u), \ \ {\rm Sl}^1_4(u), \ \ {\rm Sl}^3_2(u), \ \ {\rm Sl}^3_4(u), \ \
        {\rm Sl}^1_{\overline{4}2}(u), \ \ {\rm Sl}^3_{\overline{4}2}(u), \\
        u\in\{x_2^k, \ell_2x_2^k, \ell_4x_2^k, \ell_2\ell_4x_2^k \colon k\le 1\};   \\
        {\rm Sl}^i_{\overline{h},h+1}(u), \qquad i,h\in\{1,2,3,4\}, \ \  u\in\{1,\ell_1,\ell_3,\ell_1\ell_3,\ell_2,\ell_4,\ell_2\ell_4\}.
        \end{align*}
\end{enumerate}
\end{thm}

\begin{rmk}
\rm To make the statement neat, we allow the (E) and (G) family to overlap; this does not matter, since the (E) family will not be taken into account.
\end{rmk}

Use {\it ${\rm L}_i$-relators} (resp. {\it ${\rm G}_i$-relators}) to name the sliding relators given in the $i$-th row of the (L) family (resp. the (G) family). Let $\mathcal{L}_i$ (resp. $\mathcal{G}_i$) denote the set of ${\rm L}_i$-relators (resp. ${\rm G}_i$-relators).
When $\mathfrak{f}$ is a ${\rm L}_i$-relator (resp. ${\rm G}_i$-relator), call $\mathfrak{f}=0$ a {\it ${\rm L}_i$-relation} (resp. {\it ${\rm G}_i$-relation}).


Given $a\in\mathcal{S}(U)$, say $a=\sum_i\lambda_i\mathbf{l}_i$, with $0\ne\lambda_i\in R$ and $\mathbf{l}_i$ being a link, put
$|a|=\max_i|\mathbf{l}_i|$.

For $a,a'\in\mathcal{S}(U)$, write $a\equiv a'$ if $a-a'$ equals a $\mathcal{R}$-linear combination of links of degree smaller than $\min\{|a|,|a'|\}$; write $a\equiv 0$ if $a$ equals a $\mathcal{R}$-linear combination of links of degree smaller than $|a|$.

\bigskip

To efficiently study $\mathcal{S}(U)$, we propose a technique named ``embedding surfaces".
Take an embedding $\iota:\Sigma_{0,n+1}\times[0,1]\hookrightarrow U$.
Given $f\in\mathcal{S}(U)$, isotope it into ${\rm Im}(\iota)$ whenever possible,
then use relations induced from $\mathcal{S}_n$ to ``normalize" $f$, in certain sense.

\begin{figure}[h]
  \centering
  \includegraphics[width=11cm]{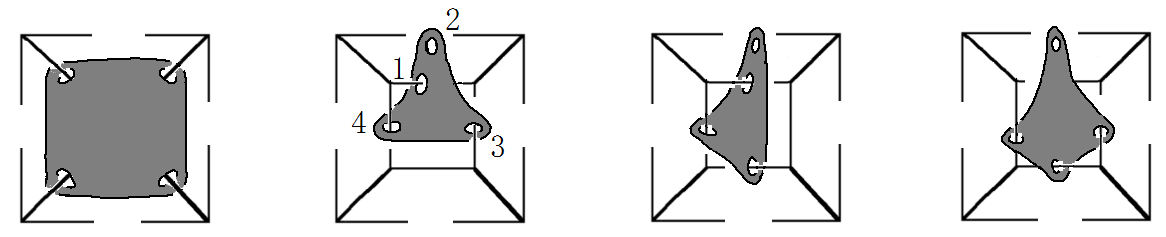}\\
  \caption{Embeddings $\Sigma_{0,5}\hookrightarrow U$. The first one is $\iota_0$, and the second is $\iota_1$.}\label{fig:embedding}
\end{figure}

\begin{figure}[h]
  \centering
  \includegraphics[width=11.5cm]{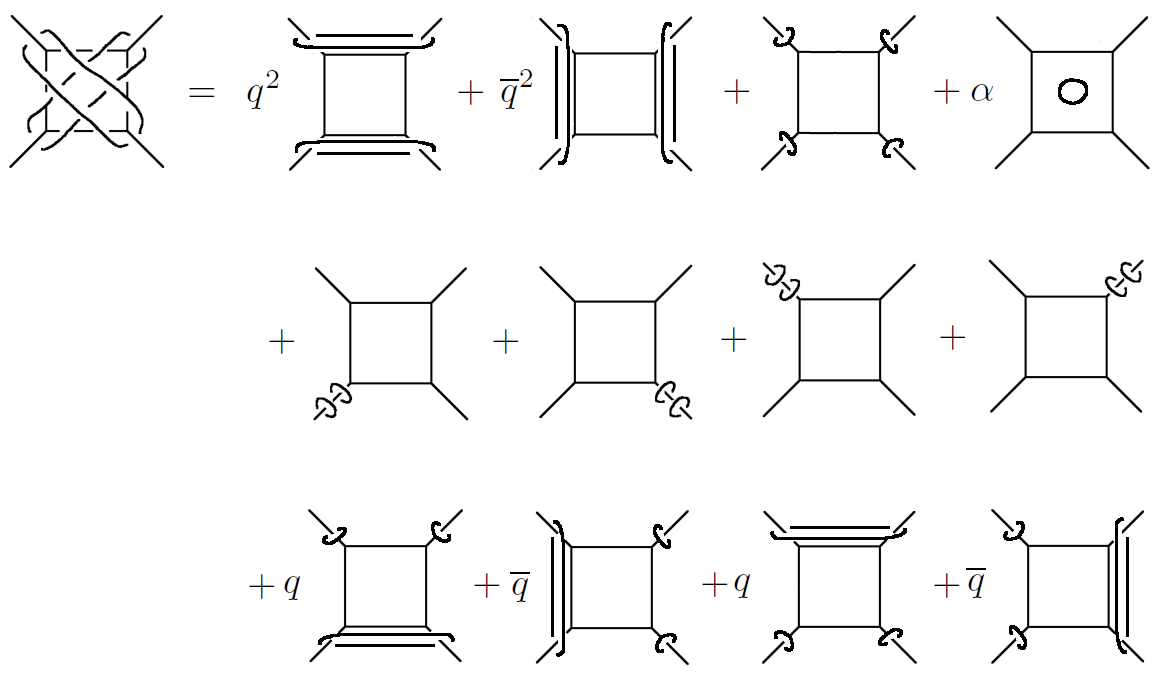}\\
  \caption{An application of the formula in Figure \ref{fig:13-times-24}.}\label{fig:y1-times-y2}
\end{figure}

As an example, taking $\iota_0$ as shown in the first part of Figure \ref{fig:embedding} and applying the formula given by Figure \ref{fig:13-times-24}, we deduce an equation in $\mathcal{S}(U)$; see Figure \ref{fig:y1-times-y2}.

Given $\varepsilon_1,\varepsilon_2,\varepsilon_3,\varepsilon_4\in\mathbb{Z}_{\ge0}$, let
$\widetilde{\mathcal{K}}_{\varepsilon_1\varepsilon_2\varepsilon_3\varepsilon_4}$ denote the $R$-submodule of $\iota_\ast(\mathcal{S}_n)$ generated by $\mathbf{l}$ with $|\mathbf{l}|_i\le \varepsilon_i$ for $1\le i\le 4$, and let $\mathcal{K}_{\varepsilon_1\varepsilon_2\varepsilon_3\varepsilon_4}$ denote the quotient of $\widetilde{\mathcal{K}}_{\varepsilon_1\varepsilon_2\varepsilon_3\varepsilon_4}$ by its submodule generated by $\mathbf{l}$ with $|\mathbf{l}|<\sum_{i=1}^4\varepsilon_i$.

Take an embedding $\iota_1:\Sigma_{0,5}\hookrightarrow U$ as illustrated by the second part of Figure \ref{fig:embedding}, where the hole labeled with $k$ is the image of the $k$-th hole of $\Sigma_{0,5}$.
For $t_{i_1\cdots i_r}\in\mathcal{S}_4$, let $t'_{i'_1\cdots i'_r}\in\mathcal{S}(U)$ denote its image under $(\iota_1)_\ast$, if the hole labeled with $k$ is the $i'_k$-th hole. Here we adopt the convention that the upper-middle hole is the $0$-th one.
By Lemma \ref{lem:basic-reduce},
\begin{align*}
\alpha t'_{0241}&\equiv t'_{04}t'_{12}-\overline{q}^2t'_{01}t'_{24}-q^2t'_{14}t'_{02},  \\
t'_{104}t'_{12}&\equiv\overline{q}^2t'_{124}t'_{10}+\overline{q}^2t'_{14}t'_{102}.
\end{align*}
Note that $t'_{12},t'_{14}\in\mathcal{R}$, so $t'_{0241}\equiv -\overline{q}^2\alpha^{-1}t'_{01}t'_{24}$, implying $$\dim_{\mathcal{R}}\mathcal{K}_{1101}=1,$$
by which we mean that the nonzero elements of $\mathcal{K}_{1101}$ are $\mathcal{R}$-multiples of each other. Moreover, $t'_{124}t'_{10}\equiv 0$, implying
$$\mathcal{K}_{2101}=0.$$
Use the same embedding. By Lemma \ref{lem:basic-reduce},
$t'_{04}t'_{124}\equiv q^2t'_{14}t'_{024}+\overline{q}^2t'_{24}t'_{104}$,
so
$$\dim_{\mathcal{R}}\mathcal{K}_{1102}=1.$$

Similarly, using the other embeddings in Figure \ref{fig:embedding}, we can show
$$\dim_{\mathcal{R}}\mathcal{K}_{1011}=\dim_{\mathcal{R}}\mathcal{K}_{0111}=\dim_{\mathcal{R}}\mathcal{K}_{2011}
=\dim_{\mathcal{R}}\mathcal{K}_{0112}=\dim_{\mathcal{R}}\mathcal{K}_{1021}=1.$$

\bigskip

Call ${\rm Sl}^i_{j_1\cdots j_r}(u)$ {\it redundant} if it equals a $\mathcal{R}$-linear combination of the sliding relators given in Theorem \ref{thm:basic-relator}. Call ${\rm Sl}^i_{j_1\cdots j_r}(u)$ {\it reducible} if it equals a $\mathcal{R}$-linear combination of
sliding relators of smaller degrees. In particular, ${\rm Sl}^i_{j_1\cdots j_r}(u)$ is reducible if $j_s\in\{j_t,\overline{j_t}\}$ for some $1\le s<t\le r$.

\begin{lem}\label{lem:reduce-0}
{\rm(a)} 
${\rm Sl}^1_{\overline{3}41}(x_2)\equiv -q^2{\rm Sl}^2_{\overline{2}3}(\ell_4)$;
{\rm(b)} ${\rm Sl}^1_{\overline{1}3}(1)\equiv{\rm Sl}^4_{\overline{3}4}(1)+{\rm Sl}^2_{\overline{2}3}(1)$.
\end{lem}

\begin{proof}
(a) As clarified by Figure \ref{fig:automatic}, we can write
\begin{align*}
{\rm Sl}^1_{\overline{3}41}(x_2)=-{\rm Sl}^2_{\overline{3}42}(x_2)+\beta_1{\rm Sl}^3_{\overline{3}43}(x_2)+\beta_2{\rm Sl}^4_{\overline{3}44}(x_2)
\end{align*}
for some $\beta_1,\beta_2$. 
The sliding relators ${\rm Sl}^3_{\overline{3}43}(x_2)$ and ${\rm Sl}^4_{\overline{3}44}(x_2)$ are reducible.
As shown in Figure \ref{fig:reduce-a782},
${\rm Sl}^2_{\overline{3}42}(x_2)\equiv -{\rm Sl}^2_{\overline{3}2}(\ell_4)$.
Thus, ${\rm Sl}^1_{\overline{3}41}(x_2)\equiv -q^2{\rm Sl}^2_{\overline{2}3}(\ell_4)$.

(b) Similarly as in (a), we can write
$${\rm Sl}^1_{\overline{1}3}(1)={\rm Sl}^2_{\overline{2}3}(1)+{\rm Sl}^3(1)+{\rm Sl}^4_{\overline{3}4}(1).$$
Since ${\rm Sl}^3(1)$ is reducible, the assertion follows.

\begin{figure}[h]
  \centering
  \includegraphics[width=10.5cm]{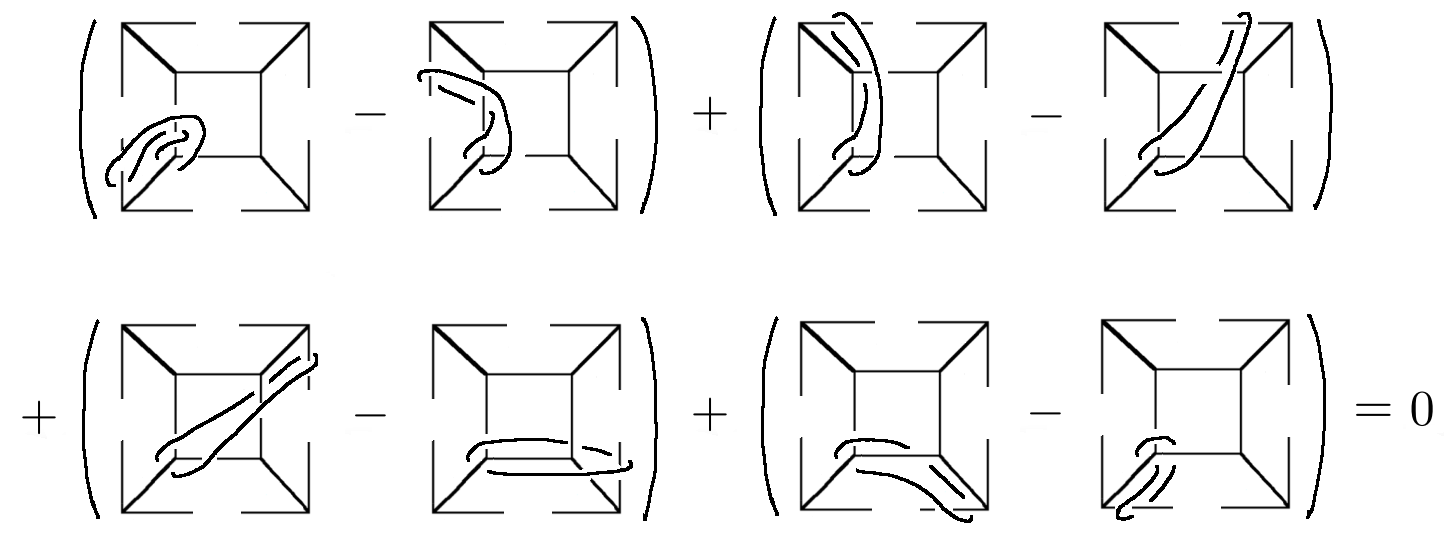}\\
  \caption{}\label{fig:automatic}
\end{figure}

\begin{figure}[h]
  \centering
  \includegraphics[width=11.5cm]{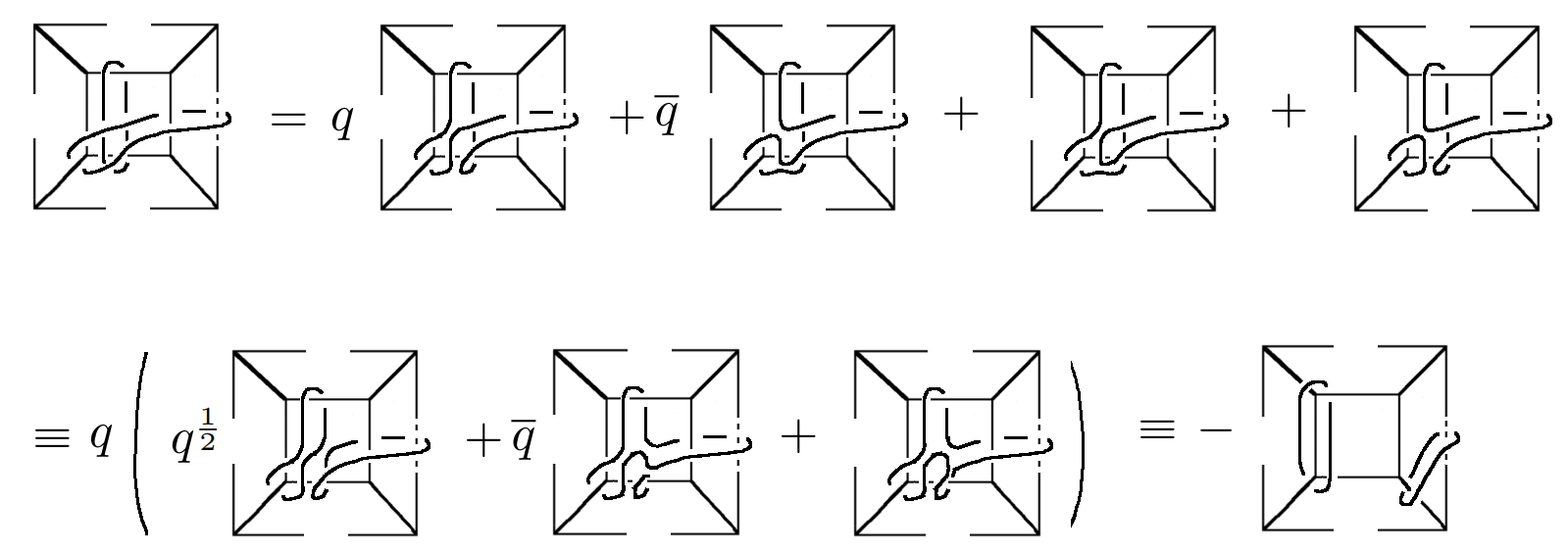}\\
  \caption{}\label{fig:reduce-a782}
\end{figure}

\end{proof}

\begin{figure}[h]
  \centering
  \includegraphics[width=11.5cm]{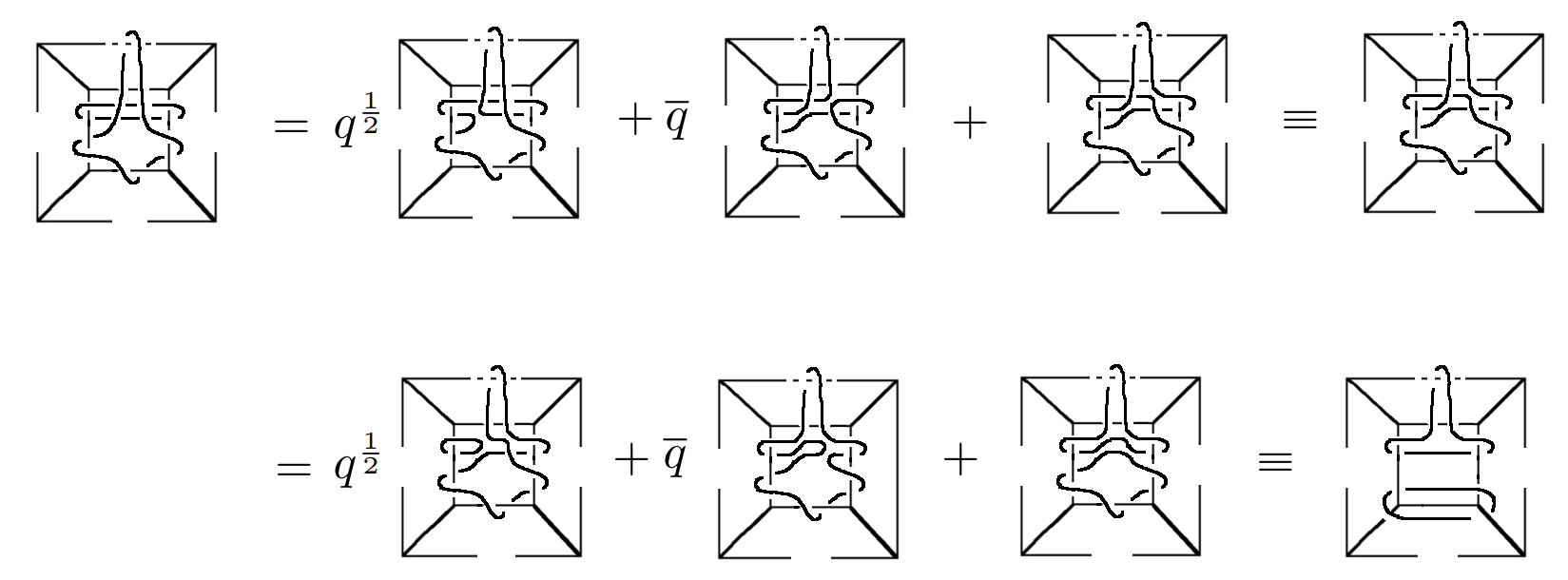}\\
  \caption{}
  \label{fig:reduce-a6-1-2}
\end{figure}

\begin{figure}[h]
  \centering
  \includegraphics[width=11.5cm]{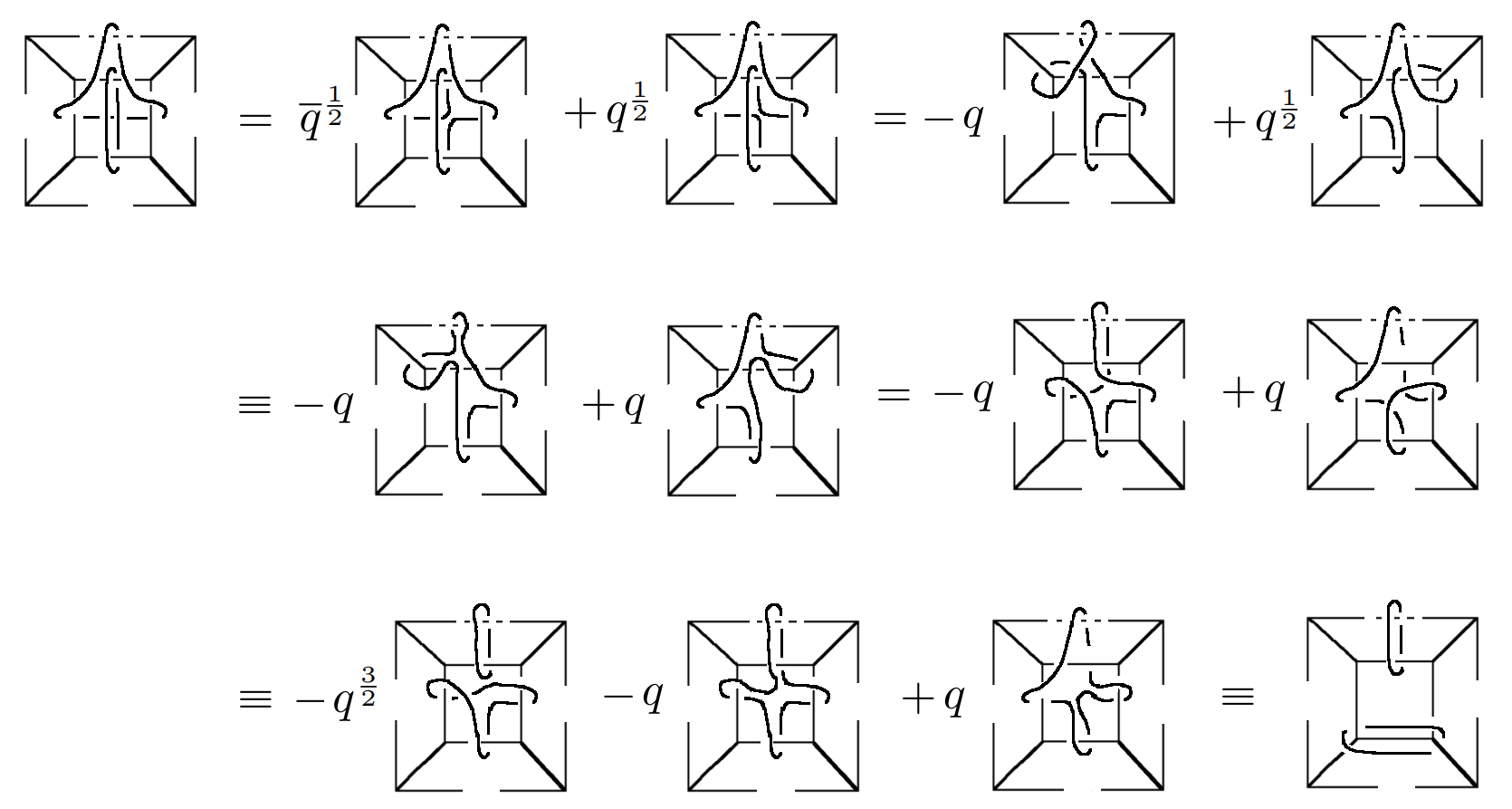}\\
  \caption{}\label{fig:reduce-a14-1}
\end{figure}

\begin{proof}[Proof of Theorem \ref{thm:basic-relator}]
By Lemma \ref{lem:bound}, it suffices to deal with ${\rm Sl}^1_{j_1\cdots j_r}(u)$ with $r=1,2$.
The proof consists of two steps. 

{\bf Step 1}. We investigate ${\rm Sl}^1_{jk}(u)$ case-by-case; the ${\rm Sl}^i_{jk}(u)$'s for $i\in\{2,3,4\}$ can be handled in a parallel way.

By Lemma \ref{lem:reduce-0} (b), ${\rm Sl}^1_{\overline{1}3}(u)$ is redundant for any $u$.

Recall the notation ${\rm fa}$ introduced in the beginning of this section.
\begin{enumerate}
  \item When ${\rm fa}=\mathbf{b}_1$, due to $\mathcal{K}_{2101}=0$, we have ${\rm Sl}^1_1(\ell_1)\equiv 0$.
        Hence ${\rm Sl}^1_1(\ell_1x_1^k)$ is redundant for all $k\ge 0$, so is ${\rm Sl}^1_1(\ell_1\ell_3x_1^k)$.
  \item Suppose ${\rm fa}=\mathbf{b}_4$.
    \begin{enumerate}
    \item It follows from $\dim_{\mathcal{R}}\mathcal{K}_{1102}=1$ that ${\rm Sl}^1_4(\ell_1)\approx{\rm Sl}^1_{\overline{4}1}(x_1)$, by which
          we mean ${\rm Sl}^1_4(\ell_1)\equiv \beta{\rm Sl}^1_{\overline{4}1}(x_1)$ for some $0\ne\beta \in R$.
          Hence ${\rm Sl}^1_4(\ell_1x_1^k)$ is redundant for all $k\ge 0$.
          Consequently, ${\rm Sl}^1_4(\ell_1\ell_3)\approx{\rm Sl}^1_{\overline{4}1}(\ell_3x_1)$, so that
          ${\rm Sl}^1_4(\ell_1\ell_3x_1^k)$ is also redundant for all $k\ge 0$.
    \item Due to $\dim_{\mathcal{R}}\mathcal{K}_{1011}=1$, we have ${\rm Sl}^1_4(x_2)\approx{\rm Sl}^1_{\overline{3}41}(1)$.
          Hence by Lemma \ref{lem:reduce-0} (a), ${\rm Sl}^1_4(x_2^k)$, ${\rm Sl}^1_4(\ell_2x_2^k)$, ${\rm Sl}^1_4(\ell_4x_2^k)$,
          ${\rm Sl}^1_4(\ell_2\ell_4x_2^k)$ are all redundant for all $k\ge 2$.
    \end{enumerate}
        Thus, it suffices to consider ${\rm Sl}^1_4(x_1^k)$ and ${\rm Sl}^1_4(\ell_3x_1^k)$.
  \item Similarly, when ${\rm fa}=\mathbf{b}_2$, it suffices to consider ${\rm Sl}^1_2(x_1^k)$ and
        ${\rm Sl}^1_2(\ell_3x_1^k)$.
  \item Suppose ${\rm fa}=\mathbf{b}_3$ and $k\ge 2$.

        Due to $\dim_{\mathcal{R}}\mathcal{K}_{0111}=1$, we have ${\rm Sl}^1_3(x_1)\approx{\rm Sl}^1_{234}(1)$; by Figure \ref{fig:reduce-a6-1-2}, ${\rm Sl}^1_{234}(x_1)\approx{\rm Sl}^1_{\overline{4}2}(\ell_3)$, so
        ${\rm Sl}^1_3(x_1^2)\approx{\rm Sl}^1_{\overline{4}2}(\ell_3)$. Hence ${\rm Sl}^1_3(x_1^k)$ is redundant.

        Consequently, ${\rm Sl}^1_3(\ell_1x_1^k)$, ${\rm Sl}^1_3(\ell_3x_1^k)$ and ${\rm Sl}^1_3(\ell_1\ell_3x_1^k)$ are also redundant.

        Hence we only need to consider ${\rm Sl}^1_3(ux_2^k)$ for $u\in\{1,\ell_2,\ell_4,\ell_2\ell_4\}$.
  \item Suppose ${\rm fa}=\mathbf{b}_{\overline{4}1}$ and $k\ge 1$.
    \begin{enumerate}
    \item By Remark \ref{rmk:l1-times-l3}, $\ell_1^2\equiv 0$.
          So ${\rm Sl}^1_{\overline{4}1}(\ell_1x_1)\equiv q^2{\rm Sl}^1_4(\ell_1^2)\equiv 0$, implying that
          ${\rm Sl}^1_{\overline{4}1}(\ell_1x_1^k)$ is redundant.
          Consequently, ${\rm Sl}^1_{\overline{4}1}(\ell_1\ell_3x_1^k)$ is redundant.
    \item ${\rm Sl}^1_{\overline{4}1}(x_2^k)$ is redundant, since $\dim_{\mathcal{R}}\mathcal{K}_{2011}=1$ implies
          ${\rm Sl}^1_{\overline{4}1}(x_2)\approx{\rm Sl}^1_1(\ell'_4)$. Hence ${\rm Sl}^1_{\overline{4}1}(\ell_2x_2^k)$,
          ${\rm Sl}^1_{\overline{4}1}(\ell_4x_2^k)$, ${\rm Sl}^1_{\overline{4}1}(\ell_2\ell_4x_2^k)$ are redundant.
    \end{enumerate}
        Thus, we only need to consider ${\rm Sl}^1_{\overline{4}1}(x_1^k)$ and ${\rm Sl}^1_{\overline{4}1}(\ell_3x_1^k)$.
  \item Similarly, when ${\rm fa}=\mathbf{b}_{\overline{1}2}$, it suffices to consider ${\rm Sl}^1_{\overline{1}2}(x_1^k)$,
        ${\rm Sl}^1_{\overline{1}2}(\ell_3x_1^k)$.
  \item Suppose ${\rm fa}=\mathbf{b}_{\overline{3}4}$ and $k\ge 1$.
    \begin{enumerate}
    \item ${\rm Sl}^1_{\overline{3}4}(x_1^k)$ is redundant, since $\dim_{\mathcal{R}}\mathcal{K}_{0112}=1$ implies
          ${\rm Sl}^1_{\overline{3}4}(x_1)\approx{\rm Sl}^1_4(\ell'_3)$.
          Consequently, ${\rm Sl}^1_{\overline{3}4}(\ell_1x_1^k)$, ${\rm Sl}^1_{\overline{3}4}(\ell_3x_1^k)$,
          ${\rm Sl}^1_{\overline{3}4}(\ell_1\ell_3x_1^k)$ are also redundant.
    \item ${\rm Sl}^1_{\overline{3}4}(x_2^k)$ is redundant, since $\dim_{\mathcal{R}}\mathcal{K}_{1021}=1$ implies
          ${\rm Sl}^1_{\overline{3}4}(x_2)\approx{\rm Sl}^1_{\overline{3}}(\ell_4)$.
          Consequently, ${\rm Sl}^1_{\overline{3}4}(\ell_2x_2^k)$, ${\rm Sl}^1_{\overline{3}4}(\ell_4x_2^k)$,
          ${\rm Sl}^1_{\overline{3}4}(\ell_2\ell_4x_2^k)$ are also redundant.
    \end{enumerate}
        Hence ${\rm Sl}^1_{\overline{3}4}(u)$ is redundant except for $u\in\{1,\ell_1,\ell_3,\ell_1\ell_3,\ell_2,\ell_4,\ell_2\ell_4\}$.
  \item The situation for ${\rm fa}=\mathbf{b}_{\overline{2}3}$ is similar.
  \item Suppose ${\rm fa}=\mathbf{b}_{\overline{4}2}$.
    \begin{enumerate}
    \item For any $k\ge 0$, ${\rm Sl}^1_{\overline{4}2}(\ell_1x_1^k)$ is redundant. Indeed, $\dim_{\mathcal{R}}\mathcal{K}_{1101}=1$ implies
          ${\rm Sl}^1_{\overline{1}42}(1)\approx{\rm Sl}^1_1(x_1)$; by Figure \ref{fig:reduce-a14-2},
          ${\rm Sl}^1_{\overline{4}2}(\ell_1)\equiv{\rm Sl}^1_{\overline{1}42}(x_1)\approx{\rm Sl}^1_1(x_1^2)$.
          Consequently, ${\rm Sl}^1_{\overline{4}2}(\ell_1\ell_3x_1^k)$ is also redundant.
    \item When $k\ge 1$, ${\rm Sl}^1_{\overline{4}2}(x_2^k)$ is redundant, as by Figure \ref{fig:reduce-a14-1},
          ${\rm Sl}^1_{\overline{4}2}(x_2)\equiv{\rm Sl}^1_1(\ell_3)$.
          So ${\rm Sl}^1_{\overline{4}2}(\ell_2x_2^k)$, ${\rm Sl}^1_{\overline{4}2}(\ell_4x_2^k)$,
          ${\rm Sl}^1_{\overline{4}2}(\ell_2\ell_4x_2^k)$ are also redundant.
  \end{enumerate}
        Hence we only need to consider ${\rm Sl}^1_{\overline{4}2}(x_1^k)$ and ${\rm Sl}^1_{\overline{4}2}(\ell_3x_1^k)$.
\end{enumerate}

\begin{figure}[h]
  \centering
  \includegraphics[width=12.5cm]{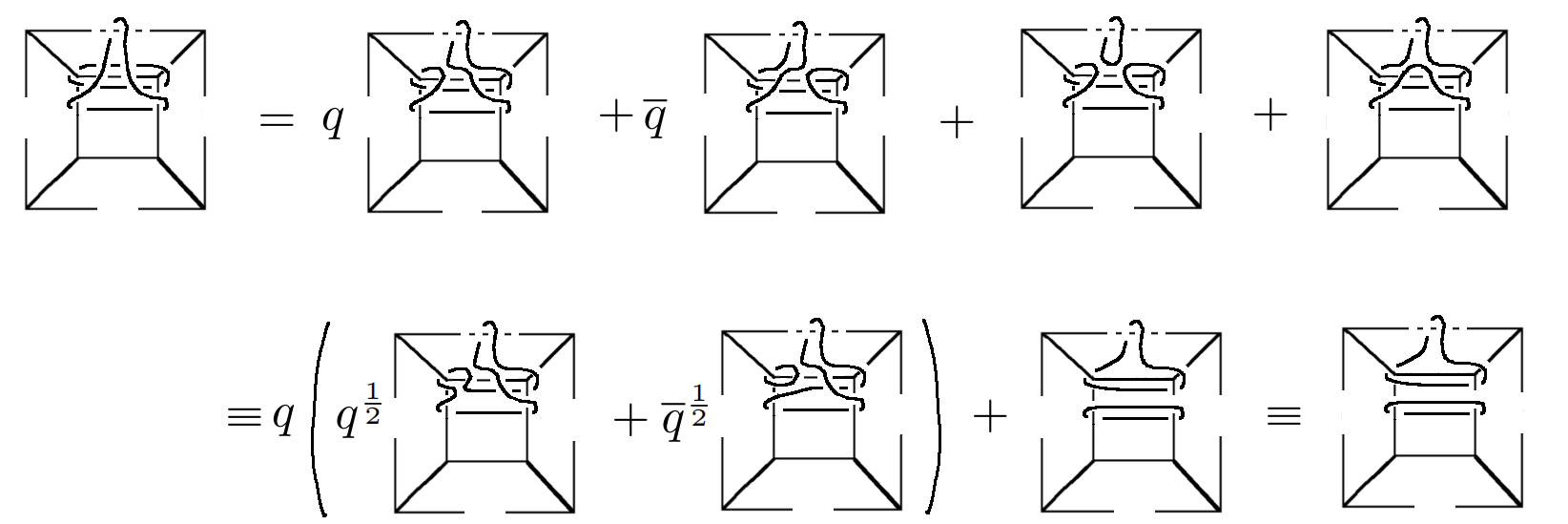}\\
  \caption{}\label{fig:reduce-a14-2}
\end{figure}

\begin{figure}[H]
  \centering
  \includegraphics[width=12cm]{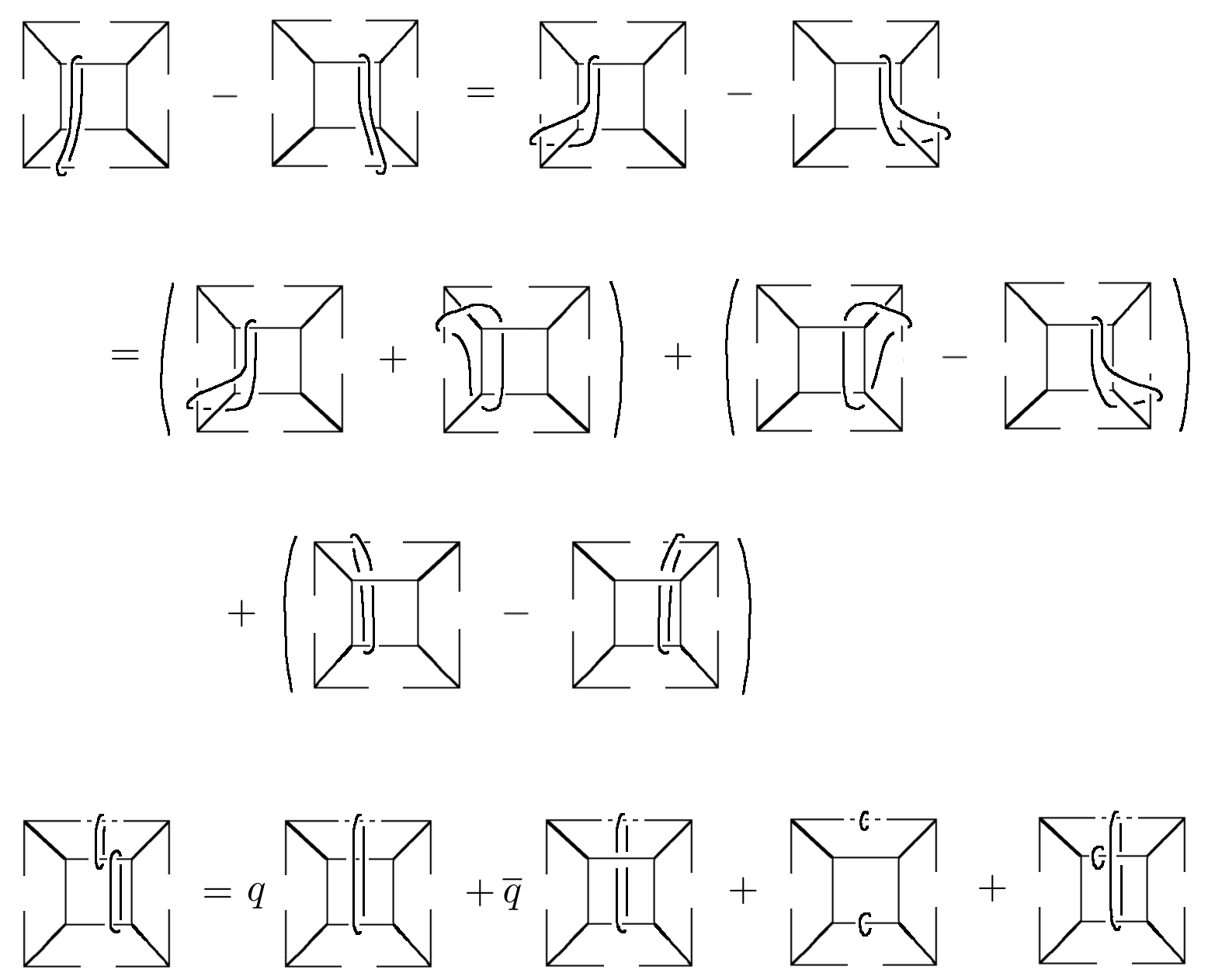}\\
  \caption{}\label{fig:reduce-0}
\end{figure}

{\bf Step 2}. From Figure \ref{fig:reduce-0} we can see that modulo ${\rm L}_2$- and ${\rm L}_4$-relators,
${\rm Sl}^3_1(1)\equiv q^2{\rm Sl}^1_3(1)-q{\rm Sl}^1_1(x_2)$, so ${\rm Sl}^3_1(u)$ is redundant.
Moreover, similarly as in the proof of Lemma \ref{lem:reduce-0} (a), we may write
$${\rm Sl}^3_3(u)=-{\rm Sl}^1_1(u)-{\rm Sl}^2_2(u)-{\rm Sl}^4_4(u),$$
hence ${\rm Sl}^3_3(u)$ is redundant for $u\in\{x_2^k,\ell_2x_2^k,\ell_4x_2^k,\ell_2\ell_4x_2^k\}$.

For similar reasons, ${\rm Sl}^4_2(u)$ is redundant for all $u$, and ${\rm Sl}^4_4(u)$ is redundant for $u\in\{x_1^k,\ell_1x_1^k,\ell_3x_1^k,\ell_1\ell_3x_1^k\}$.

The proof is completed.
\end{proof}

\section{Explicit computations}

By Lemma \ref{lem:sliding} and Theorem \ref{thm:basic-relator}, we can determine the relations in $\mathcal{S}(Y)$ as follows.
For each sliding relator $\mathfrak{f}$ in Theorem \ref{thm:basic-relator}, write $\mathfrak{f}$ as a $R$-linear combination of standard generators (given in Theorem \ref{thm:generator}), to express $\mathfrak{f}=0$ as a $\mathcal{R}$-linear dependence relations among the elements of $\mathcal{X}$.


The relations ${\rm Sl}^i(u)=0$, $1\le i\le 4$ simply give rise to $t_3=t_1$ and $t_4=t_2$.

To simplify the computation, we set $t_1=t_2$, i.e. to work in
$$\mathcal{S}_\diamond(Y)=\mathcal{S}(Y)/(t_1-t_2)\mathcal{S}(Y).$$
Let $t$ denote the common image of $t_1,t_2$ under $\mathcal{S}(Y)\twoheadrightarrow\mathcal{S}_\diamond(Y)$. Then $R[t]$ acts on $\mathcal{S}_\diamond(Y)$ in a self-evident way.

The main goal is to show that $\mathcal{S}_\diamond(Y)$ is not finitely generated as a $R[t]$-module.
Then Theorem \ref{thm:main} will follow immediately.

\subsection{Set up}

\begin{figure}[h]
  \centering
  \includegraphics[width=12.5cm]{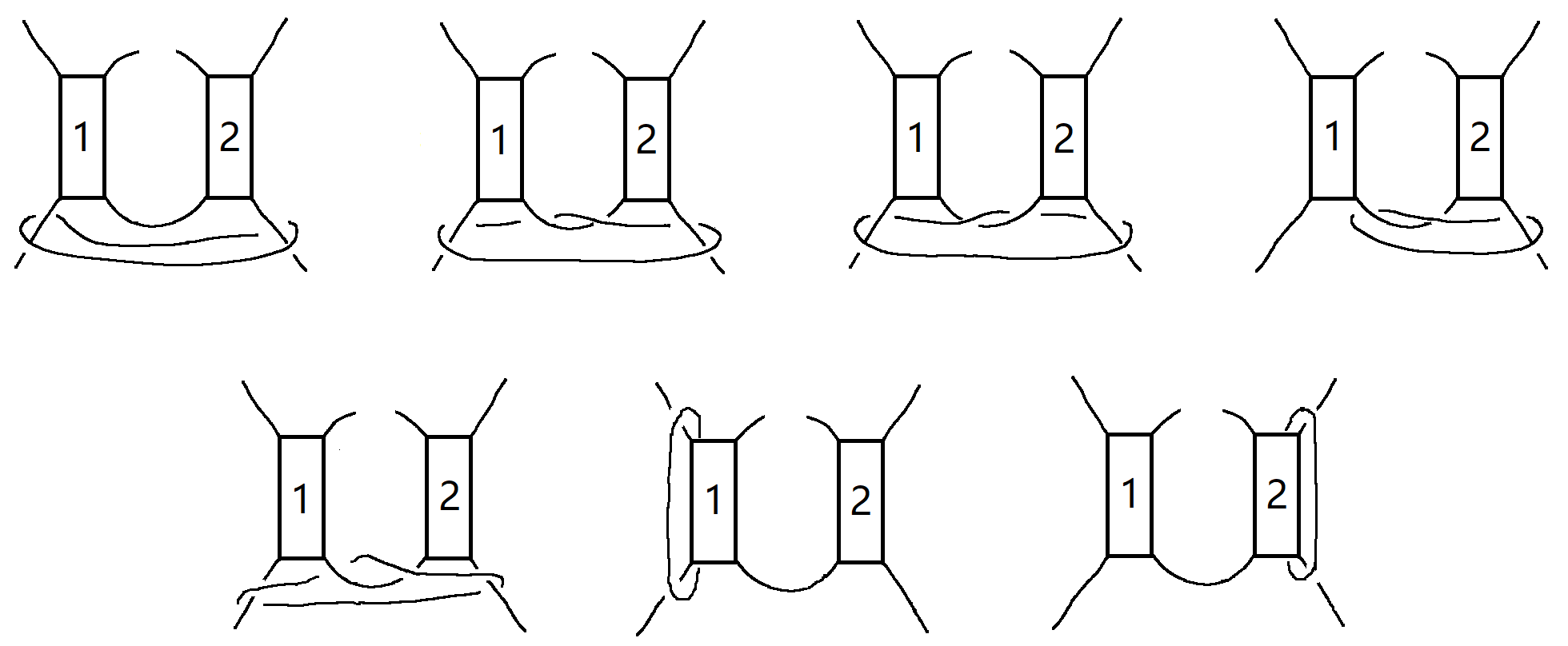}\\
  \caption{First row (from left to right): $p_{14}$, $p_{124}$, $p_{134}$, $p_{24}$. Second row: $p'_{134}$, $s_1$, $s_2$.}\label{fig:notation-2}
\end{figure}

For convenience, introduce auxiliary notations $p_{124}$, $p_{24}$ etc; see Figure \ref{fig:notation-2} for instances.
The following identities are easy to deduce:
\begin{align*}
p_{12}=r_1, \qquad p_{13}=-q^2r_1-qt^2, \qquad p_{34}=r_2, \qquad  p_{24}=-\overline{q}^2r_2-\overline{q}t^2,  \\
p_{23}=-\alpha, \qquad p_{123}=p_{234}=t, \qquad p_{14}=p_{1234}=x_1,  \\
p_{124}=-\overline{q}^2\ell_1-\overline{q}tr_1, \qquad  p_{134}=-q^2p_{124}-qtp_{14}=\ell_1+qt(r_1-x_1). 
\end{align*}

For $r=\mu+\mu^{-1}$, $n\ge 0$ and $k\in\mathbb{Z}$, put
\begin{align*}
\gamma_k(r)&=\frac{\mu^k-\mu^{-k}}{\mu-\mu^{-1}},   \\
\eta_k^n(r)
&={\sum}_{j=0}^n\tbinom{n}{j}q^{2(2j-n)+k}\gamma_{2j-n+k}(r),  \\
c_k^n(r)&=\frac{\overline{q}\eta_{k+1}^n(r)-q\eta_k^n(r)-r^n}{r-\alpha}.
\end{align*}
See Appendix A for details on these functions.

Introduce
\begin{align}
\lambda_k^n(r)=-\eta_{k+1}^n(r)+\eta_k^n(r)-c_k^n(r)t^2.
\end{align}
It follows from (\ref{eq:relation-eta-2}) and (\ref{eq:relation-c}) that
\begin{align}
\lambda_k^{n+1}(r)=q\lambda_{k+1}^n(r)+\overline{q}\lambda_{k-1}^n(r)+r^nt^2.   \label{eq:relation-lambda}
\end{align}

The following identities will be frequently used:
\begin{align}
\lambda_2^1(r)&=q(\lambda_3^0(r)-r)+1,  \label{eq:lambda-1}   \\
\lambda_2^2(r)&=q\lambda_3^1(r)+\lambda_2^0(r)+rt^2+\overline{q}(t^2-1).   \label{eq:lambda-2}
\end{align}

\begin{nota}\label{nota:simplify}
\rm Let $\mathcal{R}_\diamond=R[t,r_1,r_2,r_3,r_4]$.

Let $\mathcal{X}$ be the set of products $g_1\cdots g_m$ such that either $g_i\in\{p_{14},p_{124},p_{568}\}$ for all $i$, or $g_i\in\{p_{36},p_{346},p_{782}\}$ for all $i$. Put $|g_1\cdots g_m|^\ast=m$.

We will write elements of $\mathcal{S}_\diamond(Y)$ as $\mathcal{R}_\diamond$-linear combinations of elements of $\mathcal{X}$, and mostly be concerned with relations of the form
$\xi u\equiv \eta v$ for $\xi,\eta\in\mathcal{R}_\diamond$, $u,v\in\mathcal{X}$.
Here $\xi u\equiv \eta v$ means $\xi u-\eta v=\sum_i\kappa_iw_i$ in $\mathcal{S}_\diamond(Y)$ for some $\kappa_i\in\mathcal{R}_\diamond$ and $w_i\in\mathcal{X}$ with $|w_i|^\ast<\min\{|u|^\ast,|v|^\ast\}$.
\end{nota}

\begin{figure}[h]
  \centering
  \includegraphics[width=11cm]{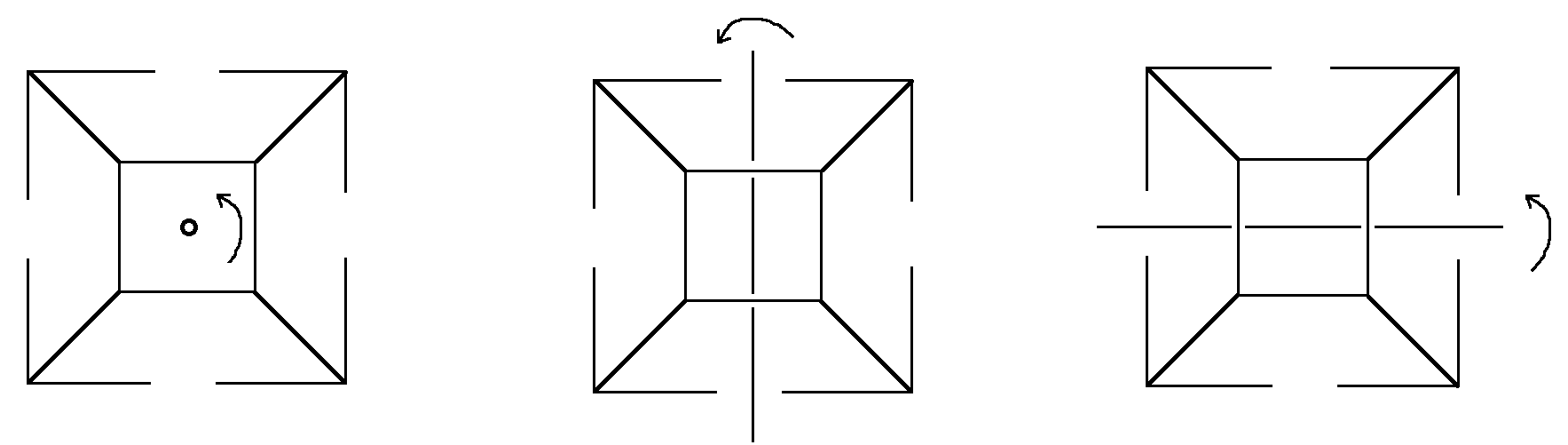}\\
  \caption{Left: $(\pi/2)$-rotation along the axis perpendicular to the plane. Middle: $\pi$-rotation along the vertical axis. 
  Right: $\pi$-rotation along the horizontal axis.}\label{fig:SP}
\end{figure}

To further simplify the computation, we adopt the {\it symmetry principle} (SP for short), in several perspectives, as provided by the symmetries of $U$ shown in Figure \ref{fig:SP}. For instance, The four tangles are the same, so each assertion holding for one tangle also holds for another.

In particular, the introduction of $p'_{134}$ is guided by SP; it is easy to see
$$p'_{134}=p_{124}+\overline{q}t(r_1-r_2).$$

For the main purpose, we can ignore ${\rm L}_2$-, ${\rm L}_4$-, ${\rm G}_1$-relations, and also those given by family (E).
This will be clear in Section \ref{sec:conclusion}.
So we only investigate ${\rm L}_1$-, ${\rm L}_3$- and ${\rm G}_2$-relations.



\subsection{${\rm L}_1$- and ${\rm L}_3$-relations}

\begin{figure}[h]
  \centering
  \includegraphics[width=11.5cm]{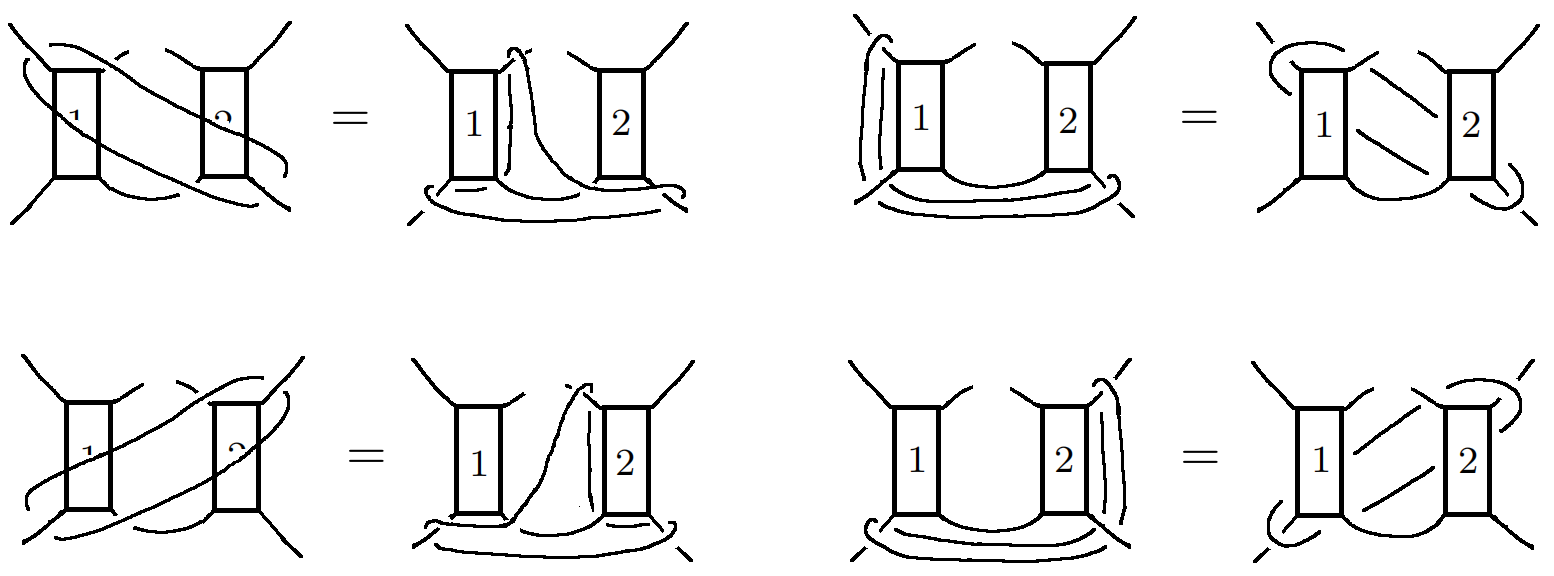}\\
  \caption{First row: $\grave{a}_{12}$, $\grave{b}_{12}$. Second row: $\acute{a}_{12}$, $\acute{b}_{12}$.}
  \label{fig:b=a-1}
\end{figure}

Introduced in Figure \ref{fig:b=a-1} are $\grave{a}_{12}$, $\grave{b}_{12}$, $\acute{a}_{12}$, $\acute{b}_{12}$.

For all $n\ge 0$ and $j\in\{1,2,3,4\}$, we have
\begin{align*}
r_j^ns_j=\lambda_3^n(r_j)+\lambda_2^n(r_j).
\end{align*}


\begin{lem}\label{lem:free}
Given $\vartheta_0,\vartheta_1,\vartheta_2\in R[t,r_1]$, if $\sum_{i=0}^2\vartheta_ir_2^i=0$ is a ${\rm L}_1$-relation, then $\vartheta_i=0$ for all $i$.
\end{lem}

See Appendix A for the proof. Consequently, for any $\xi,\zeta,\zeta'\in R[t,r_1,r_2]$, whenever $\xi x_1=\zeta$ and $\xi x_1=\zeta'$ are 
${\rm L}_1$-relations, actually $\zeta=\zeta'$, so we can just write $\xi x_1\equiv 0$, ignoring $\zeta$. This will benefit us a lot, since the expression of the ``constant" $\zeta$ is usually complicated.

Use $\xi\sim\xi'$ to denote $\xi x_1\equiv \xi'x_1$ (so that $\xi x_1^k\equiv \xi'x_1^k$ for each $k\ge 2$).

\begin{lem}\label{lem:L1}
The ${\rm L}_1$-relations are equivalent to the following:
\begin{enumerate}
  \item $c_2^0(r_j)t\ell_1\equiv -q\eta_3^n(r_j)x_1$, $\eta_3^n(r_j)\ell_1\equiv qc_3^n(r_j)tx_1$ for $j=1,2$.
  \item For $k\ge 3$, there exits $\varphi_k\in R[t][r]$ whose highest-order term is $\beta_kr^k$ for some $0\ne \beta_k\in R$
        such that $\varphi_k(r_j)\sim 0$ for $j=1,2$; in particular,
        $$\varphi_3=\lambda_3^0(r)-r, \ \ \ \varphi_4=\lambda_3^1(r)-r^2-\overline{q}^2+1,  \ \ \   \varphi_5=\lambda_3^2(r)-r^3-(\overline{q}^4-1)r.$$
  \item $r_1\sim r_2$, $r_1^2\sim r_2^2$, $r_1r_2\sim r_2^2+\overline{q}^2-1$, and
        \begin{align*}
        r_1(r_2^2-1)\sim (r_1^2-1)r_2&\sim \overline{q}(1-t^2)r_1r_2-\overline{q}^2t^2r_2-\overline{q}^3, \\ 
        (r_1^2-1)(r_2^2-1)&\sim-\overline{q}^6t^2c_2^0(r_1)c_3^0(r_2).
        \end{align*}
\end{enumerate}
\end{lem}

In parallel, we have
\begin{lem}\label{lem:L3}
The ${\rm L}_3$-relations are equivalent to the following:
\begin{enumerate}
  \item $c_2^0(r_j)t\ell_3\equiv -q\eta_3^n(r_j)x_1$, $\eta_3^n(r_j)\ell_3\equiv qc_3^n(r_j)tx_1$ for $j=3,4$.
  \item For $k\ge 3$, there exits $\varphi_k\in R[t][r]$ whose highest-order term is $\beta_kr^k$ for some $0\ne \beta_k\in R$
        such that $\varphi_k(r_j)\sim 0$ for $j=3,4$; in particular,
        $$\varphi_3=\lambda_3^0(r)-r, \ \ \ \varphi_4=\lambda_3^1(r)-r^2-\overline{q}^2+1,  \ \ \   \varphi_5=\lambda_3^2(r)-r^3-(\overline{q}^4-1)r.$$
  \item $r_4\sim r_3$, $r_4^2\sim r_3^2$, $r_3r_4\sim r_3^2+\overline{q}^2-1$, and
        \begin{align*}
        r_3(r_4^2-1)\sim (r_4^2-1)r_3&\sim \overline{q}(1-t^2)r_3r_4-\overline{q}^2t^2r_3-\overline{q}^3, \\ 
        (r_3^2-1)(r_4^2-1)&\sim-\overline{q}^6t^2c_2^0(r_4)c_3^0(r_3).
        \end{align*}
\end{enumerate}
\end{lem}

Let $\widetilde{U}=U\cup H_1\cup H_3$. The following is straightforward:
\begin{cor}\label{cor:module}
The $R[t]$-module $\mathcal{S}_\diamond(\widetilde{U})$ is generated by
\begin{align*}
&\{r_2^mr_3^nx_1^k\colon k\ge 1,\ 0\le m,n\le 2\}\cup\{r_1^{i_1}r_2^{i_2}r_3^{i_3}r_4^{i_4}\colon \ 0\le i_1,i_4\le 2\} \\
\cup\ &\{r_3^n\ell_1x_1^k,\ r_2r_3^n\ell_1x_1^k\colon k\ge 1, \ 0\le n\le 2\}  \\
&\cup\{r_3^nr_4^h\ell_1,\ r_2r_3^nr_4^h\ell_1\colon h\ge0, \ 0\le n\le 2\}  \\
\cup\ &\{r_2^m\ell_3x_1^k,\ r_2^mr_3\ell_3x_1^k\colon k\ge 1, \ 0\le m\le 2\}  \\
&\cup\{r_1^hr_2^m\ell_3,\ r_1^hr_2^mr_3\ell_3\colon h\ge 0,\ 0\le m\le 2\}  \\
\cup\ &\{r_2^mr_3^n\ell_1\ell_3x_1^k\colon k\ge 0,\ 0\le m,n\le 1\},
\end{align*}
subject to the following relations:
\begin{align*}
t\cdot r_2r_3^n\ell_1u&\equiv -\overline{q}tr_3^n\ell_1u-q^2(r_2^2-1)r_3^nux_1, \qquad u=x_1^k, r_4^h,   \\
t\cdot r_2^mr_3\ell_3u&\equiv -\overline{q}tr_2^m\ell_3u-q^2r_2^m(r_3^2-1)ux_1, \qquad u=x_1^k, r_1^h,   \\
t\cdot r_3\ell_1\ell_3x_1^k&\equiv -\overline{q}t\ell_1\ell_3x_1^k-q^2(r_3^2-1)\ell_1x_1^{k+1},  \\
t\cdot r_2r_3^n\ell_1\ell_3x_1^k&\equiv -\overline{q}tr_3^n\ell_1\ell_3x_1^k-q^2(r_2^2-1)r_3^n\ell_3x_1^{k+1}.
\end{align*}
\end{cor}

\subsection{${\rm G}_2$-relations}

In this subsection, we work in $\widetilde{U}$.
We shall investigate ${\rm Sl}^2_2(r_2^mr_3^nux_1^{k-1})$ and ${\rm Sl}^2_4(r_2^mr_3^nux_1^{k-2})$
for $k\ge 2$ and $u\in\{x_1^2,\ell_1x_1,\ell_3x_1,\ell_1\ell_3\}$, considering them as in ${\rm Sl}(\widetilde{U})$.

By Remark \ref{rmk:L1}, $r_1^n\grave{a}_{12}\equiv q\lambda_3^n(r_1)x_1$, $r_2^n\acute{a}_{12}\equiv \overline{q}\lambda_2^n(r_2)x_1$.
By SP, we have $r_3^n\grave{a}_{34}\equiv q\lambda_3^n(r_3)x_1$, $r_4^n\acute{a}_{34}\equiv \overline{q}\lambda_2^n(r_4)x_1$.
Also by SP, $r_3^n\grave{b}_{34}\equiv\overline{q}\lambda_2^n(r_2)x_1$.

\begin{figure}[H]
  \centering
  \includegraphics[width=12cm]{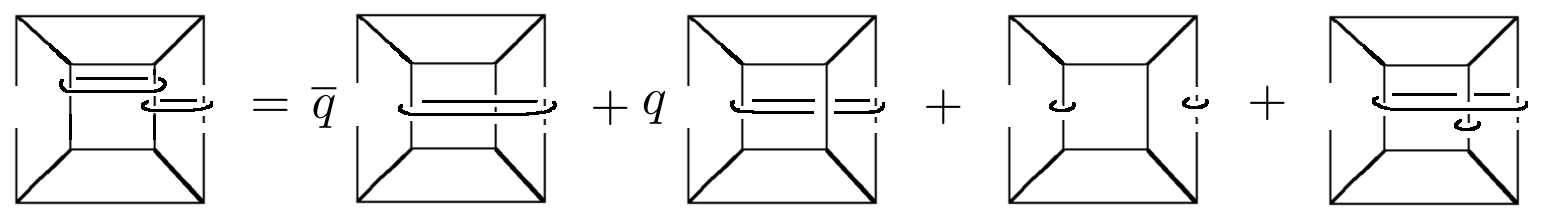}\\
  \caption{}\label{fig:back}
\end{figure}

From Figure \ref{fig:back} (which is similar to the last row of Figure \ref{fig:reduce-0}) we see that
${\rm Sl}^2_2(x_1^{k+1})$ can be replaced by $x_1^k(\acute{b}_{12}-\grave{b}_{34})$.

Given $\xi,\xi'\in \mathcal{R}_\diamond$, write $\xi\sim \xi'$ if $\xi x_1^{k+1}\equiv \xi'x_1^{k+1}$.

By Lemma \ref{lem:L1}, Lemma \ref{lem:L3} and (\ref{eq:lambda-1}), (\ref{eq:lambda-2}), for $j=1,2,3,4$, we have
\begin{align}
\lambda_2^1(r_j)&\sim 1,  \label{eq:equiv-1}  \\
\lambda_2^2(r_j)
&\sim(1-\overline{q}^2)\lambda_2^0(r_j)+r_j.  \label{eq:equiv-2}
\end{align}

\begin{lem}\label{lem:G2-redundant}
Let $\mathcal{E}_k\subset\mathcal{S}(\widetilde{U})$ be the $R$-submodule generated by ${\rm Sl}^2_2(r_2^mr_3^nx_1^{h})$,
${\rm Sl}^2_{4}(r_2^mr_3^nx_1^{h-1})$ for $0\le m,n\le 2$ and $h\le k$. Then
\begin{align*}
&R\{r_2^mr_3^nx_1^k(\acute{b}_{12}-\grave{b}_{34})\colon 0\le m,n\le 2\}  \\
=\ &R\{r_2^mr_3^nx_1^k(\acute{a}_{12}-\grave{a}_{34})\colon 0\le m,n\le 2\}+\mathcal{E}_k.
\end{align*}
\end{lem}

\begin{figure}[H]
  \centering
  \includegraphics[width=13cm]{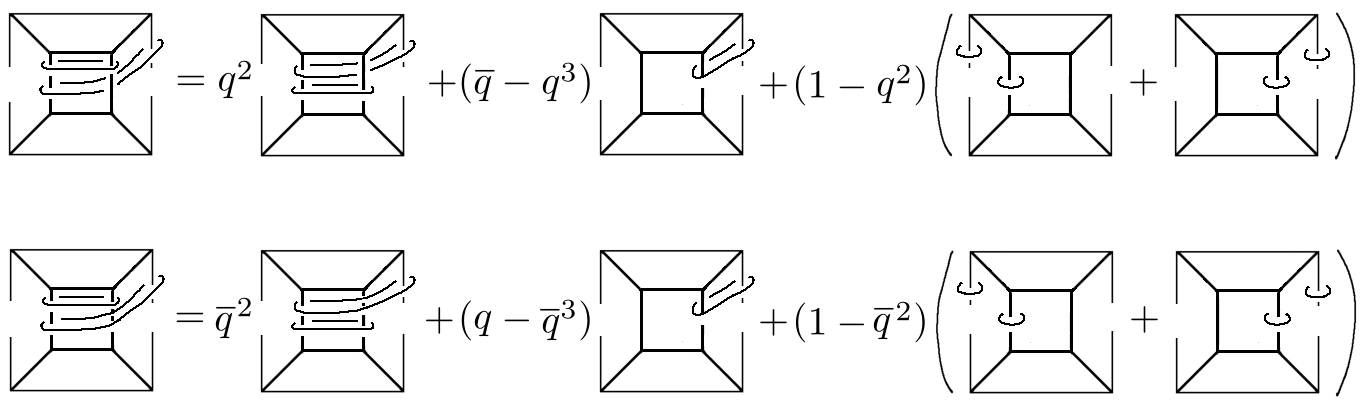}\\
  \caption{}\label{fig:commute-1}
\end{figure}

\begin{proof}
From Figure \ref{fig:commute-1} we see $x_1\acute{b}_{12}\equiv q^2\acute{b}_{12}x_1$ and $x_1\acute{a}_{12}\equiv\overline{q}^2\acute{a}_{12}x_1$.
So
\begin{align}
x_1^k(\acute{b}_{12}-\grave{b}_{34})&\equiv q^{2k-2}x_1\acute{b}_{12}x_1^{k-1}-x_1\grave{b}_{34}x_1^{k-1}
=x_1(q^{2k-2}\acute{b}_{12}-\grave{b}_{34})x_1^{k-1},  \nonumber  \\
x_1^k(\grave{a}_{34}-\acute{a}_{12})&\equiv (\grave{a}_{34}-\overline{q}^{2k}\acute{a}_{12})x_1^{k}
\equiv x_3(\overline{q}^{2}\grave{a}_{34}-\overline{q}^{2k}\acute{a}_{12})x_1^{k-1}.  \label{eq:equivalence-0}
\end{align}

From the expression of $\lambda_2^0(r_2)$ and (\ref{eq:equiv-1}), (\ref{eq:equiv-2})
we see
\begin{align}
&R\{r_2^mx_1(\acute{b}_{12}-\overline{q}^{2k-2}\grave{b}_{34})x_1^{k-1}\colon 0\le m\le 2\} \nonumber  \\
=\ &R\{r_2^i\acute{a}_{12}(\acute{b}_{12}-\overline{q}^{2k-2}\grave{b}_{34})x_1^{k-1}\colon 0\le i\le 2\}+\mathcal{E}_k.  \label{eq:equivalence-1}
\end{align}
Similarly,
\begin{align}
&R\{r_3^nx_3(\grave{a}_{34}-\overline{q}^{2k-2}\acute{a}_{12})x_1^{k-1}\colon 0\le n\le 2\}  \nonumber \\
=\ &R\{r_3^j\grave{b}_{34}(\grave{a}_{34}-\overline{q}^{2k-2}\acute{a}_{12})x_1^{k-1}\colon 0\le j\le 2\}+\mathcal{E}_k. \label{eq:equivalence-2}
\end{align}

By the formula in Figure \ref{fig:y1-times-y2},
\begin{align*}
\acute{a}_{12}\acute{b}_{12}=\acute{a}_{12}\grave{a}_{12}&=q^2x_1^2+2qt^2x_1+\overline{q}^2s_2^2+2\overline{q}t^2s_2+t^4+4t^2+\alpha^2,  \\
\grave{b}_{34}\grave{a}_{34}=\acute{a}_{34}\grave{a}_{34}&=q^2x_3^2+2qt^2x_3+\overline{q}^2s_3^2+2\overline{q}t^2s_3+t^4+4t^2+\alpha^2.
\end{align*}
Hence
\begin{align*}
(\acute{a}_{12}\acute{b}_{12}-\grave{b}_{34}\grave{a}_{34})x_1^{k-1}
=(\overline{q}^2(s_2^2-s_3^2)+2\overline{q}t^2(s_2-s_3))x_1^{k-1}\equiv 0.
\end{align*}
Noticing $\grave{b}_{34}\acute{a}_{12}=\acute{a}_{12}\grave{b}_{34}$ 
and comparing (\ref{eq:equivalence-1}) with (\ref{eq:equivalence-2}) completes the proof.
\end{proof}

\begin{rmk}
\rm Consequently, ${\rm Sl}^2_2(x_1^{k+1})$ is redundant.
Similar arguments can show that ${\rm Sl}^2_2(\ell_1x_1^k)$, ${\rm Sl}^2_2(\ell_3x_1^k)$, ${\rm Sl}^2_2(\ell_1\ell_3x_1^{k-1})$ are all redundant.
\end{rmk}

In virtue of (\ref{eq:equivalence-0}), (\ref{eq:equivalence-2}), the relations ${\rm Sl}^2_4(r_2^ir_2^jx_1^k)$ are equivalent to $$r_2^mr_3^n\grave{b}_{34}(\acute{a}_{12}-q^{2k-2}\grave{a}_{34})x_1^{k-1}\equiv 0, \qquad 0\le m,n\le 2;$$
since $r_3^n\grave{b}_{34}\equiv\overline{q}\lambda_2^n(r_3)x_1$, $r_2^m\acute{a}_{12}\equiv\overline{q}\lambda_2^m(r_2)x_1$, and $\grave{b}_{34}\grave{a}_{34}\equiv q^2x_1^2$, they have the form
$$\lambda_2^m(r_2)\lambda_2^n(r_3)\sim q^{2k+2}r_2^mr_3^n, \qquad  0\le m,n\le 2,$$
into which we shall take a deep dive.
\begin{enumerate}
  \item The case $m=n=0$ is simply
        $\lambda_2^0(r_2)\lambda_2^0(r_3)\sim q^{2k+2}.$
  \item The case $m=0$, $n=1$ reads $\lambda_2^0(r_2)\lambda_2^1(r_3)\sim q^{2k+2}r_3$, which together with (\ref{eq:equiv-1}) for $j=3$
        yields $\lambda_2^0(r_2)\sim q^{2k+2}r_3$, or equivalently,
        \begin{align*}
        \gamma_3(r_2)\sim\overline{q}(1-t^2)r_2-q^{2k-1}r_3-\overline{q}^2t^2.
        \end{align*}
  \item In parallel, $m=1$, $n=0$ leads to 
        \begin{align*}
        \gamma_3(r_3)\sim\overline{q}(1-t^2)r_3-q^{2k-1}r_2-\overline{q}^2t^2.
        \end{align*}
  \item $m=1,n=1$ leads to $q^{2k+2}r_2r_3\sim\lambda_2^1(r_2)\lambda_2^1(r_3)\sim 1$, so
        $r_2r_3\sim \overline{q}^{2k+2}.$
  \item $m=0,n=2$ leads to
        \begin{align*}
        q^{2k+2}r_3^2&\sim\lambda_2^0(r_2)\lambda_2^2(r_3)\stackrel{(\ref{eq:equiv-2})}\sim \lambda_2^0(r_2)\big((1-\overline{q}^2)\lambda_2^0(r_3)+r_3\big)  \\
        &\sim (1-\overline{q}^2)q^{2k+2}+\lambda_2^0(r_2)r_3,
        \end{align*}
        implying
        $\lambda_2^0(r_2)r_3\sim q^{2k+2}\gamma_3(r_3)+q^{2k}.$
  \item In parallel, $m=2,n=0$ leads to
        $r_2\lambda_2^0(r_3)\sim q^{2k+2}\gamma_3(r_2)+q^{2k}.$
  \item $m=1,n=2$ leads to
        $$q^{2k+2}r_2r_3^2\sim\lambda_2^2(r_3)\stackrel{(\ref{eq:equiv-2})}\sim(1-\overline{q}^2)\lambda_2^0(r_3)+r_3
        \sim(1-\overline{q}^2)q^{2k+2}r_2+r_3,$$
        implying
        $r_2\gamma_3(r_3)\sim\overline{q}^{2k+2}r_3-\overline{q}^2r_2.$
  \item In parallel, $m=2,n=1$ leads to
        $\gamma_3(r_2)r_3\sim\overline{q}^{2k+2}r_2-\overline{q}^2r_3.$
  \item Finally, $m=2,n=2$ leads to
        \begin{align*}
        &q^{2k+2}r_2^2r_3^2\sim\lambda_2^2(r_2)\lambda_2^2(r_3)
        \sim((1-\overline{q}^2)\lambda_2^0(r_2)+r_2)((1-\overline{q}^2)\lambda_2^0(r_3)+r_3)  \\
        \sim\ &(1-\overline{q}^2)^2q^{2k+2}+(1-\overline{q}^2)\big(q^{2k+2}(\gamma_3(r_2)+\gamma_3(r_3))+2q^{2k}\big)+\overline{q}^{2k+2},
        \end{align*}
        implying
        $(\gamma_3(r_2)+\overline{q}^2)(\gamma_3(r_3)+\overline{q}^2)\sim \overline{q}^{4k+4}.$
\end{enumerate}

\begin{nota}
\rm Let $\alpha_k=q^k+\overline{q}^k$. Sometimes denote $\theta$ for $t^2-1$.

For $1\le i\le 9$, use ``Eq.$i$" to denote the last equation obtained in Case $i$.
\end{nota}

Examining Eq.6 and using Eq.2, Eq.4, Eq.7, we obtain
\begin{align*}
q^{2k}&\sim r_2\lambda_2^0(r_3)-q^{2k+2}\gamma_3(r_2)
=r_2(-q^3\gamma_3(r_3)-q^2\theta r_3-qt^2)-q^{2k+2}\gamma_3(r_2)   \\
&\sim q^3(\overline{q}^2r_2-\overline{q}^{2k+2}r_3)-\overline{q}^{2k}\theta-qt^2r_2+q^{2k+2}(\overline{q}\theta r_2+q^{2k-1}r_3+\overline{q}^2t^2)  \\
&\sim q(q^{2k}-1)(\theta r_2+(\alpha_k^2-1)r_3+\overline{q}^{k+1}\alpha_k\theta)+q^{2k}.
\end{align*}
Hence
\begin{align}
\rho_2:=(t^2-1)r_2+(\alpha_k^2-1)r_3+\overline{q}^{k+1}\alpha_k(t^2-1)\sim 0.  \label{eq:sim-7'}
\end{align}
In a parallel way, we have
\begin{align}
\rho_3:=(t^2-1)r_3+(\alpha_k^2-1)r_2+\overline{q}^{k+1}\alpha_k(t^2-1)\sim 0.   \label{eq:sim-8'}
\end{align}
Examining Eq.1 and using Eq.2--Eq.4, Eq.7--Eq.9, we have
\begin{align*}
q^{2k+2}&\sim \lambda_2^0(r_2)\lambda_2^0(r_3)=(q^3\gamma_3(r_2)+q^2\theta r_2+qt^2)(q^3\gamma_3(r_3)+q^2\theta r_3+qt^2)   \\
&=q^6(\gamma_3(r_2)+\overline{q}^2)(\gamma_3(r_3)+\overline{q}^2)+q^5\theta(r_2\gamma_3(r_3)+\gamma_3(r_2)r_3)  \\
&\ \ \ +q^4\theta(\gamma_3(r_2)+\gamma_3(r_3))+q^4\theta^2r_2r_3+q^3t^2\theta(r_2+r_3)+q^2t^4-q^2  \\
&\sim q^{2-4k}+q^3\theta(\overline{q}^{2k}-1)(r_2+r_3)-q^4\theta\big(\overline{q}(\theta+q^{2k})(r_2+r_3)+2\overline{q}^2t^2\big) \\
&\ \ \ +q^{2-2k}\theta^2+q^3t^2\theta(r_2+r_3)+q^2t^4-q^2  \\
&=q^2(\overline{q}^{2k}-1)(q^{k+1}\alpha_k(t^2-1)(r_2+r_3)+t^4-2t^2+\alpha_k^2)+q^{2k+2}.
\end{align*}
Hence
$$q^{k+1}\alpha_k(t^2-1)(r_2+r_3)+t^4-2t^2+\alpha_k^2\sim 0;$$
on the other hand, summing (\ref{eq:sim-7'}) and (\ref{eq:sim-8'}) yields
$$(t^2+\alpha^2_k-2)(r_2+r_3)+2\overline{q}^{k+1}\alpha_k(t^2-1)\sim 0.$$
Eliminating $r_2+r_3$, we obtain
\begin{align}
\rho_1:=(t^2-\alpha^2_k)(t^2-2+\alpha_k)(t^2-2-\alpha_k)\sim 0.  \label{eq:sim-9'}
\end{align}

Eq.2, Eq.3, Eq.4, Eq.7, Eq.8, Eq.9 can be rewritten respectively as
\begin{align*}
\rho_4:&=r_2^2+\overline{q}(t^2-1)r_2+q^{2k-1}r_3+\overline{q}^2t^2-1\sim 0,  \\
\rho_5:&=r_3^2+\overline{q}(t^2-1)r_3+q^{2k-1}r_2+\overline{q}^2t^2-1\sim 0,  \\
\rho_6:&=r_2r_3-\overline{q}^{2k+2}\sim 0,   \\
\rho_7:&=r_2r_3^2-\overline{q}^{2k+2}r_3+(\overline{q}^2-1)r_2\sim 0,  \\
\rho_8:&=r_2^2r_3-\overline{q}^{2k+2}r_2+(\overline{q}^2-1)r_3\sim 0,  \\
\rho_9:&=(r_2^2-1+\overline{q}^2)(r_3^2-1+\overline{q}^2)-\overline{q}^{4k+4}\sim 0.
\end{align*}

\begin{rmk}
\rm We can regard (\ref{eq:sim-7'})--(\ref{eq:sim-9'}) as ``quantum effects", since they are deduced from dividing by $q^{2k}-1$.

It is worth seeking the topological meaning of $\rho_1$.
\end{rmk}

\medskip

\begin{figure}[h]
  \centering
  \includegraphics[width=13cm]{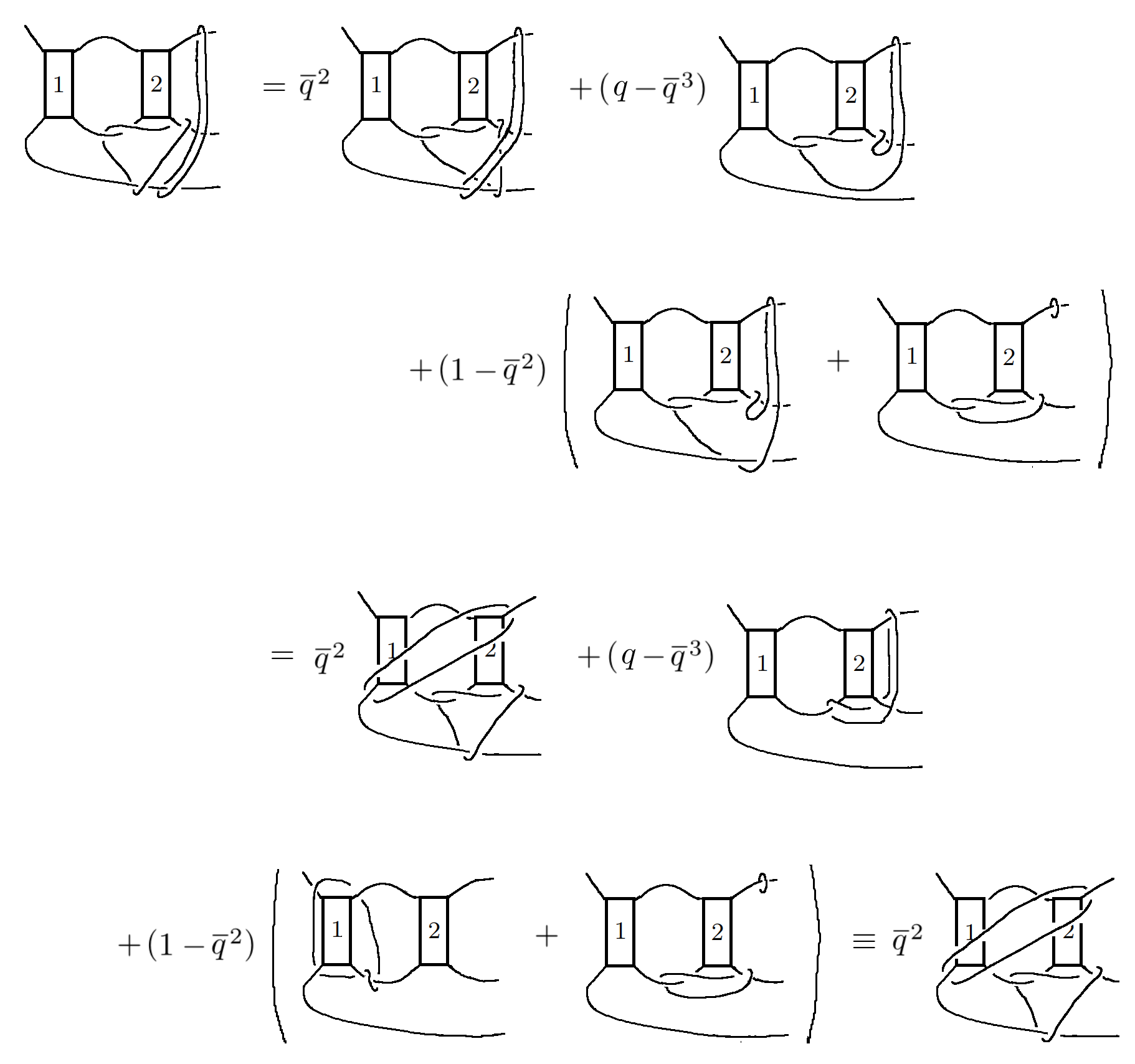}\\
  \caption{}\label{fig:commute-2}
\end{figure}

Now consider ${\rm Sl}^2_4(u)$ when $\ell_1\mid u$ or $\ell_3\mid u$.

For simplicity, we do not consider ${\rm Sl}^2_{4}(u)$ for $\ell_1\mid u$ or $\ell_3\mid u$ but $t\nmid u$.
Use $G'_2$-relator to name a $G_2$-relator not of such form.

It is important to note that $t$ times a $G_2$-relator always equals a $R$-linear combination of $G'_2$-relators.

From Figure \ref{fig:commute-2} we see $p_{124}\acute{a}_{12}\equiv\overline{q}^2\acute{a}_{12}p_{124}$; on the other hand,
\begin{align}
\acute{a}_{12}tp_{124}&\equiv\overline{q}^2\big(q\eta_2^0(r_2)p_{14}+c_2^0(r_2)tp_{134}\big)tp_{124} \nonumber  \\
&\equiv \overline{q}^2\big(q\eta_2^0(r_2)tp_{14}p_{124}+c_2^0(r_2)t^2p_{14}^2\big)   \nonumber  \\
&\equiv\overline{q}^2(c_2^0(r_2)t^2+\eta_3^0(r_2))p_{14}-tp_{124}.   \label{eq:final-1}
\end{align}
Similarly, $p_{568}\grave{a}_{34}\equiv q^2\grave{a}_{34}p_{568}$, and (using $q\eta_4=\eta_3^1-\overline{q}c_2^0+1$)
\begin{align}
\grave{a}_{34}tp_{568}&\equiv \big(q^2c_3^0(r_3)tp_{568}-q\eta_4^0(r_3)p_{58}\big)tp_{568}  \nonumber   \\
&\equiv -c_3^0(r_3)t^2p_{58}^2-q(c_3^0(r_3)t^2+\eta_4^0(r_3))tp_{58}p_{568}   \nonumber  \\
&\equiv q\lambda_3^0(r_3)tp_{58}p_{568}\equiv \overline{q}^2\eta_3^0(r_3)p_{58}^2-tp_{58}p_{568}.  \label{eq:final-2}
\end{align}

\begin{figure}[h]
  \centering
  \includegraphics[width=11cm]{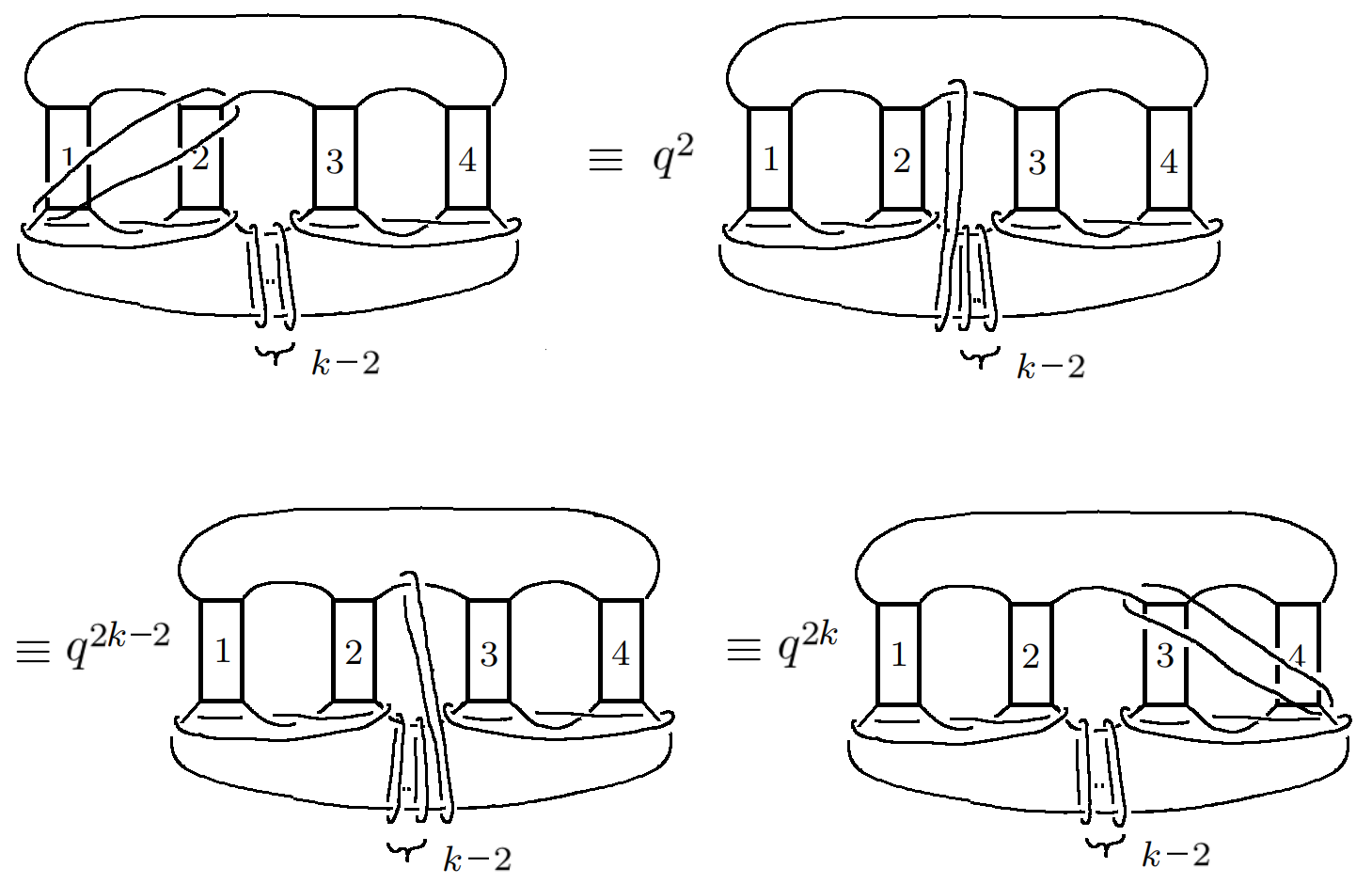}\\
  \caption{}\label{fig:commute-final}
\end{figure}

By Corollary \ref{cor:module}, it suffices to consider ${\rm Sl}^2_{4}(r_3^nt\ell_1x_1^k)=0$ for $0\le n\le 2$,
${\rm Sl}^2_{4}(r_2^mt\ell_3x_1^k)=0$ for $0\le m\le 2$, and ${\rm Sl}^2_{4}(t\ell_1\ell_3x_1^{k-2})=0$.
\begin{enumerate}
  \item ${\rm Sl}^2_4(r_3^nt\ell_1x_1^{k-1})=0$ takes the form
        $$r_3^n\acute{a}_{12}tp_{124}p_{14}^{k-1}\equiv q^{2k}r_3^ntp_{124}p_{14}^{k-1}\grave{a}_{34},$$
        which, by (\ref{eq:final-1}), is equivalent to
        \begin{align*}
        (q^{2k+1}\lambda_3^n(r_3)+r_3^n)tp_{14}^kp_{124}&\equiv\overline{q}^2(c_2^0(r_2)t^2+\eta_3^0(r_2))r_3^np_{14}^{k+1}.
        \end{align*}
        Denote $\xi\simeq\xi'$ if $(\xi-\xi')tp_{124}p_{14}^k\equiv \kappa p_{14}^{k+1}$ for some $\kappa\in\mathcal{R}_\diamond$.

        The cases $n=0$ and $n=1$ respectively lead to
        \begin{align*}
        r_3+\overline{q}^{2k+1}\simeq 0,  \qquad       r_3^2+\overline{q}^2-1+\overline{q}^{2k+1}r_3\simeq 0,
        \end{align*}
        and then $n=2$ gives rise to
        \begin{align*}
        0&\simeq (\overline{q}^{2k+1}-\overline{q}\theta)r_3^2+(1-\overline{q}^2t^2)r_3+(\overline{q}-\overline{q}^3)t^2-\overline{q}  \\
        &\simeq (1-\overline{q}^{2k})\overline{q}^{2k+3}(t^2-\alpha_k^2)
        \end{align*}
        Thus, ${\rm Sl}^2_4(r_3^ntx_1^{k-1}\ell_1)=0$ for $0\le n\le2$ exactly reduce $r_3^2tp_{124}p_{14}^k$, $r_3tp_{124}p_{14}^k$ and $t^3p_{124}p_{14}^k$.
  \item ${\rm Sl}^2_4(r_2^mt\ell_3x_3^{k-1})=0$ takes the form
        $$r_2^mtp_{568}p_{58}^{k-1}\acute{a}_{12}\equiv q^{2k}r_2^m\grave{a}_{34}tp_{568}p_{58}^{k-1},$$
        which, by (\ref{eq:final-2}), is equivalent to
        $$(\overline{q}^{2k+1}\lambda_2^m(r_2)+r_2^m)tp_{14}^kp_{568}\equiv \overline{q}^2r_2^m\eta_3^0(r_3)p_{58}^{k+1}.$$
        Similarly as above, ${\rm Sl}^2_4(r_2^mtx_3^{k-1}\ell_3)=0$ for $0\le m\le2$ exactly reduce $r_2^2tp_{568}p_{58}^k$, $r_2tp_{568}p_{58}^k$ and $t^3p_{568}p_{58}^k$.
  \item 
        From Figure \ref{fig:commute-final} we see that ${\rm Sl}^2_4(t\ell_1\ell_3x_1^{k-2})=0$ takes the form
        $$\acute{a}_{12}tp_{14}^{k-2}p_{124}p_{568}\equiv q^{2k}\grave{a}_{34}tp_{14}^{k-2}p_{124}p_{568},$$
        which, by (\ref{eq:final-1}), (\ref{eq:final-2}), reads
        \begin{align*}
        &\big(\overline{q}^2(c_2^0(r_2)t^2+\eta_3^0(r_2))p_{14}^k-tp_{14}^{k-1}p_{124}\big)p_{568}  \\
        \equiv\ &q^{2k}(\overline{q}^2\eta_3^0(r_3)p_{14}^k-tp_{14}^{k-1}p_{568})p_{124}.
        \end{align*}
        Hence $tp_{14}^{k-1}p_{124}p_{568}\equiv \kappa_1p_{14}^kp_{568}+\kappa_2p_{14}^kp_{124}$ for some $\kappa_1,\kappa_2\in\mathcal{R}_\diamond$.
\end{enumerate}

\subsection{Conclusion}\label{sec:conclusion}

Decompose the free $\mathcal{R}_\diamond$-module $\mathcal{S}_\diamond(U)$ as $\mathcal{M}_0\oplus\mathcal{M}_1$, where
$$\mathcal{M}_0=\mathcal{R}_\diamond\{u\in\mathcal{X}\colon x_1^3,\ell_1x_1^2,\ell_3x_1^2,\ell_1\ell_3x_1\nmid u\},$$
(recalling Theorem \ref{thm:generator} for $\mathcal{X}$), and $\mathcal{M}_1$ is its direct summand.

Note that ${\rm L}_2$-, ${\rm L}_4$-, ${\rm G}_1$-relators and those in the (E) family all belong to $\mathcal{M}_0$.
Let $\mathcal{K}\subset\mathcal{M}_1$ (resp. $\mathcal{K}'$) denote the image of ${\rm L}_1$-, ${\rm L}_3$- and ${\rm G}_2$- (resp. ${\rm G}'_2$-) relators under the projection $\mathcal{S}_\diamond(U)\twoheadrightarrow\mathcal{M}_1$. Note that $t\cdot\mathcal{K}\subseteq\mathcal{K}'$.

The result of the previous subsections is summarized as
\begin{lem}
As a $R$-module, $\mathcal{M}_1/\mathcal{K}'\cong \bigoplus_{k\ge 2}R\mathcal{B}_k$, with
\begin{align*}
\mathcal{B}_k=\ &\{t^ix_1^{k+1}\colon 0\le i\le 5\}\cup\{t^ir_2x_1^{k+1}\colon 0\le i\le 3\}  \\
&\cup\{r_3^n\ell_1x_1^{k},r_2r_3^n\ell_1x_1^{k}\colon 0\le n\le 2\}\cup\{t\ell_1x_1^{k},t^2\ell_1x_1^{k}\}  \\
&\cup\{r_2^n\ell_3x_1^{k},r_2^nr_3\ell_3x_1^{k}\colon 0\le n\le 2\}\cup\{t\ell_3x_1^{k},t^2\ell_3x_1^{k}\}  \\
&\cup\{r_2^mr_3^n\ell_1\ell_3 x_1^{k-1}\colon 0\le m,n\le 1\}.
\end{align*}
Moreover, $t\cdot\mathcal{B}_j\subset \bigoplus_{k=2}^jR\mathcal{B}_k$ for each $j$.
\end{lem}

\begin{proof}
For the coefficients of $x_1^{k+1}$, just note that as a $R$-module,
\begin{align*}
&R[t]\{r_2^mr_3^n\colon 0\le m,n\le 2\}/R[t]\{\rho_1,\ldots,\rho_9\} \\
\cong\ &R\{t_i\colon 0\le i\le 5\}\oplus R\{t^ir_2\colon 0\le i\le 3\}.
\end{align*}
Indeed, if $\rho_2\sim 0$ is used to substitute $r_3$, then $\rho_3\sim 0$ can be replaced by 
$$\big(\theta^2-(q^{2k}+1+\overline{q}^{2k})^2\big)r_2\sim \overline{q}^{k+1}\alpha_k\theta(q^{2k}+1+\overline{q}^{2k}-\theta).$$

The remaining thing is easy.
\end{proof}

Let $\lambda$ denote the following composite of $R[t]$-module morphisms
$$\mathcal{S}_\diamond(Y)\cong\mathcal{S}_\diamond(U)/{\rm Sl}_\diamond(U)
\twoheadrightarrow\mathcal{M}_1/\mathcal{K}\stackrel{t\cdot}\to\mathcal{M}_1/\mathcal{K}'.$$
The most important point is that $\lambda(x_1^{h+1})=tx_1^{h+1}\in R\mathcal{B}_h$ does not vanish.

\begin{proof}[Proof of Theorem \ref{thm:main}]
Given any finite subset $\mathcal{N}\subset\mathcal{S}(Y)$, denote its image under the canonical quotient map $\mathcal{S}(Y)\twoheadrightarrow\mathcal{S}_\diamond(Y)$ by $\mathcal{N}_\diamond$.

Take $d\ge 2$ such that $\lambda(\mathcal{N}_\diamond)\subset\bigoplus_{k=2}^dR\mathcal{B}_k$.
Then $x_1^{h+1}\notin R[t]\mathcal{N}_\diamond$ for any $h>d$: otherwise we would obtain
$$\lambda(x_1^{h+1})\in \lambda(R[t]\mathcal{N}_\diamond)\subseteq{\bigoplus}_{k=2}^dR\mathcal{B}_k,$$
contradicting that $\lambda(x_1^{h+1})\in R\mathcal{B}_h$ is nonzero.
Thus, $x_1^{h+1}\notin R[t_1,t_2]\mathcal{N}$.
\end{proof}

\section*{Appendix A: Proof of Lemma \ref{lem:free}}

In this section, let $\mathcal{S}'(M)=\mathcal{S}(M;\mathbb{Z}[q^{\pm\frac{1}{2}}])$ for a $3$-manifold $M$.

\begin{figure}[H]
  \centering
  \includegraphics[width=8.5cm]{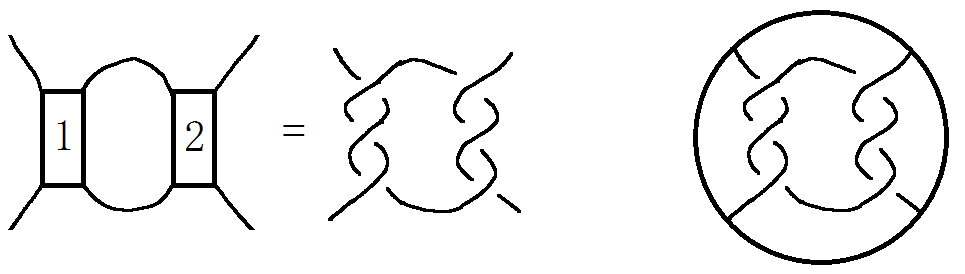}\\
  \caption{Left: the tangle $T$. Right: $T$ is put in a $3$-ball $B$ with $B\cap T=\partial T$.}\label{fig:T33}
\end{figure}

Let $T$ denote the tangle given in Figure \ref{fig:T33}. Take a $3$-ball $B$ containing $T$ with $B\cap T=\partial T$, and let $X$ be the complement of a tubular neighborhood of $T$ in $B$.
What we are really going to show is
\begin{lem}\label{lem:free'}
Given $\vartheta_0,\vartheta_1,\vartheta_2\in \mathbb{Z}[q^{\pm\frac{1}{2}}][t_1,t_2,r_1]$, if $\vartheta_i\ne 0$ for at least one $i$, then $\sum_{i=0}^2\vartheta_ir_2^i\ne0$ in $\mathcal{S}'_\diamond(X):=\mathcal{S}'(X)/(t_1-t_2)\mathcal{S}'(X)$.
\end{lem}

\begin{figure}[h]
  \centering
  \includegraphics[width=11cm]{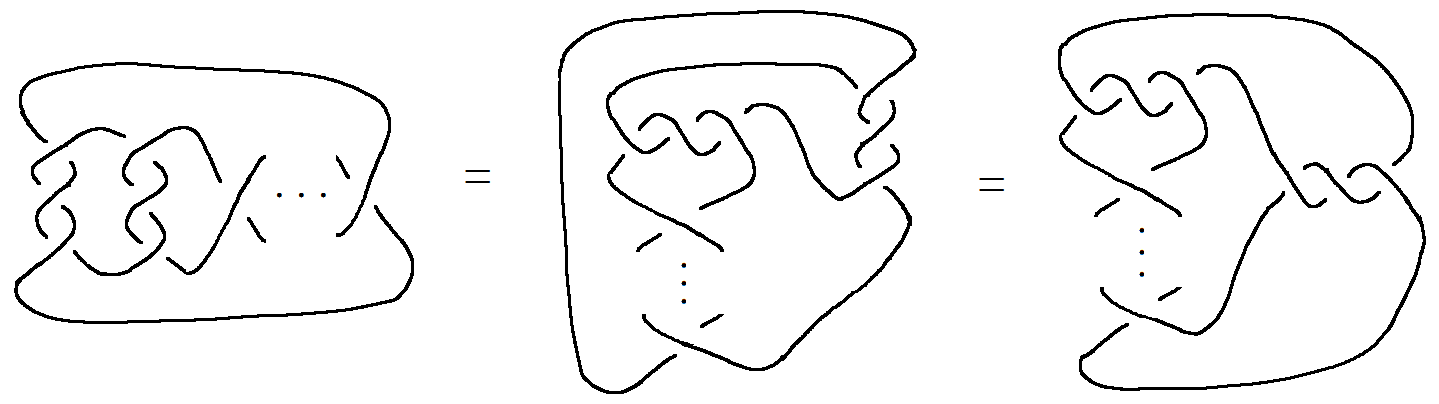}\\
  \caption{The tangle $T$ can be enlarged to the 2-bridge knot $\mathfrak{b}((9n+6)/(3n+1))$, by joining with $-n$ half-twists.}\label{fig:enlarge}
\end{figure}

\begin{figure}[h]
  \centering
  \includegraphics[width=4.5cm]{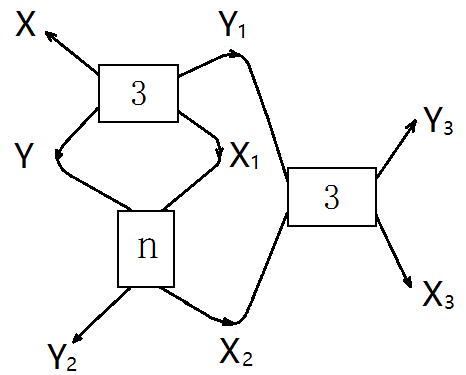}\\
  \caption{The rational tangle $[[3],[n],[3]]$. Connecting its northwest- and southwest end respectively to northeast- and southeast end yields $\mathfrak{b}([[3,n,3]])$.}\label{fig:rational}
\end{figure}

Take a sufficiently large odd integer $n$. Join $T$ with $-n$ half-twists and complete it to a knot, which is the 2-bridge knot $\mathfrak{b}([[3,n,3]])$, as shown in Figure \ref{fig:enlarge}. Let $N$ denote its exterior in $S^3$. The inclusion $X\subset N$ induces a morphism $\mathcal{S}'(X)\to\mathcal{S}'(N)$ which descends to $\mathcal{S}'_\diamond(X)\to\mathcal{S}'(N)$.

By \cite{Le06}, $\mathcal{S}'(N)$ is freely generated by $t^ir_1^j$ for $i\ge 0$ and $0\le j\le(9n+5)/2$;
in particular, $\mathcal{S}'(N)$ is torsion-free.

As elucidated in \cite[Page 4--5]{Ch18}, a trace-free ${\rm SL}(2,\mathbb{C})$-representation of a general 2-bridge link
$\mathfrak{b}([[k_1,\ldots,k_h]])$ can be constructed as follows.
For $a\in\mathbb{C}^\times$, let
$$A(a)=\frac{1}{2}\left(\begin{array}{cc}  (a+a^{-1})\sqrt{-1} &  a-a^{-1} \\  a-a^{-1} & -(a+a^{-1})\sqrt{-1} \end{array}\right).$$
Put $\mathbf{x}=A(1)$, $\mathbf{y}=A(b_0)$, and $\mathbf{y}_{i}=A(b_i)$ for $1\le i\le h$, where $b_0=b$, $b_1=-b^{-k_1}$ and recursively, $b_i=b_{i-1}^{-k_i}b_{i-2}$.

In the present case, $k_1=3$, $k_2=n$, $k_3=3$, so
$$b_1=-b^{-3}, \qquad  b_2=-b^{3n+1}, \qquad  b_3=b^{-9n-6}.$$
The condition for $b$ to define a representation of $\mathfrak{b}([[3,n,3]])$ is $\mathbf{x}^{-1}=\mathbf{y}_{3}=A(b_3)$,
which is equivalent to $(-b)^{9n+6}=1$.
So $b=-\zeta^k$ with $0\le k<9n+6$, where $\zeta=e^{2\pi\sqrt{-1}/(9n+6)}$. Let $\phi_j=\zeta^j+\zeta^{-j}$.

Note that ${\rm tr}(\mathbf{x}\mathbf{y}_{i})=-b_i-b_i^{-1}$. Hence
\begin{align*}
r_1=-{\rm tr}(\mathbf{x}\mathbf{y})=-\phi_k,  \qquad
r_2=-{\rm tr}(\mathbf{x}\mathbf{y}_{2})=-\phi_{k(3n+1)}.
\end{align*}

Assume $\sum_{i=0}^2\vartheta_ir_2^i=0$ in $\mathcal{S}'(N)$ for some $\vartheta_i\in \mathbb{Z}[q^{\pm\frac{1}{2}}][t,r_1]$.
Due to the freeness of $\mathcal{S}'(N)$, simultaneously dividing out a power of $t$ if necessary, we may assume that $t\nmid\vartheta_i$ for at least one $i$. Let $\sigma_i(r_1)=\vartheta_i(0,r_1)$.

Note that
\begin{align*}
{\sum}_{i=0}^{2}\sigma_ir_2^i=
(\sigma_0+2\sigma_2)-\sigma_1\phi_{k(3n+1)}+\sigma_2\phi_{k(3n+4)}=f_n(\phi_k),
\end{align*}
where
$$f_n(x)=\sigma_0(x)+2\sigma_2(x)-\sigma_1(x)\theta_{3n+1}(x)+\sigma_2\theta_{3n+4}(x).$$
Write $9n+5=2m$.
Observe that $\det\big((\phi_{kv})_{k,v=0}^m\big)\ne 0$, similarly as the non-vanishingness of the Vandermonde determinant.
It follows from
$$f_n(\phi_k)=0, \qquad  0\le k\le m$$
that $f_n(x)=0$.
When $n$ is sufficiently large, the coefficients in $\sigma_0,\sigma_1,\sigma_2$ do not cancel each other in $f_n(x)$,
contradicting the assumption.

Thus, actually $\sum_{i=0}^2\vartheta_ir_2^i\ne0$ in $\mathcal{S}'_\diamond(X)$.

\section*{Appendix B: Some useful formulas}

For $k,n\ge 0$ and $r=\mu+\mu^{-1}$, let $\lambda_{\pm}=q^{\pm2}\mu+q^{\mp2}\mu^{-1}$, and put
\begin{align}
\eta_k^n(r)=q^k\frac{\mu^{k}\lambda_+^n-\mu^{-k}\lambda_-^n}{\mu-\mu^{-1}}={\sum}_{j=0}^n\tbinom{n}{j}q^{2(2j-n)+k}\gamma_{2j-n+k}(r).
\end{align}

It is obvious that the highest-order term of $\eta_k^n(r)$ is
$$\begin{cases}
q^{2n+k}r^{n+k-1}, &k>0 \\
(q^{2n}-\overline{q}^{2n})r^{n-1},&k=0
\\ -\overline{q}^{2n+|k|}r^{n+|k|-1}, &k<0
\end{cases}.$$

The following identities are easy to verify:
\begin{align}
\eta_k^0(r)&=q^k\gamma_k(r), \nonumber \\
\eta_k^1(r)&=q^k(q^2\gamma_{k+1}(r)+\overline{q}^{2}\gamma_{k-1}(r)), \nonumber \\
r\eta_k^n(r)&=\overline{q}\eta_{k+1}^n(r)+q\eta_{k-1}^n(r), \label{eq:relation-eta-1} \\
\eta_k^{n+1}(r)&=q\eta_{k+1}^n(r)+\overline{q}\eta_{k-1}^n(r);  \label{eq:relation-eta-2}  \\
c_k^n(r)&=c_{k-1}^n(r)+\eta_k^n(r), \nonumber \\
rc_k^n(r)&=\overline{q}c_{k+1}^n(r)+qc_{k-1}^n(r)-r^n, \nonumber  \\
c_k^{n+1}(r)&=qc_{k+1}^n(r)+\overline{q}c_{k-1}^n(r)-r^n.  \label{eq:relation-c}
\end{align}

Here are the values of $c_k^n(r)$ for some small $k,n$:
\begin{align*}
c_0^0(r)=0, \qquad c_0^1(r)=q^{2}-1, \qquad  c_1^0(r)=q, \qquad c_1^1(r)=q^3r+q^2-1,  \\ 
c_2^0(r)=q^2r+q, \qquad c_2^1(r)=q^4(r^2-1)+q^3r+q^2.
\end{align*}
The following identities can be directly verified and will be useful:
\begin{align}
c_3^0(r)&=\overline{q}c_2^1(r),  \label{eq:coincide-1} \\
c_3^1(r)&=\overline{q}c_2^2(r)+(\overline{q}^3-\overline{q})c_2^0(r),  \label{eq:coincide-2}  \\
c_3^2(r)&=\overline{q}c_2^3(r)+(\overline{q}^5-\overline{q})c_2^1(r).  \label{eq:coincide-3}
\end{align}

\medskip

Formally put
$$\left(\begin{array}{cc} \check{t}_{13} \\ \check{t}_{23} \end{array}\right)=\left(\begin{array}{cc} t_{13} \\ t_{23} \end{array}\right)
+\frac{1}{\alpha^2-t_{12}^2}\left(\begin{array}{cc} t_{12} & \alpha \\ \alpha & t_{12} \end{array}\right)
\left(\begin{array}{cc}  t_2t_3+t_1t_{123} \\ t_1t_3+t_2t_{123} \end{array}\right).$$
Then
$$\left(\begin{array}{cc} \check{t}_{13} \\ \check{t}_{23} \end{array}\right)t_{12}
=\left(\begin{array}{cc} q^2t_{12} & \overline{q}-q^3 \\ q-\overline{q}^{3} & \overline{q}^{2}t_{12} \end{array}\right)
\left(\begin{array}{cc} \check{t}_{13} \\ \check{t}_{23} \end{array}\right).$$
The very definition of $\eta_k^n$ leads us to
\begin{align}
\left(\begin{array}{cc} \check{t}_{13} \\ \check{t}_{23} \end{array}\right)t_{12}^n
=\left(\begin{array}{cc} \overline{q}\eta_1^n(t_{12}) & -q\eta_0^n(t_{12})  \\  \overline{q}\eta_0^n(t_{12})  & -q\eta_{-1}^n(t_{12}) \end{array}\right)
\left(\begin{array}{cc} \check{t}_{13} \\ \check{t}_{23} \end{array}\right).  \label{eq:multiply}
\end{align}

\begin{figure}[h]
  \centering
  \includegraphics[width=11.5cm]{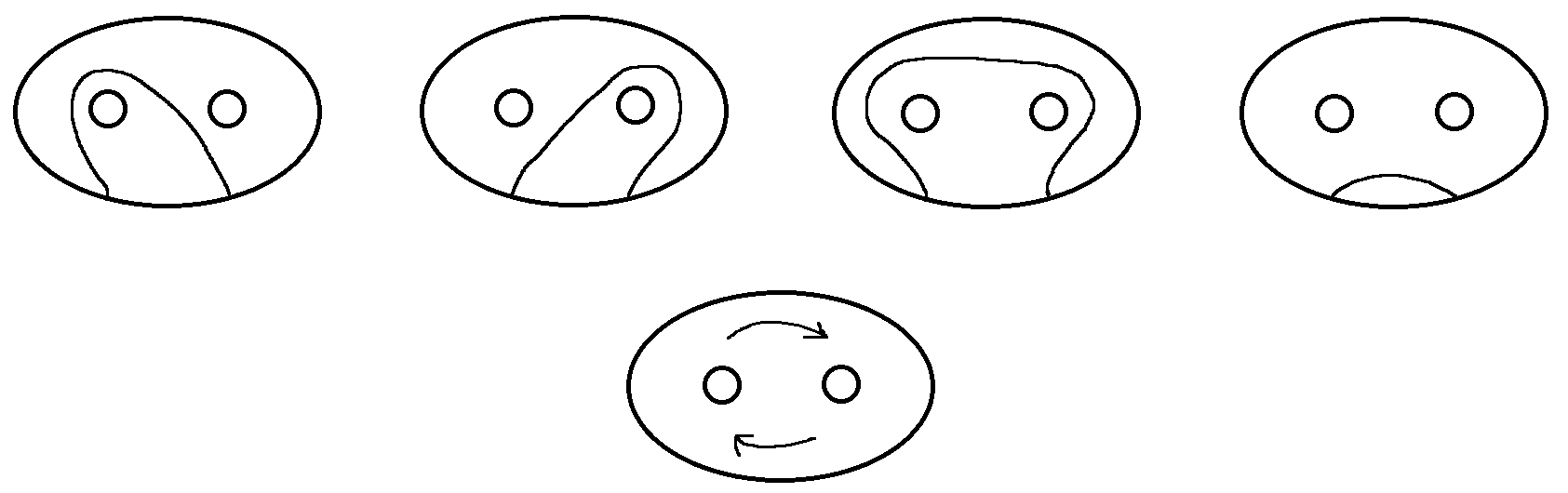}\\
  \caption{First row (from left to right): $\mathbf{x}_1$, $\mathbf{x}_2$, $\mathbf{x}_{12}$, $\mathbf{e}$.
  Second row: $\sigma$.}\label{fig:Sigma2}
\end{figure}

\begin{figure}[h]
  \centering
  \includegraphics[width=9.5cm]{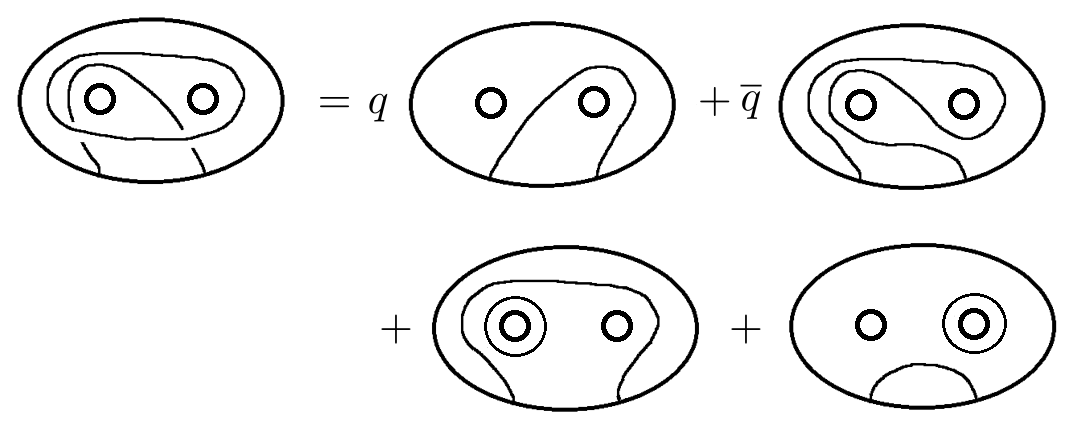}\\
  \caption{$t_{12}\mathbf{x}_1=q\mathbf{x}_2+\overline{q}\sigma(\mathbf{x}_1)+t_2\mathbf{e}+t_1\mathbf{x}_{12}$.}\label{fig:action}
\end{figure}

Fix $\mathsf{a},\mathsf{b}\in \Sigma_{0,3}\times\{\frac{1}{2}\}$.
Stacking makes $\mathcal{S}(\mathsf{a},\mathsf{b})$ into a $\mathcal{S}_2$-$\mathcal{S}_2$-bimodule.

Let $\mathbf{x}_1,\mathbf{x}_2,\mathbf{x}_{12},\mathbf{e}$ denote the elements of $\mathcal{S}(\mathsf{a},\mathsf{b})$ shown in the first row of Figure \ref{fig:Sigma2}. Let $\sigma:\Sigma_{0,3}\to\Sigma_{0,3}$ denote the diffeomorphism illustrated in the second row; we require it to fix the outer boundary of $\Sigma_{0,3}$ pointwise.

Obviously, $\sigma(\mathbf{x}_2)=\mathbf{x}_1$. From Figure \ref{fig:action} we see
\begin{align}
\sigma(\mathbf{x}_1)=qt_{12}\mathbf{x}_1-q^2\mathbf{x}_2-q(t_1\mathbf{x}_{12}+t_2\mathbf{e}).  \label{eq:action}
\end{align}
Let $r=t_{12}$, set $t_1=t_2=t$, and formally set
\begin{align*}
\check{\mathbf{x}}_i=\mathbf{x}_i-\frac{t(\mathbf{x}_{12}+\mathbf{e})}{r-\alpha}, \qquad i=1,2.
\end{align*}
Then
$$\sigma^k\left(\begin{array}{cc} \check{\mathbf{x}}_1 \\ \check{\mathbf{x}}_2 \end{array}\right)
=\left(\begin{array}{cc} qr & -q^2 \\ 1 & 0 \end{array}\right)^k
\left(\begin{array}{cc} \check{\mathbf{x}}_1 \\ \check{\mathbf{x}}_2 \end{array}\right)
=q^k\left(\begin{array}{cc} \gamma_{k+1}(r) & -q\gamma_k(r) \\ \overline{q}\gamma_k(r) & -\gamma_{k-1}(r) \end{array}\right)
\left(\begin{array}{cc} \check{\mathbf{x}}_1 \\ \check{\mathbf{x}}_2 \end{array}\right).$$

Moreover, similarly as (\ref{eq:multiply}), we can deduce
$$\left(\begin{array}{cc} \check{\mathbf{x}}_1 \\ \check{\mathbf{x}}_2 \end{array}\right)r^n=
\left(\begin{array}{cc} \overline{q}\eta_1^n(r) & -q\eta_0^n(r) \\ \overline{q}\eta_0^n(r) & -q\eta_{-1}^n(r) \end{array}\right)
\left(\begin{array}{cc} \check{\mathbf{x}}_1 \\ \check{\mathbf{x}}_2 \end{array}\right).$$
Hence
$$\sigma^k\left(\begin{array}{cc} \check{\mathbf{x}}_1 \\ \check{\mathbf{x}}_2 \end{array}\right)r^n=
\left(\begin{array}{cc} \overline{q}\eta_{k+1}^n(r) & -q\eta_k^n(r) \\ \overline{q}\eta_k^n(r) & -q\eta_{k-1}^n(r) \end{array}\right)
\left(\begin{array}{cc} \check{\mathbf{x}}_1 \\ \check{\mathbf{x}}_2 \end{array}\right).$$
Substituting $\check{\mathbf{x}}_i$ back, we obtain
\begin{align}
\sigma^k(\mathbf{x}_1)r^n&=\overline{q}\eta_{k+1}^n(r)\mathbf{x}_1-q\eta_k^n(r)\mathbf{x}_2-c_k^n(r)t(\mathbf{x}_{12}+\mathbf{e}),  \label{eq:expression-1}  \\
\sigma^k(\mathbf{x}_2)r^n&=\overline{q}\eta_k^n(r)\mathbf{x}_1-q\eta_{k-1}^n(r)\mathbf{x}_2-c_{k-1}^n(r)t(\mathbf{x}_{12}+\mathbf{e}).  \label{eq:expression-2}
\end{align}

\begin{figure}[H]
  \centering
  \includegraphics[width=12.4cm]{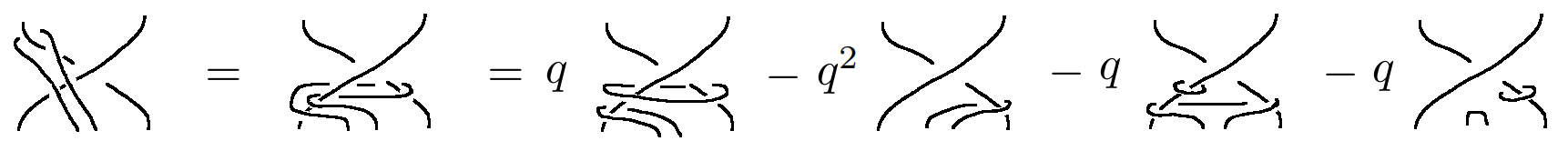}\\
  \caption{The second equality is due to (\ref{eq:action}).}\label{fig:technique-1}
\end{figure}

\begin{figure}[H]
  \centering
  \includegraphics[width=13.2cm]{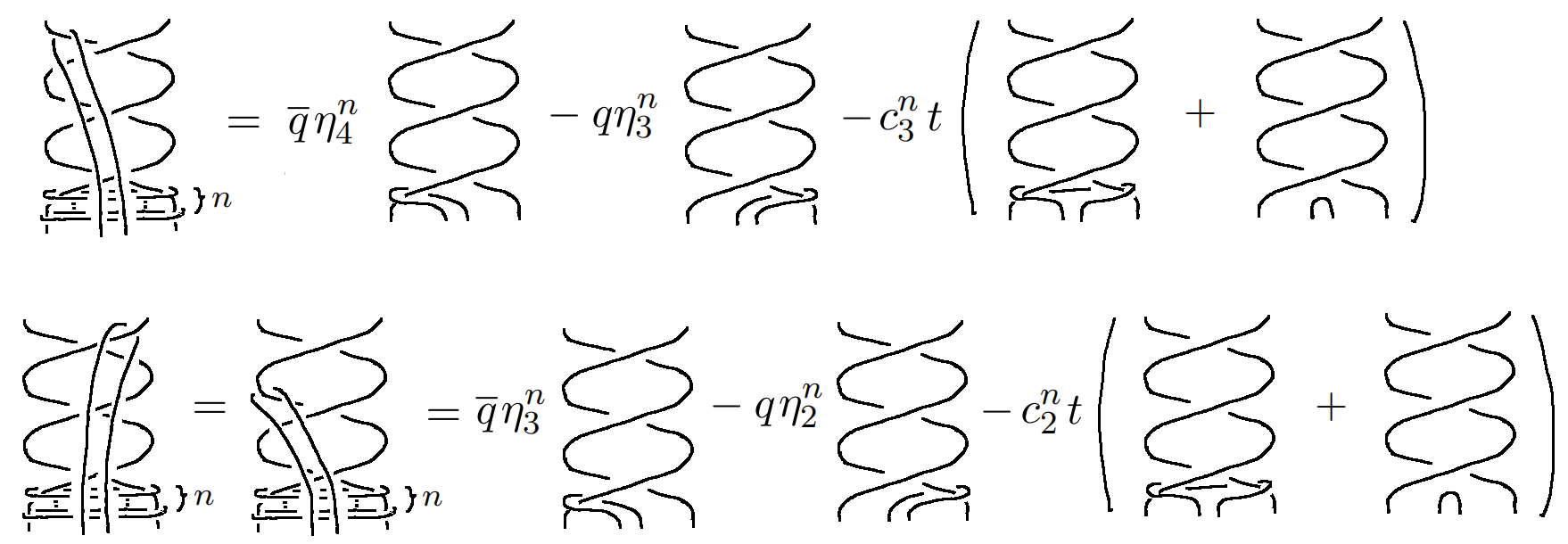}\\
  \caption{These are obtained by applying the formula in Figure \ref{fig:technique-1}. Here $\eta_k^n=\eta_k^n(r_j)$ and $c_k^n=c_k^n(r_j)$, for $j\in\{1,2,3,4\}$.}\label{fig:technique-2}
\end{figure}

Based on (\ref{eq:action}), (\ref{eq:expression-1}), (\ref{eq:expression-2}), we clarify the steps for computing ${\rm L}_1$-relators, as shown in Figure \ref{fig:technique-1}, \ref{fig:technique-2}.

\section*{Appendix C: Explicit computations for ${\rm L}_1$-relators}

\begin{figure}[H]
  \centering
  \includegraphics[width=12.7cm]{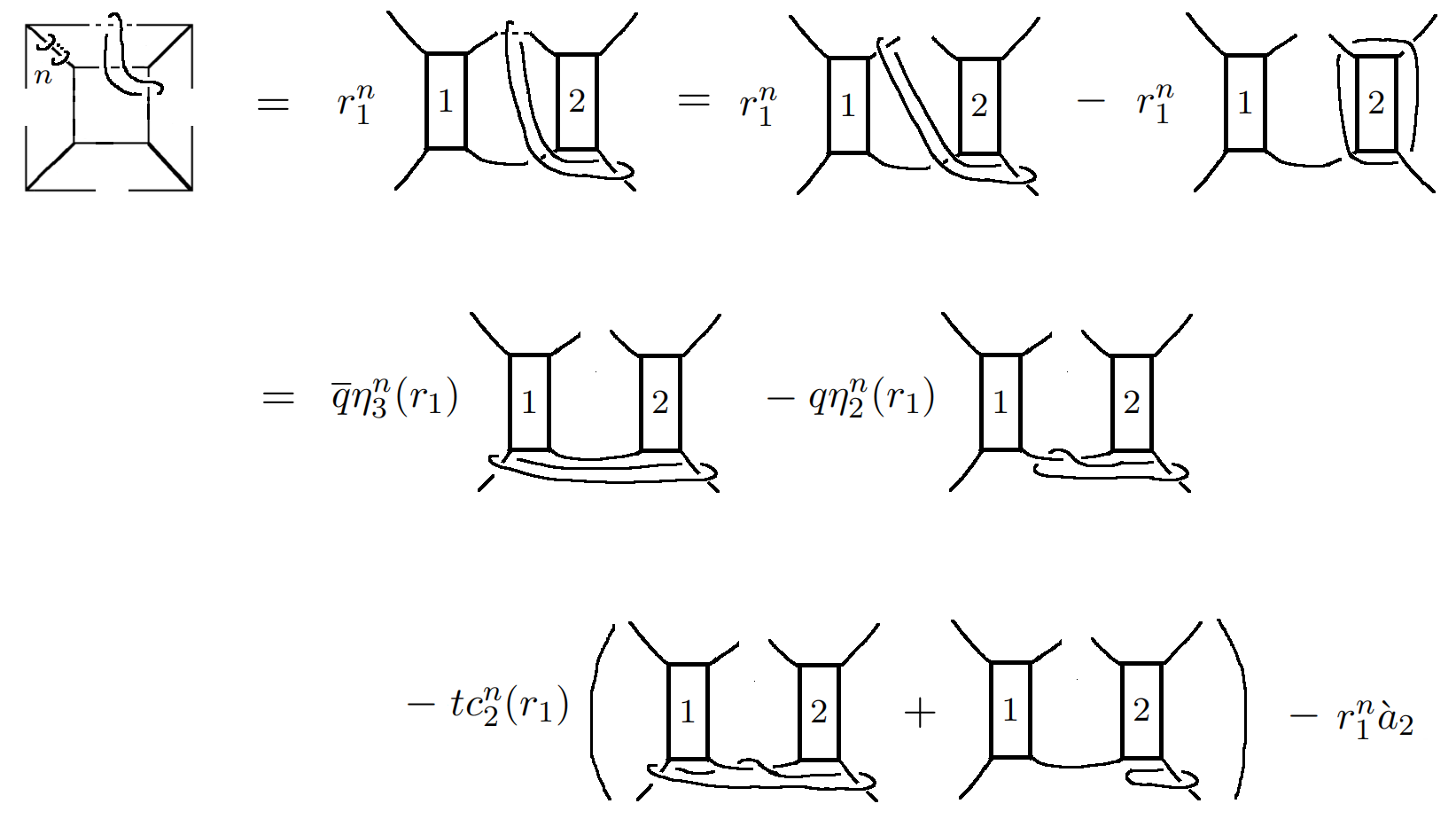}\\
  \caption{}\label{fig:relator-1}
\end{figure}

\begin{figure}[H]
  \centering
  \includegraphics[width=12.7cm]{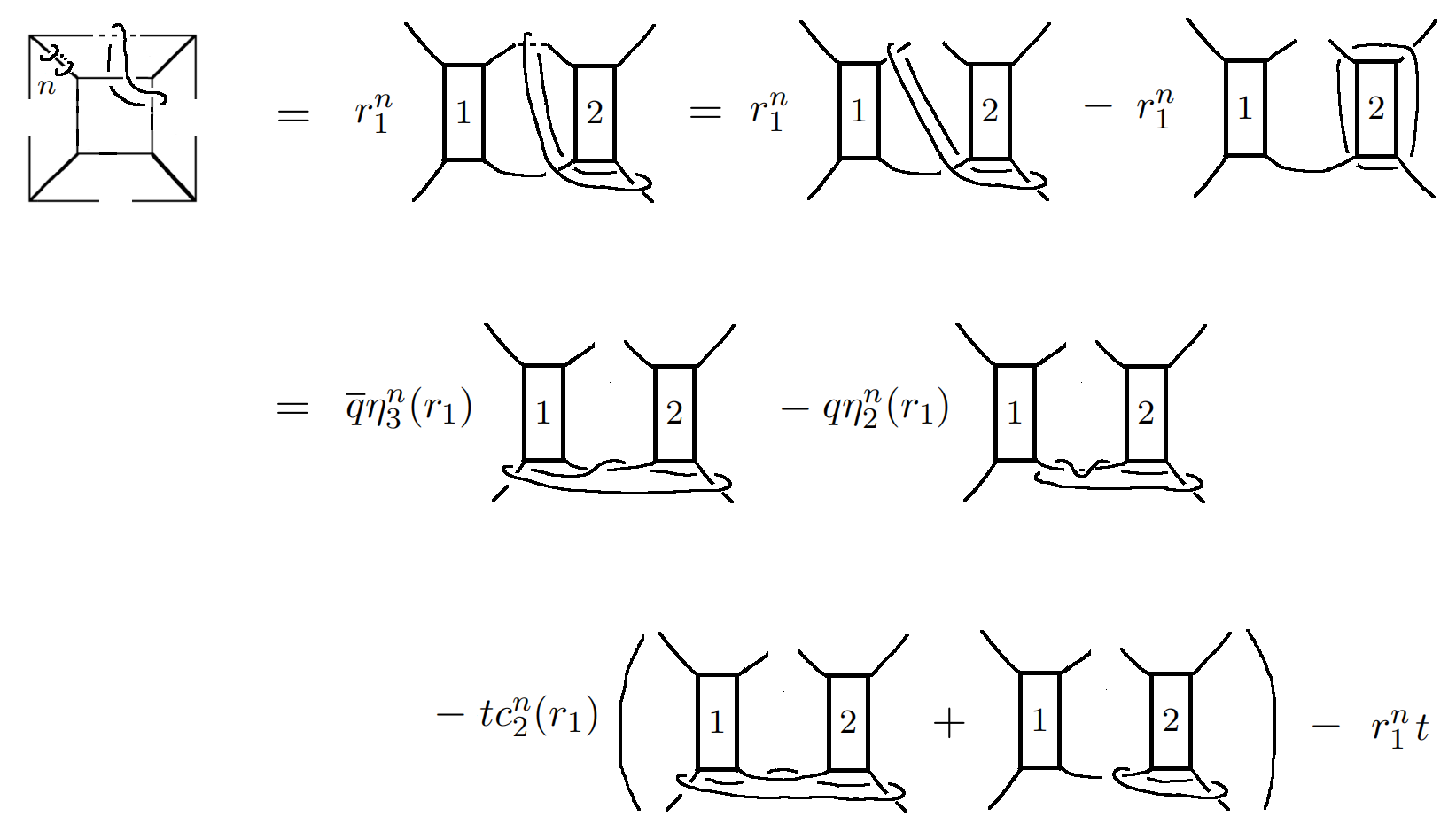}\\
  \caption{}\label{fig:relator-2}
\end{figure}

Clearly,
\begin{align}
{\rm Sl}^1_1(r_1^mr_2^n)=(\lambda_2^m(r_1)+\lambda_3^m(r_1))r_2^n-r_1^m(\lambda_2^n(r_2)+\lambda_3^n(r_2)).  \label{eq:L1-1}
\end{align}
As illustrated in Figure \ref{fig:relator-1}, Figure \ref{fig:relator-2},
\begin{align}
{\rm Sl}^1_2(r_1^n)&\approx r_1^n\grave{a}_2-\big(\overline{q}\eta_3^n(r_1)p_{14}-q\eta_2^n(r_1)p_{24}-tc_2^n(r_1)(p_{124}+t)\big)  \nonumber  \\
&\equiv c_2^n(r_1)tp_{124}-\overline{q}\eta_3^n(r_1)p_{14}.  \label{eq:L1-2}  \\
{\rm Sl}^1_{\overline{1}2}(r_1^n)&\approx r_1^nt-\big(\overline{q}\eta_3^n(r_1)p_{134}-q\eta_2^n(r_1)t-tc_2^n(r_1)(p_{14}+r_2)\big) \nonumber\\
&\equiv q\eta_3^n(r_1)p_{124}+tc_3^n(r_1)p_{14}.  \label{eq:L1-3}
\end{align}
By SP,
\begin{align}
{\rm Sl}^1_4(r_2^n)&\approx c_2^n(r_2)tp'_{134}-\overline{q}\eta_3^n(r_2)p_{14}
\approx c_2^n(r_2)tp_{124}-\overline{q}\eta_3^n(r_2)p_{14},  \label{eq:L1-4}   \\
{\rm Sl}^1_{\overline{4}1}(r_2^n)&\approx \eta_3^n(r_2)p'_{134}+\overline{q}tc_3^n(r_2)p_{14}
\approx \eta_3^n(r_2)p_{124}+\overline{q}tc_3^n(r_2)p_{14}.    \label{eq:L1-5}
\end{align}

\begin{figure}[H]
  \centering
  \includegraphics[width=12cm]{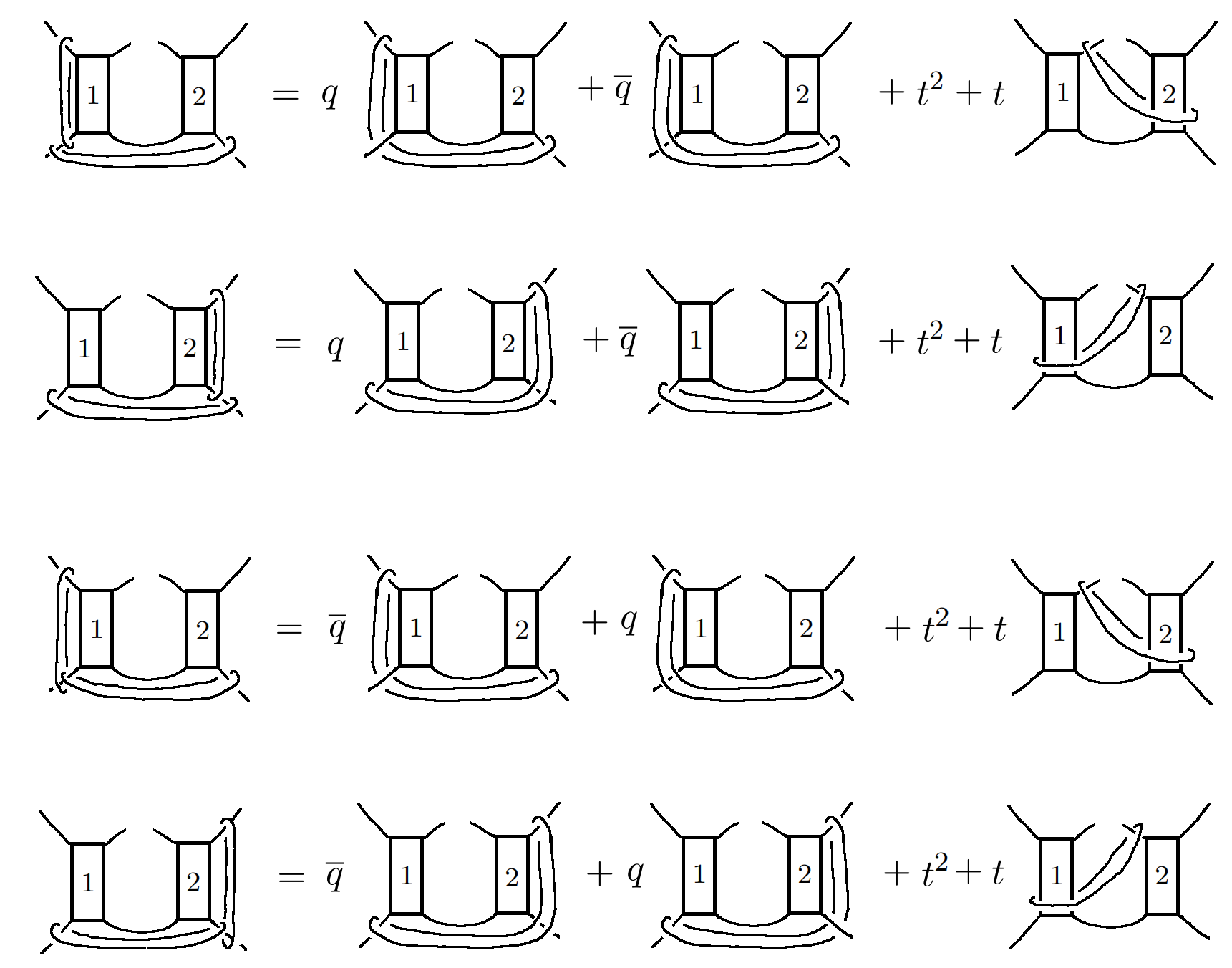}\\
  \caption{}\label{fig:b=a-2}
\end{figure}

From Figure \ref{fig:b=a-2} we see that ${\rm Sl}^1_1(r_1^mr_2^nx_1)$ can be replaced by $\acute{b}_{12}-\grave{a}_{12}$ and $\grave{b}_{12}-\acute{a}_{12}$.

\begin{figure}[H]
  \centering
  \includegraphics[width=12cm]{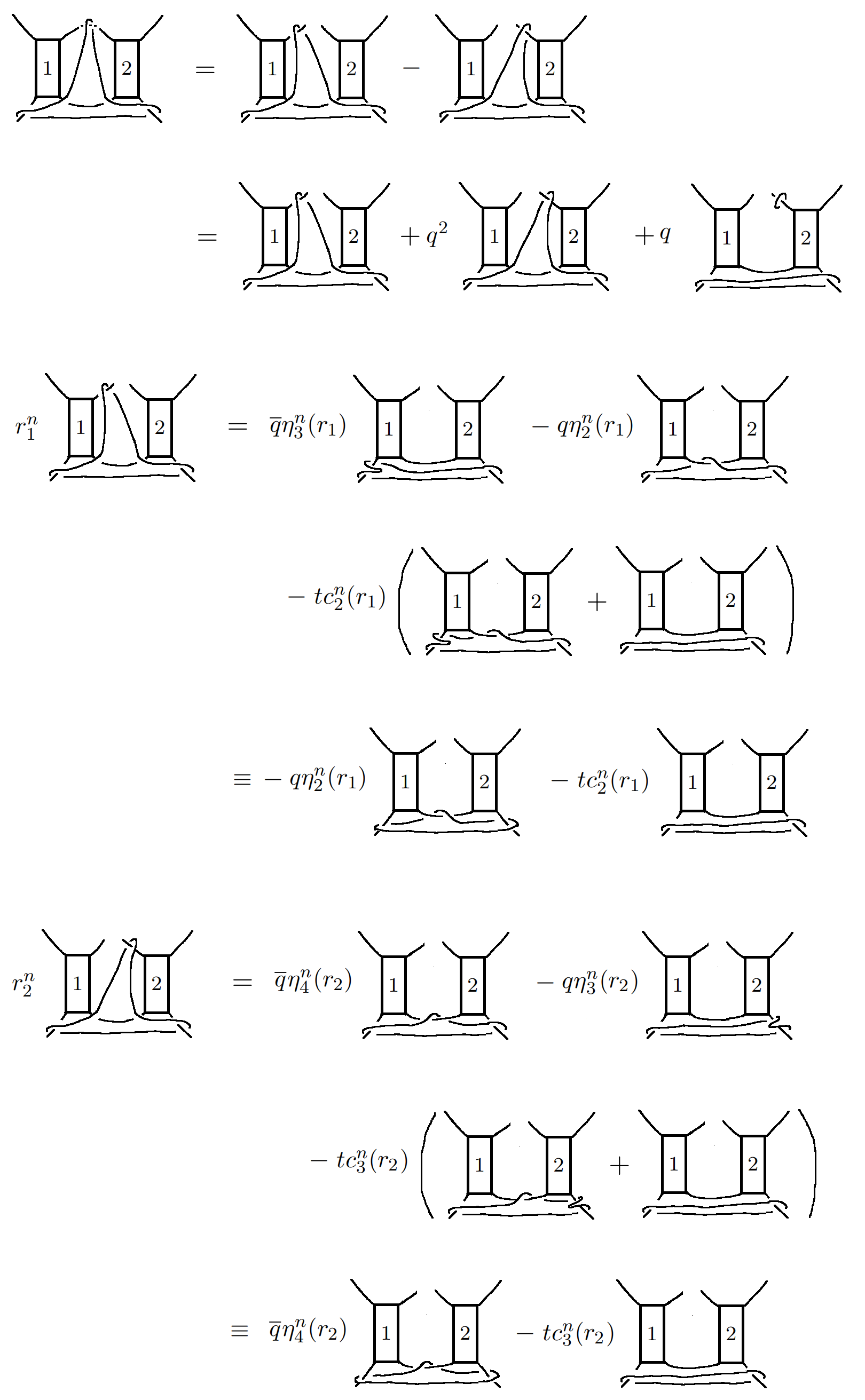}\\
  \caption{}\label{fig:relator-5}
\end{figure}

From Figure \ref{fig:relator-5} we see
\begin{align*}
&{\rm Sl}^1_{\overline{4}2}(r_1^mr_2^n)\approx -\eta_2^m(r_1)r_2^np_{124}+qr_1^m\eta_4^n(r_2)p_{134}  \\
&\hspace{23mm} +(r_1^mr_2^n-\overline{q}c_2^m(r_1)r_2^n-qr_1^mc_3^n(r_2))tp_{14}  \\
\equiv\ &\big(r_1^mr_2^n-\overline{q}c_2^m(r_1)r_2^n-qr_1^mc_4^n(r_2)\big)tp_{14}-(\eta_2^m(r_1)r_2^n+q^2r_1^m\eta_4^n(r_2))p_{124}.
\end{align*}
Modulo ${\rm Sl}^1_{\overline{4}1}(r_2^{n+1})$,
$$q^2\eta_4^n(r_2))p_{124}=q\eta_3^{n+1}(r_2)p_{124}-\eta_2^n(r_2)p_{124}\equiv -c_3^{n+1}(r_2)tp_{14}-\eta_2^n(r_2)p_{124},$$
so ${\rm Sl}^1_{\overline{4}2}(r_1^mr_2^n)$ can be replaced by
\begin{align}
\overline{q}(r_1^mc_2^n(r_2)-c_2^m(r_1)r_2^n)tp_{14}+(r_1^m\eta_2^n(r_2)-\eta_2^m(r_1)r_2^n)p_{124}.
\end{align}

\begin{rmk}\label{rmk:L1}
\rm Observe that
\begin{align*}
r_1^n\grave{a}_{12}-q\lambda_3^n(r_1)p_{14}&\equiv -q\eta_4^n(r_1)p_{14}+q^2(c_2^n(r_1)+\eta_3^n(r_1))tp_{124}-q\lambda_3^n(r_1)p_{14}  \\
&\equiv -q^2{\rm Sl}^1_2(r_1^n)-q^2t{\rm Sl}^1_{\overline{1}2}(r_1^n),   \\
r_2^n\acute{a}_{12}-\overline{q}\lambda_2^n(r_2)p_{14}&\equiv \overline{q}\eta_2^n(r_2)p_{14}-\overline{q}^2c_2^n(r_2)t(q^2p_{124}+qtp_{14})-\overline{q}\lambda_2^n(r_2)p_{14}  \\
&\equiv {\rm Sl}^1_4(r_2^n).
\end{align*}
Hence, ${\rm Sl}^1_2(r_1^n)=0$, ${\rm Sl}^1_4(r_2^n)=0$ can respectively be replaced by
\begin{align}
r_1^n\grave{a}_{12}\equiv q\lambda_3^n(r_1)x_1, \qquad  r_2^n\acute{a}_{12}\equiv \overline{q}\lambda_2^n(r_2)x_1.   \label{eq:grave-acute}
\end{align}
\end{rmk}

\begin{lem}
For $k\ge 1$,
$${\rm Sl}^1_2(r_1^k)\approx (r_1^{k+2}+\psi)p_{14}+\zeta{\rm Sl}^1_2(1)
+{\rm Sl}^1_{\overline{1}2}(\varphi)$$
for some $\zeta\in R[t]$, $\psi,\varphi\in R[t][r_1]$ with $\deg_{r_1}\psi<k+2$ and $\deg_{r_1}\varphi<k$.
\end{lem}

\begin{proof}
The case $k=1,2,3$ has been established. Suppose $k\ge 4$ and the assertion holds for smaller $k$.
Since the leading term of $\lambda_3^{k-3}(r)$ is $-q^{2k-2}r^k$, we have
$r_1^k=-\overline{q}^{2k-2}r_1^{k-3}s_1+\varrho$ for some $\varrho\in R[t][r_1]$ with $\deg_{r_1}\varrho<k$.
By the inductive hypothesis, it suffices to deal with ${\rm Sl}^1_2(r_1^{k-3}s_1)\equiv {\rm Sl}^1_2(r_1^{k-3})s_1$.

To this end, note that ${\rm Sl}^1_{\overline{1}2}(\varphi)s_1\equiv{\rm Sl}^1_{\overline{1}2}(\varphi s_1)$, and for each $j<k$,
\begin{align*}
r_1^jp_{14}s_1&\equiv q^2r_1^js_1p_{14}+(\overline{q}-q^3)r_1^j\grave{b}_{12}  \qquad \qquad  (\text{by\ Figure\ \ref{fig:b=a-2}})  \\
&\equiv q^2r_1^js_1p_{14}+\big(\overline{q}-q^3)(\overline{q}\lambda_2^j(r_1)p_{14}+{\rm Sl}^1_2(r_1^j)\big).
\end{align*}

Furthermore, we can write
$${\rm Sl}^1_2(1)s_1\approx{\rm Sl}^1_2(s_1)\approx (r_1^5+\psi')p_{14}+\zeta'{\rm Sl}^1_2(1)+{\rm Sl}^1_{\overline{1}2}(\varphi'),$$
for some $\zeta'\in R[t]$ and $\varphi',\psi'\in R[t][r_1]$ with $\deg_{r_1}\psi'<5$ and $\deg_{r_1}\varphi'<3$.

\end{proof}

Consequently, when ${\rm Sl}^1_2(1)$ and ${\rm Sl}^1_{\overline{1}2}(r_1^n)$'s are accepted, the ${\rm Sl}^1_2(r_1^k)$'s for $k\ge 1$ exactly allow us to reduce $r_1^np_{14}$ for $n\ge 3$.

\begin{figure}[H]
  \centering
  \includegraphics[width=12cm]{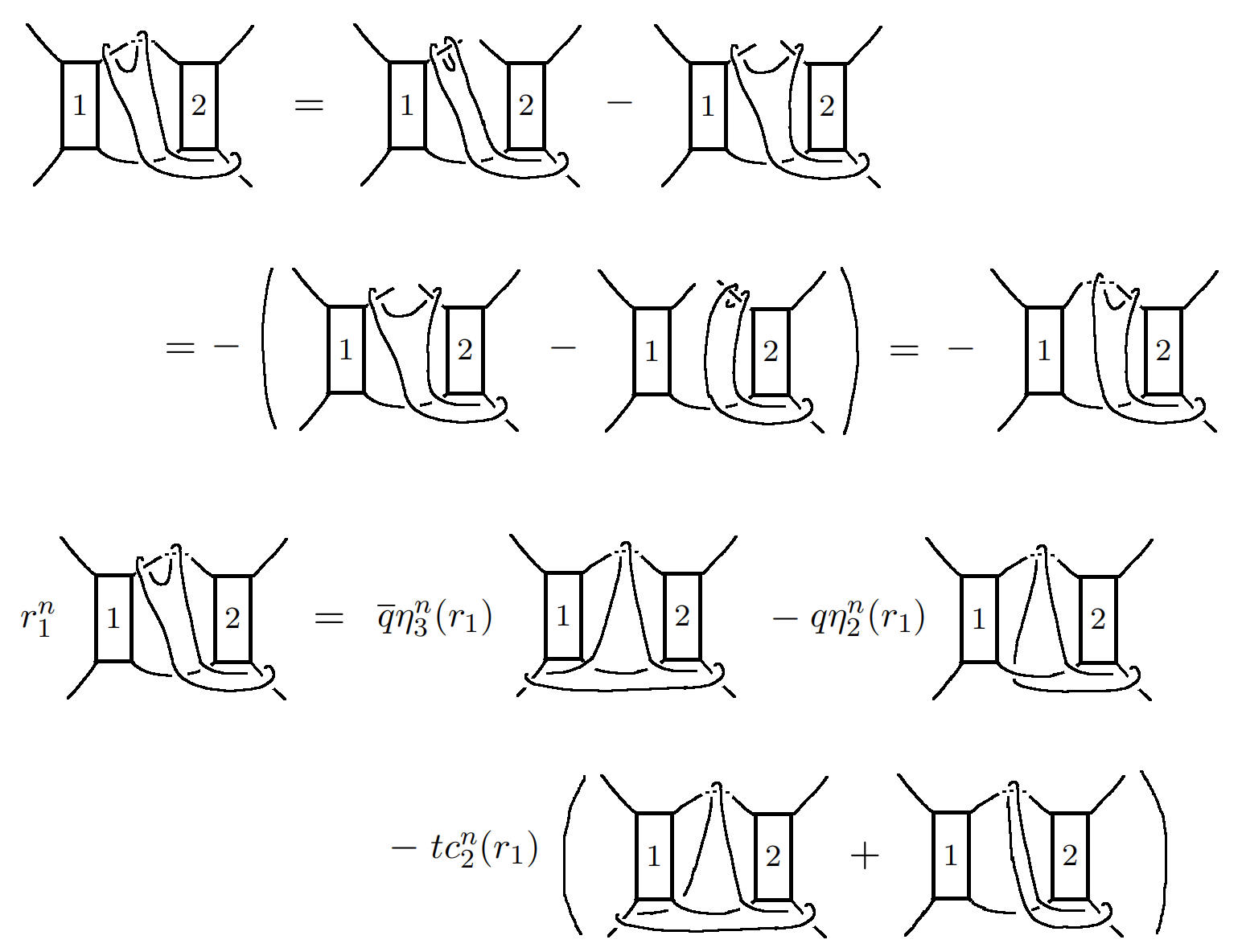}\\
  \caption{}\label{fig:redundance-1}
\end{figure}

\begin{figure}[H]
  \centering
  \includegraphics[width=12cm]{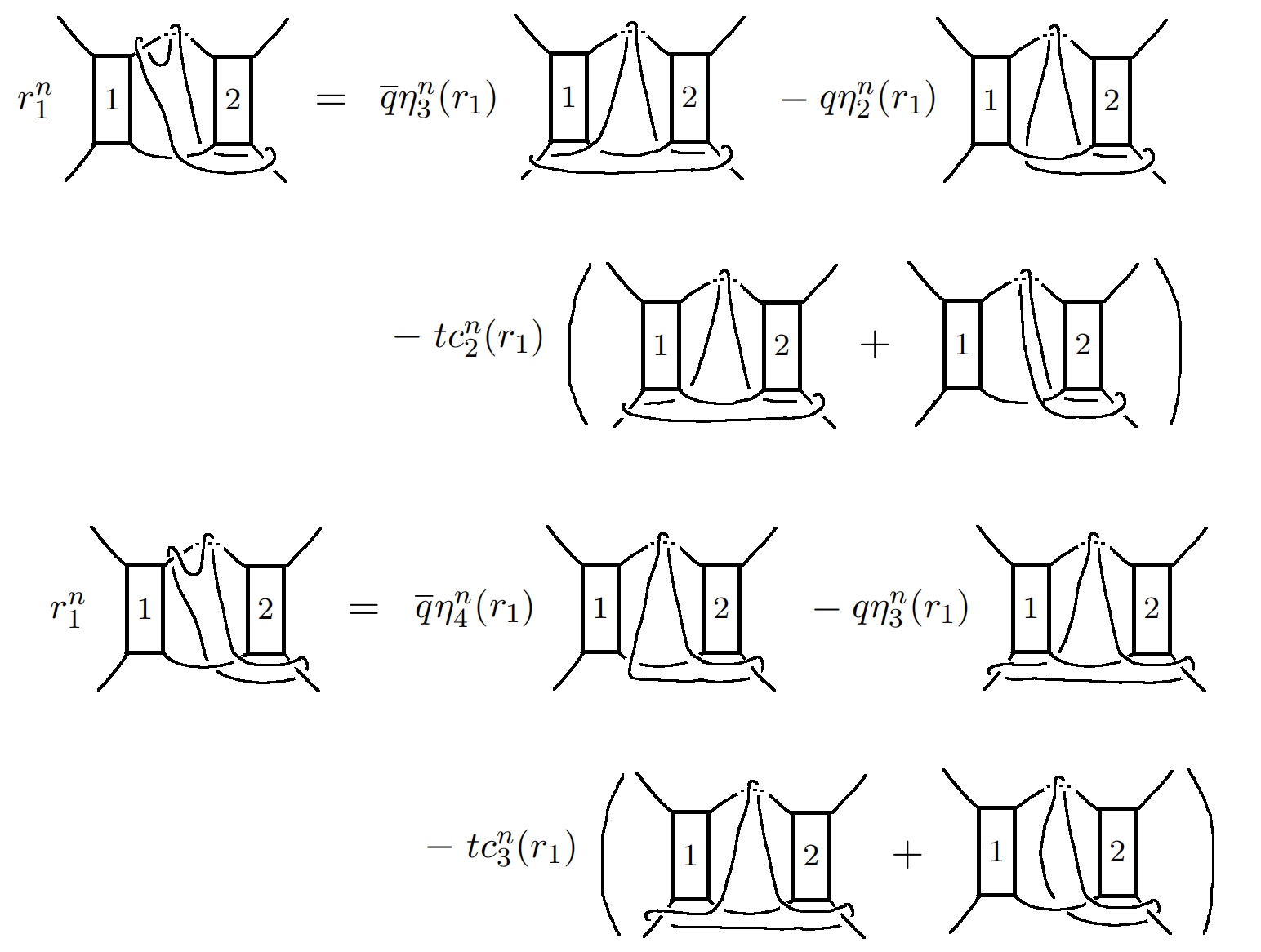}\\
  \caption{}\label{fig:redundance-2}
\end{figure}

As shown in Figure \ref{fig:redundance-1},
$\overline{q}{\rm Sl}^1_{4\overline{2}}(\eta_3^n(r_1))-t{\rm Sl}^1_{\overline{1}4\overline{2}}(c_2^n(r_1))$ equals a $R[t]$-linear combination of relators of the form ${\rm Sl}^1_{2}(u)$, ${\rm Sl}^1_{\overline{1}2}(u)$ and ${\rm Sl}^1_1(u)$. So ${\rm Sl}^1_{4\overline{2}}(r_1^2r_2^n)$
is redundant.

Shown in Figure \ref{fig:redundance-2} are another two relations among relators, showing that $r_1^2r_2^n(\grave{a}_{12}-\acute{b}_{12})$, $r_1^2r_2^n(\grave{b}_{12}-\acute{a}_{12})$ are redundant.

By SP, ${\rm Sl}^1_{4\overline{2}}(r_1^mr_2^2)$, $r_1^mr_2^2(\grave{a}_{12}-\acute{b}_{12})$, $r_1^mr_2^2(\grave{b}_{12}-\acute{a}_{12})$ 
are also redundant. 

Thus, it suffices to consider 
${\rm Sl}^1_{4\overline{2}}(r_1^mr_2^2)$, $r_1^mr_2^2(\grave{a}_{12}-\acute{b}_{12})$, $r_1^mr_2^2(\grave{b}_{12}-\acute{a}_{12})$ 
for $0\le m,n\le 1$.

\begin{proof}[Proof of Lemma \ref{lem:L1}]
Assuming ${\rm Sl}^1_2(r_1^n)=0$, ${\rm Sl}^1_{\overline{1}2}(r_1^n)=0$, ${\rm Sl}^1_4(r_2^n)=0$, ${\rm Sl}^1_{\overline{4}1}(r_2^n)=0$, we have
\begin{align*}
0={\rm Sl}^1_{\overline{4}2}(1)&\approx (r_1-r_2)p_{124}+\overline{q}t(r_1-r_2)p_{14},  \\
0={\rm Sl}^1_{\overline{4}2}(r_1)&\approx (\eta_2^1(r_1)-q^2r_1r_2)p_{124}+\overline{q}t(c_2^1(r_1)-r_1c_2^0(r_2))p_{14}  \\
&\equiv (1-q^2r_1r_2)p_{124}+t(r_1-qr_1r_2)p_{14},  \\
0={\rm Sl}^1_{\overline{4}2}(r_2)&\approx (1-q^2r_1r_2)p_{124}+t(r_2-qr_1r_2)p_{14},
\end{align*}
which can be replaced by
$$(r_1-r_2)p_{124}\equiv 0, \quad (r_1r_2-\overline{q}^2)p_{124}+t(\overline{q}r_1r_2+\overline{q}^2r_2)p_{14}\equiv 0, \quad 
t(r_1-r_2)p_{14}\equiv 0.$$
As can be verified easily, ${\rm Sl}^1_{\overline{4}2}(r_1r_2)=0$ holds automatically.

Next, consider
\begin{align*}
0=r_1^mr_2^n(\grave{a}_{12}-\acute{b}_{12})\equiv (\lambda_3^m(r_1)r_2^n-r_1^m\lambda_3^n(r_2))p_{14},  \\
0=r_1^mr_2^n(\grave{b}_{12}-\acute{a}_{12})\equiv (\lambda_2^m(r_1)r_2^n-r_1^m\lambda_2^n(r_2))p_{14}.
\end{align*}
\begin{enumerate}
  \item The case $m=n=0$ is equivalent to $r_1^i\sim r_2^i$ for $i=1,2$.
  \item The case $m=1,n=0$ gives rise to $r_1r_2\sim\lambda_3^1(r_1)\sim r_1^2+\overline{q}^2-1\sim r_2^2+\overline{q}^2-1$, and
        $1\sim\lambda_2^1(r_1)\sim r_1\lambda_2^0(r_2)$ which is equivalent to
        $$r_1(r_2^2-1)\sim\overline{q}(1-t^2)r_1r_2-\overline{q}^2t^2r_2-\overline{q}^3.$$
  \item By SP, the case $m=0,n=1$ gives rise to $r_1r_2\sim r_2^2+\overline{q}^2-1$ and
        $$(r_1^2-1)r_2\sim\overline{q}(1-t^2)r_1r_2-\overline{q}^2t^2r_2-\overline{q}^3.$$
  \item Then the case $m=n=1$ becomes trivial, as can be verified.
\end{enumerate}

Finally,
$$\eta_3^0(r_1)\eta_3^0(r_2)p_{14}\equiv q\eta_3^0(r_1)c_2^0(r_2)tp_{124}\equiv -c_3^0(r_1)c_2^0(r_2)t^2p_{14},$$
yielding $(r_1^2-1)(r_2^2-1)\sim -\overline{q}^6t^2c_3^0(r_1)c_2^0(r_2)$.
\end{proof}

\bigskip

\noindent
Haimiao Chen (orcid: 0000-0001-8194-1264)\ \ \  \emph{chenhm@math.pku.edu.cn} \\
Department of Mathematics, Beijing Technology and Business University, \\
Liangxiang Higher Education Park, Fangshan District, Beijing, China.


\begin{thebibliography}{}


\bibitem{AF22}
J.R. Aranda and N. Ferguson,
Generating sets for the Kauffman skein module of a family of Seifert fibered spaces,
{\it New York J. Math.} 28 (2022), 44--68.



\bibitem{BP22}
R.P. Bakshi and J.H. Przytycki,
Kauffman bracket skein module of the connected sum of handlebodies: a counterexample,
{\it Manuscripta Math.} 167 (2022), 809--820.



\bibitem{Bu97}
D. Bullock,
Rings of ${\rm SL}_2(\mathbb{C})$-characters and the Kauffman bracket skein module,
{\it Comment. Math. Helv.} 72 (1997), 521--542.



\bibitem{Ca17}
A. Carrega,
Nine generators of the skein space of the 3-torus,
{\it Algebr. Geom. Topol.} 17 (2017), no. 6, 3449--3460.


\bibitem{Ch18}
H.-M. Chen,
Trace-free ${\rm SL}(2,\mathbb{C})$-representations of Montesinos links,
{\it J. Knot Theory Ramifications} 27 (2018), no. 8, 1850050, 10 pp.


\bibitem{Ch22-0}
H.-M. Chen,
The ${\rm SL}(2,\mathbb{C})$-character variety of a Montesinos knot,
{\it Stud. Sci. Math. Hung.} 59 (2022), no. 1, 75--91.


\bibitem{Ch22}
H.-M. Chen,
On skein algebras of planar surfaces,
arXiv: 2206.07856.


\bibitem{Ch24}
H.-M. Chen,
Kauffman bracket skein algebra of the $4$-punctured disk,
arXiv: 2411.15829.




\bibitem{De21}
R. Detcherry,
Infinite families of hyperbolic 3-manifolds with finite-dimensional skein modules,
{\it J. Lond. Math. Soc.} (2) 103 (2021), no. 4, 1363--1376.


\bibitem{DKS23}
R. Detcherry, E. Kalfagiannt and A.S. Sikora,
Skein modules and character varieties of Seifert manifolds,
arXiv: 2311.03258.



\bibitem{DKS25}
R. Detcherry, E. Kalfagiannt and A.S. Sikora,
Kauffman bracket skein modules of small 3-manifolds,
{\it Adv. Math.} 467 (2025), 110169, 45 pp.



\bibitem{DW21}
R. Detcherry and M. Wolff,
A basis for the Kauffman skein module of the product of a surface and a circle,
{\it Algebr. Geom. Topol.} 21 (2021), no. 6, 2959--2993.


\bibitem{Gi18}
P.M. Gilmer,
On the Kauffman bracket skein module of the 3-torus,
{\it Indiana Univ. Math. J.} 67 (2018), no. 3, 993--998.


\bibitem{GM19}
P.M. Gilmer and G. Masbaum,
On the skein module of the product of a surface and a circle,
{\it Proc. Amer. Math. Soc.} 147 (2019), no. 9, 4091--4106.


\bibitem{GJS23}
S. Gunningham, D. Jordan and P. Safronov,
The finiteness conjecture for skein modules,
{\it Invent. Math.} 232 (2023), no. 1, 301--363.


\bibitem{Ki97}
R. Kirby,
{\it Problems in low-dimensional topology},
AMS/IP Stud. Adv. Math., 2.2, Geometric topology (Athens, GA, 1993), 35--473, Amer. Math. Soc., Providence, RI, 1997.


\bibitem{Ko25}
J. Korinman,
Skein modules of closed 3-manifolds define line bundles over character varieties,
arXiv: 2501.02617.


\bibitem{Le06}
Thang T.Q. L\^e,
The colored Jones polynomial and the A-polynomial of knots,
{\it Adv. Math.} 207 (2006), no. 2, 782--804.

\bibitem{LT14}
Thang T.Q. L\^e and A.T. Tran,
The Kauffman bracket skein module of two-bridge links,
{\it Proc. Amer. Math. Soc.} 142 (2014), no. 3, 1045--1056.



\bibitem{Ma10}
J. March\'e,
The skein module of torus knots,
{\it Quantum Topol.} 1 (2010), no. 4, 413--421.


\bibitem{Pr99}
J.H. Przytycki,
Fundamentals of Kauffman bracket skein modules,
{\it Kobe Math. J.} 16 (1999), 45--66.


\bibitem{PS00}
J.H. Przytycki and A.S. Sikora,
On skein algebras and ${\rm Sl}_2(\mathbb{C})$-character varieties,
{\it Topology} 39 (2000), no. 1, 115--148.


\end{thebibliography}
\end{document}